\DeclareUnicodeCharacter{2500}{}
\documentclass[11pt]{article}
\usepackage{xurl} 
\usepackage{cite} 

\usepackage[utf8]{inputenc}
\usepackage[T1]{fontenc}
\usepackage{geometry}
\usepackage{lmodern,microtype}
\usepackage{amsmath,amssymb,amsfonts,amsthm,mathtools}
\usepackage{bm,bbm}
\usepackage{graphicx}
\usepackage{booktabs}
\usepackage{enumitem}
\usepackage{algorithm}
\usepackage{algpseudocode}
\usepackage{hyperref}
\usepackage[nameinlink,capitalize]{cleveref}
\usepackage{tikz-cd}
\usepackage{subcaption} 
\usepackage[font=small,labelfont=bf]{caption}

\newcommand{\Tor}{\operatorname{Tor}}

\title{\Large Interpretability and Representability of Commutative Algebra, Algebraic Topology, and Topological Spectral Theory for Real-World Data\\[2mm]}
\author{
  Yiming Ren\thanks{%
    Department of Mathematics, Michigan State University. Email: \texttt{renyimi2@msu.edu}%
  }
  \and
  Guowei Wei\thanks{%
    Department of Mathematics, Michigan State University. Email: \texttt{weig@msu.edu}. Corresponding author.%
  }
}
\date{}

\begin{document}
\maketitle

\begin{abstract}
Recent years have witnessed a fast growth in mathematical artificial intelligence (AI). One of the most successful mathematical AI approaches is topological data analysis (TDA) via persistent homology (PH) that provides explainable AI (xAI) by extracting multiscale structural features from complex datasets. This work investigates the interpretability and representability of three foundational mathematical AI methods, PH,  persistent Laplacians (PL) derived from  spectral theory, and persistent commutative algebra (PCA) rooted in Stanley-Reisner theory. We apply these methods to a set of data, including geometric shapes, synthetic complexes, fullerene structures, and biomolecular systems to examine their geometric, topological and algebraic properties. PH captures topological invariants such as connected components, loops, and voids through persistence barcodes. PL extends PH by incorporating spectral information, quantifying topological invariants, geometric stiffness and connectivity via harmonic and non-harmonic spectra. PCA introduces algebraic invariants such as graded Betti numbers, facet persistence, and f/h-vectors, offering combinatorial, topological, geometric,  and algebraic perspectives on data over scales. Comparative analysis reveals that while PH offers computational efficiency and intuitive visualization, PL provides enhanced geometric sensitivity, and PCA delivers rich algebraic interpretability. Together, these methods form a hierarchy of mathematical representations, enabling explainable and generalizable AI for real-world data.
\end{abstract}

\tableofcontents

\section{Introduction}
Topological data analysis (TDA) has emerged as a powerful framework for extracting the hidden shapes and patterns in complex data. A key technique of TDA is persistent homology (PH), which examines how topological features such as connected components, loops, and voids form and disappear across multiple scales\cite{edelsbrunner2008surveys, zomorodian2004computing}. By constructing a multiscale filtration (growing sequence of simplicial complexes) from data, PH captures enduring homological invariants and summarizes them in persistence diagrams or barcodes \cite{ghrist2008barcodes}. This approach has become the de facto tool for quantifying shape in data and has seen wide success across scientific domains, including applications to chemistry\cite{townsend2020representation},   victories in drug design competitions \cite{nguyen2019mathematical, nguyen2020mathdl}, and the discovery viral evolution mechanisms \cite{chen2020mutations}.  
The success of TDA for molecular sciences was reviewed  \cite{wee2025review}.

Despite the success of persistent homology, it captures only the persistence of topological invariants and may overlook homotopic shape evolution in data over scales. To address these limitations, Wang et al. introduced persistent spectral theory, also called persistent  (combinatorial) Laplacians (PL) as a spectral companion to PH \cite{wang2020persistent}. Loosely speaking, a persistent Laplacian is the multiscale analogue of the combinatorial Laplace operator\cite{eckmann1944harmonische}, which generalizes the notion of graph Laplacians \cite{kirchhoff1847ueber} to higher dimensional simplicial complexes. Its harmonic eigenspaces reproduce homology (indeed, for each $q$, $\ker \mathcal{L}^{\mathrm{pers}}_{q}\cong H^{\mathrm{pers}}_{q}$), while the non‑harmonic spectrum captures additional geometric and structural information that barcodes alone do not encode. In this manner, PL fuses topological invariants with Hodge‑theoretic (spectral) features and provides a richer description of data geometry. It has been proven to outperform persistent homology, as demonstrated in tests conducted on more than 30 datasets\cite{qiu2023persistent}. Beyond the combinatorial construction, several domain‑specific persistent topological Laplacians have been proposed, including persistent sheaf Laplacians for labeled point clouds \cite{wei2025persistent1}, persistent path Laplacians for directed networks \cite{wang2023persistent}, persistent hyperdigraph Laplacians \cite{chen2023persistent}, PLs of non-branching complexes \cite{bakke2025persistent},
and persistent hypergraph Laplacians \cite{liu2021persistent}.  Computational algorithms \cite{memoli2022persistent,wang2021hermes,jones2025petls} and stability analysis \cite{liu2023algebraic} have been given to the basic persistent Laplacian. 
Persistent spectral approach has been generalized to quantum persistent homology or  persistent Dirac operators \cite{ameneyro2024quantum}. These operators have been studied across multiple topological domains with diverse applications \cite{suwayyid2024persistent,wee2023persistent}. 
A comprehensive survey of these topological spectral theory developments was given in \cite{wei2025persistent}.  Empirically, incorporating spectral features has improved modeling of geometric and structural patterns in real data, with successful PL‑based applications in areas such as protein–ligand interaction modeling \cite{meng2021persistent,chen2024multiscale}, gene‑network analysis \cite{cottrell2023plpca}, in silico deep mutational scanning \cite{chen2023topological}, and the successful prediction of COVID‑19 variant dominance about two months in advance \cite{chen2022persistent}.

In addition, algebraic topology and persistent topological spectral approaches have been extended to persistent Mayer topology \cite{shen2024persistent} and persistent interaction topology \cite{liu2025persistent}. These new approaches significantly broaden the scope of TDA. However, algebraic topology tools are designed for point cloud data or data on graphs. For a vast variety of data on manifolds and curves embedded in 3-space, one needs mathematical approaches from different fields. To this end, differential topology, i.e., evolutionary de Rham-Hodge theory \cite{chen2021evolutionary} or persistent de Rham-Hodge Laplacian, was introduced for data on smooth manifolds \cite{su2024persistent}. A similar generalization has also been discussed recently \cite{wolf2025generalized}.   
Furthermore, TDA was generalized to geometric topology, i.e., evolutionary Khovanov homology\cite{shen2024evolutionary} for curves embedded in 3-space. The reader is referred to a review \cite{su2025topological}.  

The success of algebraic topology, differential topology, geometric topology, and spectral theory approaches has motivated the exploration of persistent commutative algebra (PCA) for data. Commutative algebra is a branch of mathematics that studies commutative rings, their ideals, modules, and related algebraic 
structures \cite{miller2005combinatorial,eisenbud2013commutative}. Very recently, Suwayyid and Wei introduced persistent Stanley–Reisner theory \cite{suwayyid2025persistent} as a means to incorporate commutative algebraic invariants for multiscale algebraic data analysis (ADA). This framework utilizes the classic correspondence between simplicial complexes and commutative rings (via Stanley–Reisner ideals) to define new ADA across scales. These persistent algebraic invariants such as such as persistent graded
Betti numbers, f-vectors, h-vectors, and facet persistence barcodes are provably stable under perturbations and offer a novel perspective: they capture subtle combinatorial patterns and algebraic invariants in the data that might be invisible to homology- or Laplacian-based methods. Although the field is nascent, early applications indicate promise in molecular and biomedical settings, including protein–ligand binding \cite{feng2025caml}, genetic origins of disease \cite{wee2025commutative}, and genomic analysis \cite{suwayyid2025cakl}.

To understand the effectiveness of the aforementioned mathematical methods in data analysis, it is essential to examine why they perform so well and how their performance may vary. These methods offer several key characteristics: \\
{\it Multiscale Analysis:} These methods enable multiscale analysis, making them well-suited for addressing multiscale problems.\\
{\it Mathematical Invariants:} They leverage unique topological or algebraic invariants, which cannot be derived from alternative statistical, physical, chemical, or biological approaches.\\
{\it  Simplification and Complexity Reduction:} These methods facilitate mathematical simplification and/or complexity reduction, which is critical for handling highly complex datasets.\\
{\it High-Dimensional Representations:} They provide high-dimensional representations, such as topological/algebraic dimensions 0, 1, 2, etc., which are particularly valuable for modeling many-body interactions.\\
{\it Suitable for Machine Learning:} These methods are easily paired with machine learning (ML) or deep learning (DL) algorithms, enhancing their effectiveness for data analysis and prediction. A notable milestone in this field was the integration of TDA  with neural networks, introduced by Cang and Wei in 2017. They coined the term topological deep learning (TDL) \cite{cang2017topologynet}. \\
{\it Explainable Artificial Intelligence (xAI):} Unlike conventional deep learning methods, which often result in opaque "black-box" AI, these methods lead to explainable AI (xAI) and represent a new frontier in rational learning \cite{papamarkou2024position}. \\
The performance of these mathematical AI approaches depends on their ability to represent data effectively. The representability of algebraic topology, in particular, has been explored in the literature \cite{cang2018representability}, and earlier mathematical representations were reviewed \cite{nguyen2020review}.

In this work, we explore the interpretability and representability of PH, PL, and PCA.  The goal is to systematically evaluate how each framework encodes the shape of data via simplicial complex and to elucidate their respective strengths, differences, and limitations. 
For the sake of simplicity and direct comparison, we focus on simplicial complex, although many other topological spaces, such as cellular complex, path complex, cellular sheaf,  directed flag complex, hypergraph, etc., have been introduced in algebraic topology and topological spectral approaches\cite{su2025topological}. By applying PH, PL, and PCA to the same topological space and the same set of datasets, ranging from simple geometric shapes to high-dimensional synthetic complexes, fullerene structures, and biomolecular systems, we provide a unified view of how these methods extract topological, spectral, and algebraic features.  Each of the three frameworks targets a distinct aspect of the data information.  PH focuses on topological connectivity and cavities, capturing the birth and death of homological invariants across scales.  PL extends this perspective by incorporating geometric information, quantifying the stiffness and stability of topological structures through the spectra of multiscale Laplace operators.  PCA  analyzes topological, algebraic, and combinatorial information, encoding the evolution of ideals, syzygies, and graded Betti numbers that describe higher-order relationships among simplices. Indeed, PH has already demonstrated strong performance in molecular fingerprinting, effectively capturing topological motifs of chemical compounds and protein structures \cite{xia2014persistent}.  Building upon this foundation, PLs introduce geometric connectivity into the analysis, leading to improved modeling of networks,  landscapes, and other physical systems\cite{liu2025manifold}.  PCA, in turn, provides an additional layer of interpretability by revealing subtle combinatorial and algebraic patterns  \cite{suwayyid2025cakl}.  By jointly considering the homological, spectral, and algebraic characteristics, our study demonstrates that PH, PL, and PCA together yield a more comprehensive and robust representation of complex real-world data than any single method alone.

\section{Simplicial Complexes, Persistent Homology, and Filtrations}

Persistent homology is a multiscale extension of simplicial homology, a fundamental tool in algebraic topology that classifies simplicial complexes through their homology groups. These groups capture topological features of a simplicial complex, such as connected components, loops, and cavities. While simplicial homology provides a rigorous classification, it does not directly apply to point-cloud data, which lack an inherent topological structure unless explicitly constructed. Moreover, the outcome is highly sensitive to the choice of scale, since different scales may yield different homology groups. Persistent homology resolves these issues by introducing filtrations and persistence, thereby offering a robust multiscale characterization of the underlying data.

The simplicial complexes are built from simplices. A $k$-simplex is defined as the convex hull of $k+1$ affinely independent points $v_{0}, v_{1}, \dots, v_{k}$  
.and has dimension $k$. Examples include a point (0-simplex), a line segment (1-simplex), a triangle (2-simplex), and a tetrahedron (3-simplex). The vertices of the $k$-simplex are the points $v_{0}, v_{1}, \dots, v_{k}$, and simplices formed from any subset of these vertices are called faces.

A simplicial complex $K$ is a finite set of simplices that satisfies the following conditions: (i) every face of a simplex in $K$ is also contained in $K$, and (ii) the intersection of any two simplices in $K$ is either empty or a common face of both. The dimension of $K$ is defined as the maximum dimension of its simplices. 

Given a simplicial complex $K$, one defines the $k$-th chain group $\mathcal{C}_{k}(K)$ as the abelian group generated by all $k$-simplices of $K$, with coefficients in a field such as $\mathbb{Z}_{2}$. An oriented $k$-simplex is represented by an ordered list of its vertices $[v_{0}, v_{1}, \dots, v_{k}]$, and orientations are necessary for defining the boundary operator (though unnecessary when using $\mathbb{Z}_{2}$ coefficients, since $-1=+1$). The $k$-th boundary operator is the linear map
\[
\partial_{k} : \mathcal{C}_{k}(K) \to \mathcal{C}_{k-1}(K),
\]
defined by
\[
\partial_{k}[v_{0}, v_{1}, \dots, v_{k}] 
= \sum_{i=0}^{k} (-1)^{i}[v_{0}, \dots, \widehat{v_{i}}, \dots, v_{k}],
\]
where $\widehat{v_{i}}$ indicates omission of $v_{i}$. The boundary operators satisfy $\partial_{k}\partial_{k+1} = 0$, leading to the chain complex
\[
\cdots \xrightarrow{\partial_{k+2}} \mathcal{C}_{k+1}
\xrightarrow{\partial_{k+1}} \mathcal{C}_{k}
\xrightarrow{\partial_{k}} \mathcal{C}_{k-1}
\xrightarrow{\partial_{k-1}} \cdots.
\]

The $k$-th simplicial homology group of $K$ is then given by
\[
H_{k}(K) = \ker \partial_{k} \,/\, \mathrm{im}\,\partial_{k+1},
\]
where $\ker \partial_{k}$ is the group of $k$-cycles and $\mathrm{im}\,\partial_{k+1}$ is the group of $k$-boundaries. The homology group $H_{k}(K)$ identifies $k$-cycles that are not boundaries, i.e., the $k$-dimensional holes in the simplicial complex. Its rank, the $k$-th Betti number $\beta_{k}$, counts these features. In particular, $\beta_{0}$ gives the number of connected components, $\beta_{1}$ the number of loops or holes, and $\beta_{2}$ the number of cavities.

A filtration is a nested sequence of simplicial complexes
\[
\emptyset \subset K_{0} \subset K_{1} \subset \cdots \subset K_{m} = K,
\]
where $K$ is a maximal simplicial complex constructed from a point cloud. Common constructions include the Vietoris--Rips complex\cite{vietoris1927hoheren}, the \v{C}ech complex \cite{edelsbrunner2010computational}, and the Alpha complex\cite{edelsbrunner2011alpha}. The inclusion maps $i_{i,j} : K_{i} \hookrightarrow K_{j}$ with $i < j$ induce homomorphisms on homology groups,
\[
f_{i,j}^{k} : H_{k}(K_{i}) \to H_{k}(K_{j}),
\]
which track the evolution of $k$-dimensional features. The $k$-th persistent homology group is defined by
\[
H^{i,j}_{k} = \mathrm{im}\, f^{k}_{i,j} 
= \ker \partial^{i}_{k} \,/\, \bigl( \mathrm{im}\, \partial^{j}_{k+1} \cap \ker \partial^{i}_{k} \bigr),
\]
and its rank
\[
\beta^{i,j}_{k} = \mathrm{rank}\, H^{i,j}_{k}
\]
is the $k$-th persistent Betti number.

A topological feature is said to be born at $K_{i}$ if it appears for the first time in $K_{i}$, and it dies at $K_{j}$ if it merges into a pre-existing feature at $K_{j}$. The persistence of such a feature, denoted by $\gamma$, is its lifetime $j-i$, which is infinite if it never dies in the filtration. The collection of birth and death information leads to various representations of persistent homology, such as persistent Betti numbers, Betti curves, persistence barcodes\cite{ghrist2008barcodes}, persistence diagrams\cite{edelsbrunner2002topological}, and persistence landscapes\cite{bubenik2015statistical}. These provide a robust and multiscale characterization of the topological structure of data.

\section{Topological Spectral Theory: Persistent Laplacians}

Persistent Laplacians extend persistent homology and provide a full spectral view of filtered simplicial complexes. The harmonic spectra with zero eigenvalues recover the outputs of persistent homology. The non-harmonic spectra with positive eigenvalues record additional geometric and combinatorial information across the filtration. This extra information reflects changes in effective stiffness, reinforcement of connectivity, and reconfiguration of shape that persistent homology alone does not reveal.

Combinatorial Laplacians generalize graph Laplacians to simplicial complexes \cite{eckmann1944harmonische}. Let \(K\) be a finite simplicial complex. Write \(\mathcal{C}_{k}(K)\) for the \(k\)-chain group over a fixed field with an inner product that makes the standard simplex basis orthonormal. The boundary operator \(\partial_{k}\) maps \(\mathcal{C}_{k}(K)\) to \(\mathcal{C}_{k-1}(K)\). Its adjoint \(\partial_{k}^{*}\) maps \(\mathcal{C}_{k-1}(K)\) to \(\mathcal{C}_{k}(K)\). The \(k\)-th combinatorial Laplacian is
\[
L_{k}=\partial_{k+1}\partial_{k+1}^{*}+\partial_{k}^{*}\partial_{k},
\]
which is real, symmetric, and positive semidefinite. Hence the spectrum is contained in \([0,\infty)\). When \(K\) is a graph and \(k=0\), the operator \(L_{0}=\partial_{1}\partial_{1}^{*}\) reduces to the usual graph Laplacian. These operators yield a discrete Hodge decomposition
\[
\mathcal{C}_{k}(K)=\mathrm{im}\,\partial_{k+1}\oplus \ker L_{k}\oplus \mathrm{im}\,\partial_{k}^{*}.
\]

With the Kronecker delta inner product, basis simplices are mutually orthogonal. The matrix of \(\partial_{k}\) has size \(N_{k-1}\times N_{k}\), where \(N_{k}\) is the number of \(k\)-simplices, and the matrix of \(\partial_{k}^{*}\) is the transpose. It follows from \cite{eckmann1944harmonische} that
\[
\ker L_{k}\cong H_{k}(K),
\]
so the multiplicity of the zero eigenvalue equals the \(k\)-th Betti number. The positive eigenvalues capture geometric and combinatorial structure. For example, the first positive eigenvalue at \(k=0\) reflects how well the underlying graph is connected.

To incorporate persistence, consider a filtration \(\{K_{t}\}_{t}\) and indices \(i\le j\) with \(K_{i}\subset K_{j}\). The \(k\)-th persistent combinatorial Laplacian acts on \(\mathcal{C}_{k}(K_{i})\) and is defined by
\[
L^{i,j}_{k}
=
\partial^{i,j}_{k+1}\bigl(\partial^{i,j}_{k+1}\bigr)^{*}
+
\bigl(\partial^{i}_{k}\bigr)^{*}\partial^{i}_{k}.
\]
Here \(\partial^{i}_{k}\) is the boundary on \(K_{i}\). The map \(\partial^{i,j}_{k+1}\) is the persistent boundary obtained by restricting \(\partial^{j}_{k+1}\) on \(K_{j}\) to the largest subspace \(\mathcal{C}^{i,j}_{k+1}\subset \mathcal{C}_{k+1}(K_{j})\) whose image lies in \(\mathcal{C}_{k}(K_{i})\). The operators \((\partial^{i}_{k})^{*}\) and \((\partial^{i,j}_{k+1})^{*}\) are the adjoints with respect to the chosen inner products. When \(i=j\), the identity \(L^{i,i}_{k}=L_{k}\) holds.

The operator \(L^{i,j}_{k}\) is real, symmetric, and positive semidefinite, and it yields a persistent Hodge decomposition
\[
\mathcal{C}_{k}(K_{i})
=
\mathrm{im}\,\partial^{i,j}_{k+1}
\oplus
\ker L^{i,j}_{k}
\oplus
\mathrm{im}\,(\partial^{i}_{k})^{*}.
\]
Its kernel identifies the \(k\)-th persistent homology group
\[
\ker L^{i,j}_{k}\cong H^{i,j}_{k},
\]
so the multiplicity of the zero eigenvalue equals the \(k\)-th persistent Betti number. The harmonic spectra therefore reproduce persistent homology, while the positive spectra measure how geometry and combinatorics evolve along the filtration.

Practical implementations are available. The HERMES package \cite{wang2021hermes} builds alpha-complex filtrations and assembles a single boundary at the final Delaunay-based complex together with suitable projection matrices, which makes adjoints equal to transposes under the orthonormal simplex basis. Other algorithms for persistent Laplacians have also been proposed\cite{memoli2022persistent, jones2025petls}.

\section{Persistent Commutative Algebra}
\label{sec:PSRT}
Persistent Commutative Algebra (PCA), often referred to as Persistent Stanley–Reisner Theory, is an emerging framework that extends topological data analysis by incorporating algebraic and combinatorial structures. The fundamental idea is to represent a point cloud by a sequence of simplicial complexes—combinatorial objects constructed from vertices, edges, triangles, and higher-dimensional simplices—and then to study how these complexes evolve under a filtration. By embedding this combinatorial information into the language of commutative algebra, PCA enables the systematic tracking of structural features across scales, thereby providing an enriched algebraic perspective on the geometry and topology of data. Four key invariants are tracked through the filtration: persistent $f$–vectors, persistent $h$–vectors, graded Betti numbers, and facet ideals. Together, these invariants provide a compact algebraic–combinatorial summary of the evolving shape of data
\subsection{Persistent Stanley--Reisner Theory}

Let $\Delta$ be a finite simplicial complex on the vertex set 
$V=\{x_1,\ldots,x_n\}$, and let 
$f:\Delta\to\mathbb{R}$ be a filtration function assigning to each face 
$\sigma\in\Delta$ a real number $f(\sigma)$ interpreted as the birth scale of $\sigma$.
The function $f$ induces an increasing family of subcomplexes
\[
\Delta^{f} \;=\; \{\Delta^{\varepsilon}\}_{\varepsilon\in\mathbb{R}},
\qquad
\Delta^{\varepsilon} \;:=\; \{\sigma\in\Delta \mid f(\sigma)\le \varepsilon\},
\]
which satisfies the monotonicity condition 
$\Delta^{\varepsilon_1}\subseteq \Delta^{\varepsilon_2}$ whenever 
$\varepsilon_1\le \varepsilon_2$. 

Fix a field $k$, and consider the standard graded polynomial ring
\[
S \;:=\; k[x_1,\ldots,x_n],
\]
graded by total degree.  For each scale $\varepsilon$, the subcomplex $\Delta^{\varepsilon}$
determines its Stanley--Reisner ideal
\begin{equation}
I^{\varepsilon} \;:=\; I(\Delta^{\varepsilon})
\;=\;
\bigl(
x_{i_1}\cdots x_{i_r} \ \big|\ \{x_{i_1},\ldots,x_{i_r}\}\notin \Delta^{\varepsilon}
\bigr)\ \subseteq\ S,
\label{eq:SR_ideal_eps_long}
\end{equation}
and the associated Stanley--Reisner ring
\[
k[\Delta^{\varepsilon}] \;:=\; S/I^{\varepsilon}.
\]
Because $\Delta^{\varepsilon}$ enlarges with~$\varepsilon$, the nonface set shrinks, and hence the
ideals form a descending chain
\[
I^{\varepsilon_0}\ \supseteq\ I^{\varepsilon_1}\ \supseteq\ I^{\varepsilon_2}\ \supseteq\ \cdots,
\qquad
\varepsilon_0<\varepsilon_1<\varepsilon_2<\cdots.
\]

Each $I^{\varepsilon}$ admits a canonical decomposition as an intersection of
prime monomial ideals indexed by the facets of $\Delta^{\varepsilon}$:
\[
I^{\varepsilon} \;=\; \bigcap_{\sigma\in \mathcal{F}(\Delta^{\varepsilon})}
P_{\sigma},
\qquad
P_{\sigma} \;:=\; (x_i \mid x_i\notin \sigma)\ \subseteq\ S,
\]
where $\mathcal{F}(\Delta^{\varepsilon})$ denotes the set of maximal faces of $\Delta^{\varepsilon}$.
The $P_{\sigma}$ are the facet ideals at scale~$\varepsilon$.  We collect them in
$
\mathcal{P}^{\varepsilon}
\;:=\;
\bigl\{P_{\sigma}\ \big|\ \sigma\in\mathcal{F}(\Delta^{\varepsilon})\bigr\}.
$
To resolve the structure by dimension, we stratify
\[
\mathcal{P}^{\varepsilon}_i
\;:=\;
\bigl\{P_{\sigma}\in \mathcal{P}^{\varepsilon}\ \big|\ \dim(\sigma)=i\bigr\},
\qquad
\mathcal{P}^{\varepsilon}
\;=\;
\bigcup_{i=0}^{\dim(\Delta^{\varepsilon})}\mathcal{P}^{\varepsilon}_i.
\]
In analogy with PH, we call a facet ideal $P_{\sigma}\in \mathcal{P}^{\varepsilon}_i$
persistent to a later scale $\varepsilon'>\varepsilon$ if
$P_{\sigma}\in \mathcal{P}^{\varepsilon'}_i$ as well.
The set of $i$–dimensional facet ideals present at both scales is
\[
\mathcal{P}^{\varepsilon,\varepsilon'}_i
\;:=\;
\mathcal{P}^{\varepsilon}_i \cap \mathcal{P}^{\varepsilon'}_i,
\]
and we define the facet persistence Betti number
\[
\beta_{i}^{\varepsilon,\varepsilon'}\;:=\;\bigl|\mathcal{P}^{\varepsilon,\varepsilon'}_i\bigr|.
\]
The family $\{\beta_{i}^{\varepsilon,\varepsilon'}\}_{i,\varepsilon,\varepsilon'}$
provides a graded combinatorial signature of how the prime (facet) decomposition of the
Stanley--Reisner ideals evolve along the filtration. It is the algebraic analogue of a
barcode for facets.

\subsection{Persistent Graded Betti Numbers}

For each $\varepsilon$, the ring $k[\Delta^{\varepsilon}]=S/I^{\varepsilon}$ is a finitely generated
graded $S$–module and admits a minimal graded free resolution
\begin{equation}
\cdots \longrightarrow
\bigoplus_{j} S(-j)^{\beta_{i,j}(k[\Delta^{\varepsilon}])}
\longrightarrow \cdots \longrightarrow
\bigoplus_{j} S(-j)^{\beta_{0,j}(k[\Delta^{\varepsilon}])}
\longrightarrow k[\Delta^{\varepsilon}] \longrightarrow 0,
\label{eq:min_free_res_long}
\end{equation}
where $S(-j)$ denotes a degree-$j$ shift.  The integers
\[
\beta_{i,j}(k[\Delta^{\varepsilon}])
\;:=\;
\dim_{k}\,\Tor^S_{i}\bigl(k[\Delta^{\varepsilon}],k\bigr)_{j}
\]
are the graded Betti numbers of $k[\Delta^{\varepsilon}]$: they count minimal syzygies in
homological degree $i$ and internal degree $j$.

For any $W\subseteq V$, let
\[
\Delta^{\varepsilon}_W \;:=\; \bigl\{\sigma\in\Delta^{\varepsilon}\ \big|\ \sigma\subseteq W\bigr\}
\]
be the induced subcomplex on $W$ at scale~$\varepsilon$.  Hochster’s formula expresses the graded
Betti numbers of $k[\Delta^{\varepsilon}]$ in purely topological terms:
\begin{equation}
\beta_{i,i+j}\bigl(k[\Delta^{\varepsilon}]\bigr)
\;=\;
\sum_{\substack{W\subseteq V\\ |W|=i+j}}
\dim_{k}\,\widetilde{H}_{j-1}\!\bigl(\Delta^{\varepsilon}_W;k\bigr),
\qquad j\ge 1,
\label{eq:hochster_long}
\end{equation}
where $\widetilde{H}_{j-1}(\Delta^{\varepsilon}_{W}; k)$ denotes the $(j-1)$–st reduced simplicial homology group of the induced subcomplex $\Delta^{\varepsilon}_{W}$, computed over the coefficient field $k$.

The formula can also be decomposed into its componentwise forms:
\begin{align}
\beta_{i,i+1}\bigl(k[\Delta^{\varepsilon}]\bigr)
&=\sum_{\substack{W\subseteq V\\ |W|=i+1}}
\Bigl(\beta_0\bigl(\Delta^{\varepsilon}_W\bigr)-1\Bigr),
\label{eq:hochster_long_dim0}\\[2pt]
\beta_{i,i+j}\bigl(k[\Delta^{\varepsilon}]\bigr)
&=\sum_{\substack{W\subseteq V\\ |W|=i+j}}
\beta_{j-1}\bigl(\Delta^{\varepsilon}_W\bigr),
\qquad j\ge 2.
\label{eq:hochster_long_dimj}
\end{align}
Equation~\eqref{eq:hochster_long_dim0} measures the lack of connectivity in induced 
$(i{+}1)$–vertex subcomplexes, while equation \eqref{eq:hochster_long_dimj} aggregates higher-dimensional homology across 
$(i{+}j)$–vertex subcomplexes.

To incorporate persistence, fix $\varepsilon\le \varepsilon'$ and note that for each $W\subseteq V$
there is an inclusion
$\Delta^{\varepsilon}_W \hookrightarrow \Delta^{\varepsilon'}_W$
inducing a homomorphism on reduced homology
$\iota^{\,j-1}(\varepsilon,\varepsilon'):\ 
\widetilde{H}_{j-1}\bigl(\Delta^{\varepsilon}_W;k\bigr)
\longrightarrow
\widetilde{H}_{j-1}\bigl(\Delta^{\varepsilon'}_W;k\bigr).$
We define the persistent graded Betti numbers by
\begin{equation}
\beta^{\,\varepsilon,\varepsilon'}_{i,i+j}\bigl(k[\Delta]\bigr)
\;:=\;
\sum_{\substack{W\subseteq V\\ |W|=i+j}}
\dim_{k}\!\Bigl(\operatorname{im}\,\iota^{\,j-1}(\varepsilon,\varepsilon')\Bigr).
\label{eq:persistent_betti_long}
\end{equation}
These graded invariants record not only the persistence of homological classes across
$[\varepsilon,\varepsilon']$, but also the combinatorial structure of their supports.  In the special case when $|W|=|V|$, one recovers the classical
persistent Betti number in homological degree $|V|-i-1$. More generally, the family
$\{\beta^{\,\varepsilon,\varepsilon'}_{i,i+j}\}_{i,j}$ interpolates between topological persistence and 
the graded algebra of $k[\Delta]$.

\subsection{Persistent \texorpdfstring{$f$}{f}- and \texorpdfstring{$h$}{h}-Vectors}

Let $d=\dim(\Delta)+1$.  For each scale $\varepsilon$, the $f$–vector and $h$–vector of
$\Delta^{\varepsilon}$ are
\[
f(\Delta^{\varepsilon}) \;=\; (f^{\varepsilon}_{-1},f^{\varepsilon}_0,\ldots,f^{\varepsilon}_{d-1}),
\qquad
h(\Delta^{\varepsilon}) \;=\; (h^{\varepsilon}_0,\ldots,h^{\varepsilon}_d),
\]
where $f^{\varepsilon}_{i}$ counts the number of $i$–faces of $\Delta^{\varepsilon}$ 
(with $f^{\varepsilon}_{-1}=1$ by convention).  They are connected by the standard linear transform
\begin{equation}
h^{\varepsilon}_m \;=\; \sum_{j=0}^m \binom{d-j}{m-j}(-1)^{m-j} f^{\varepsilon}_{j-1},
\qquad
f^{\varepsilon}_{m-1} \;=\; \sum_{i=0}^m \binom{d-i}{m-i} h^{\varepsilon}_i,
\quad m=0,\ldots,d.
\label{eq:f_h_transform_long}
\end{equation}
The first identity expresses the $h$–vector as an alternating sum of face counts. The inverse relation
reconstructs $f$ from $h$.
Between two scales $\varepsilon\le \varepsilon'$, we define the persistent $h$–vector
\begin{equation}
h^{\,\varepsilon,\varepsilon'}_m
\;:=\;
\sum_{j=0}^m \binom{n-d+m-j-1}{m-j}
\left(\ \sum_{i=0}^j (-1)^i \,\beta^{\,\varepsilon,\varepsilon'}_{i,j}\ \right),
\qquad m=0,\ldots,d,
\label{eq:persistent_h_long}
\end{equation}
and obtain the persistent $f$–vector by the same inverse transform
\begin{equation}
f^{\,\varepsilon,\varepsilon'}_{m-1}
\;:=\;
\sum_{i=0}^m \binom{d-i}{m-i}\, h^{\,\varepsilon,\varepsilon'}_i,
\qquad m=0,\ldots,d.
\label{eq:persistent_f_long}
\end{equation}
Whereas the static $f$– and $h$–vectors summarize the combinatorics of one complex
$\Delta^{\varepsilon}$, their persistent counterparts encapsulate how face counts and the algebraic
relations among faces evolve and persist across scales.  In this way, persistent $f$– and $h$–vectors
blend enumerative combinatorics with homological persistence, yielding multiscale
algebraic–combinatorial invariants within the persistent Stanley--Reisner framework.

\begin{figure}[t] \centering \includegraphics[width=.2\linewidth]{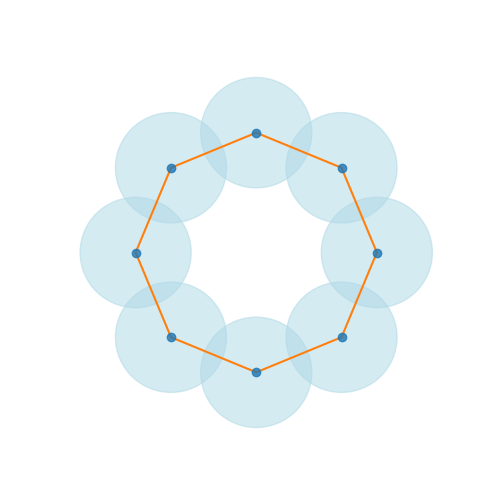} \includegraphics[width=.2\linewidth]{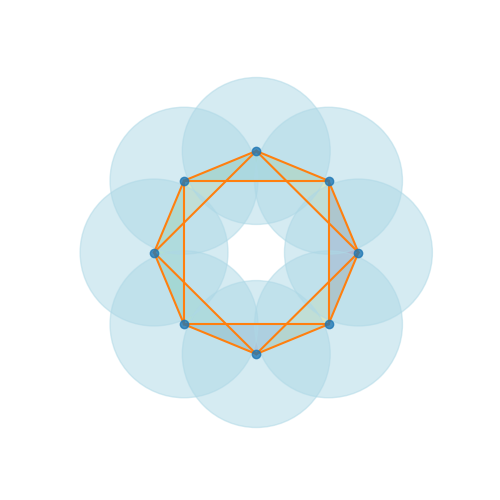} \includegraphics[width=.2\linewidth]{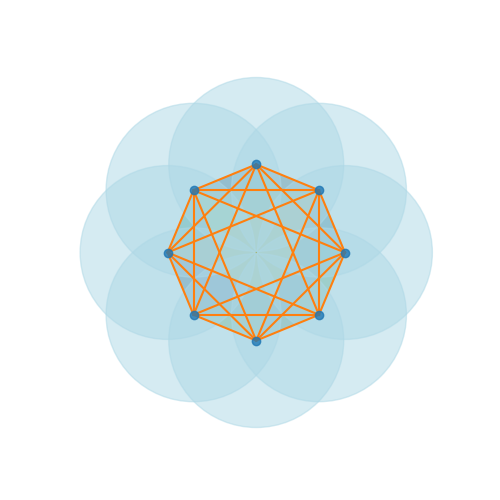} \includegraphics[width=.2\linewidth]{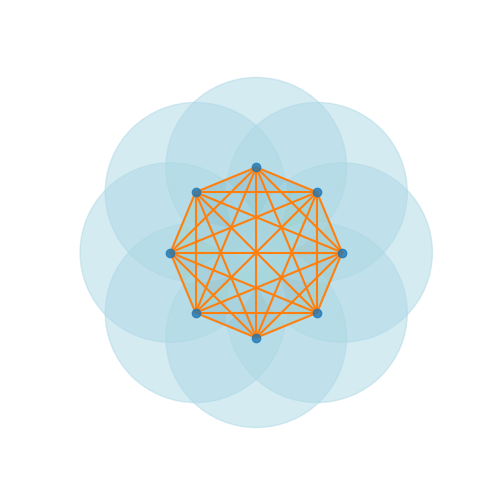} \includegraphics[width=0.8\linewidth, height=8cm]{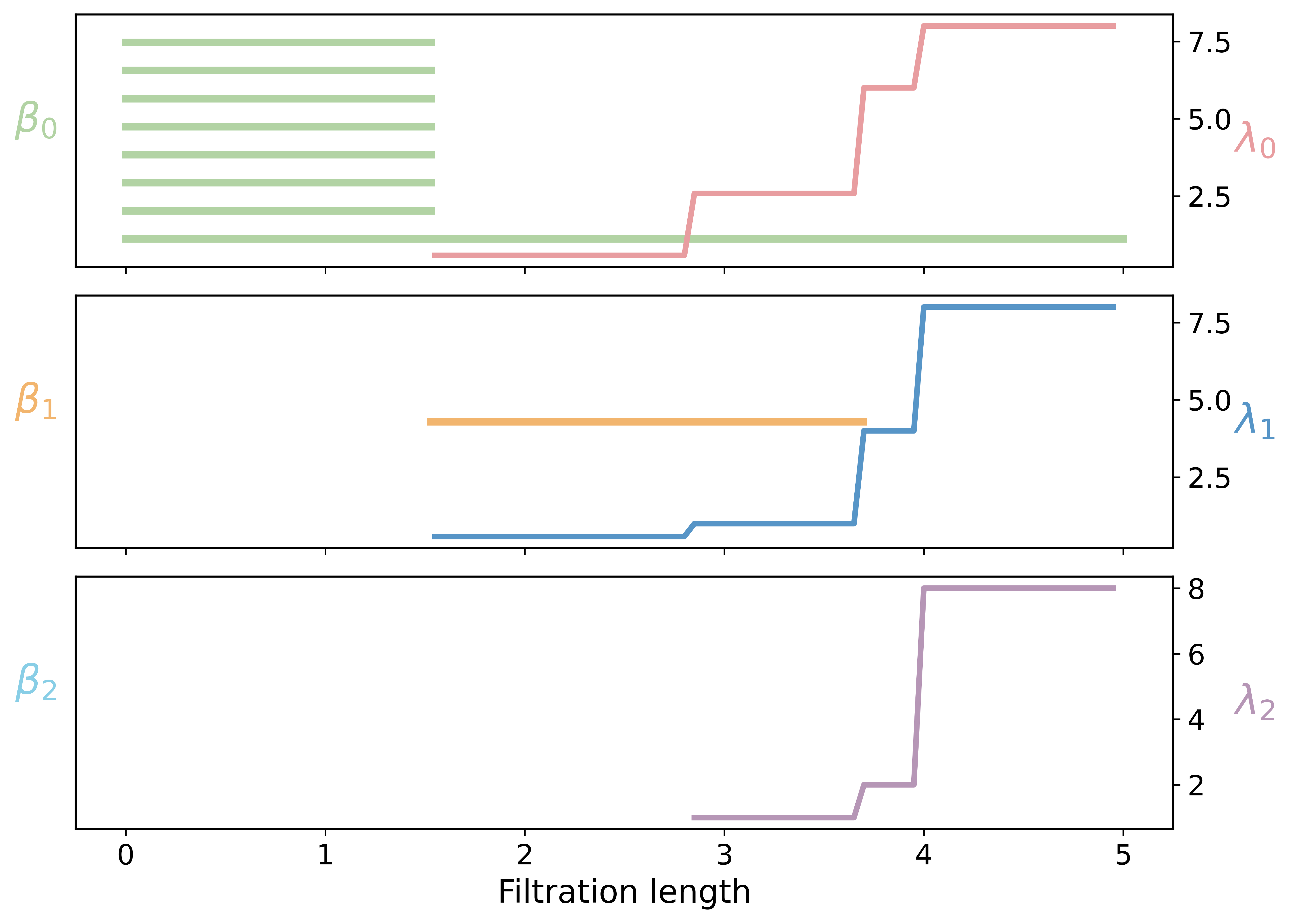} \caption{Illustration of PH and PL on the regular octagon. Top: Vietoris--Rips filtration process. Bottom: PL and PH for $0$-, $1$-, and $2$-dimensional features. } \label{fig:octagon_PL} 
\end{figure}

\begin{figure}[htbp!]
  \centering
  \begin{subfigure}{0.52\linewidth}
    \centering\includegraphics[width=\linewidth]{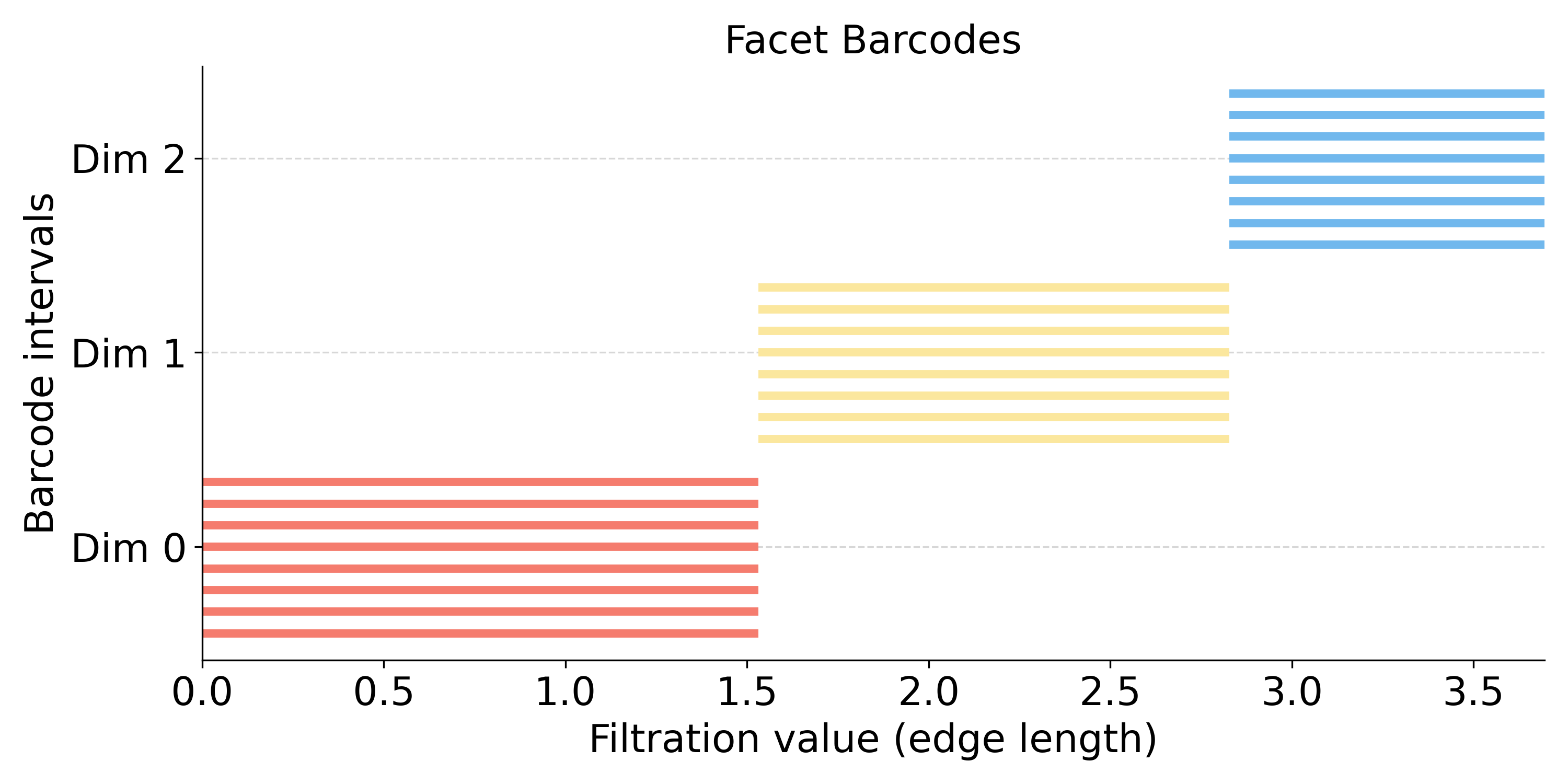}
    \subcaption*{(a) Facet persistence barcode.}
  \end{subfigure}\hfill
  \begin{subfigure}{0.48\linewidth}
    \centering\includegraphics[width=\linewidth]{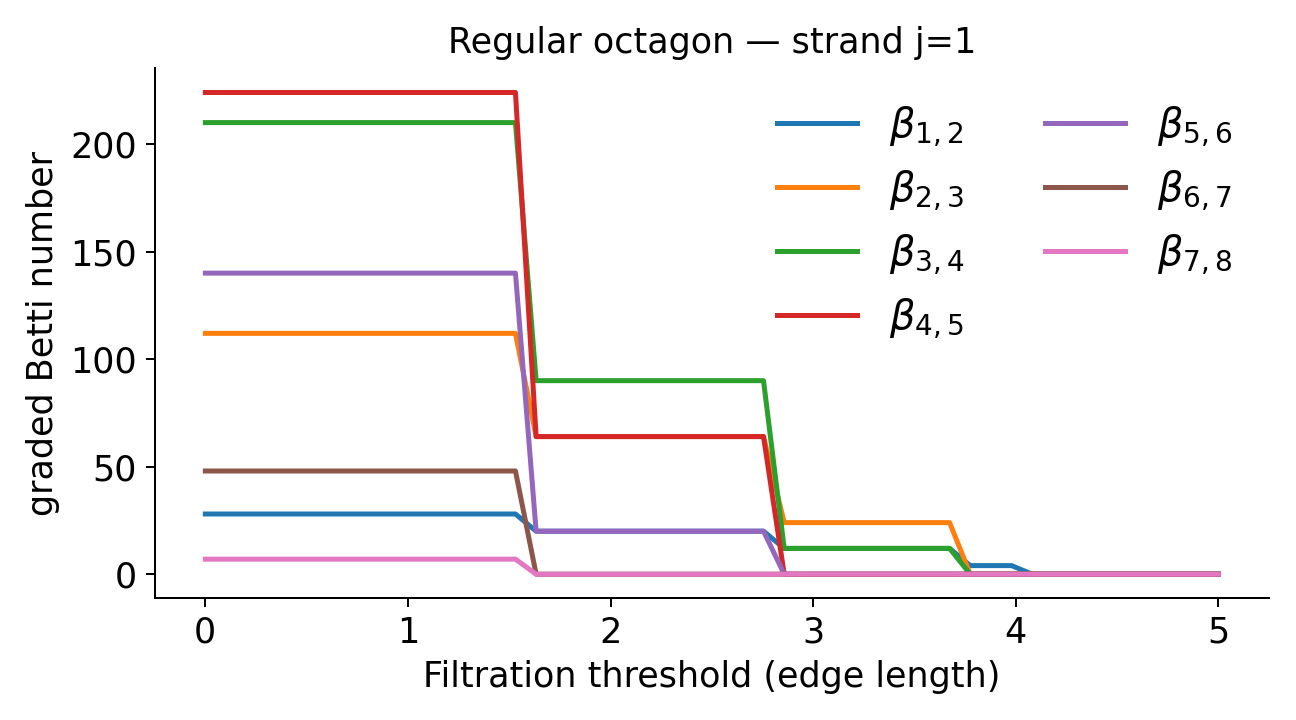}
    \subcaption*{(b) Graded Betti curves, dim$=0$ (strand $j=1$).}
  \end{subfigure}
  \vspace{0.5em}
  \begin{subfigure}{0.5\linewidth}
    \centering\includegraphics[width=\linewidth]{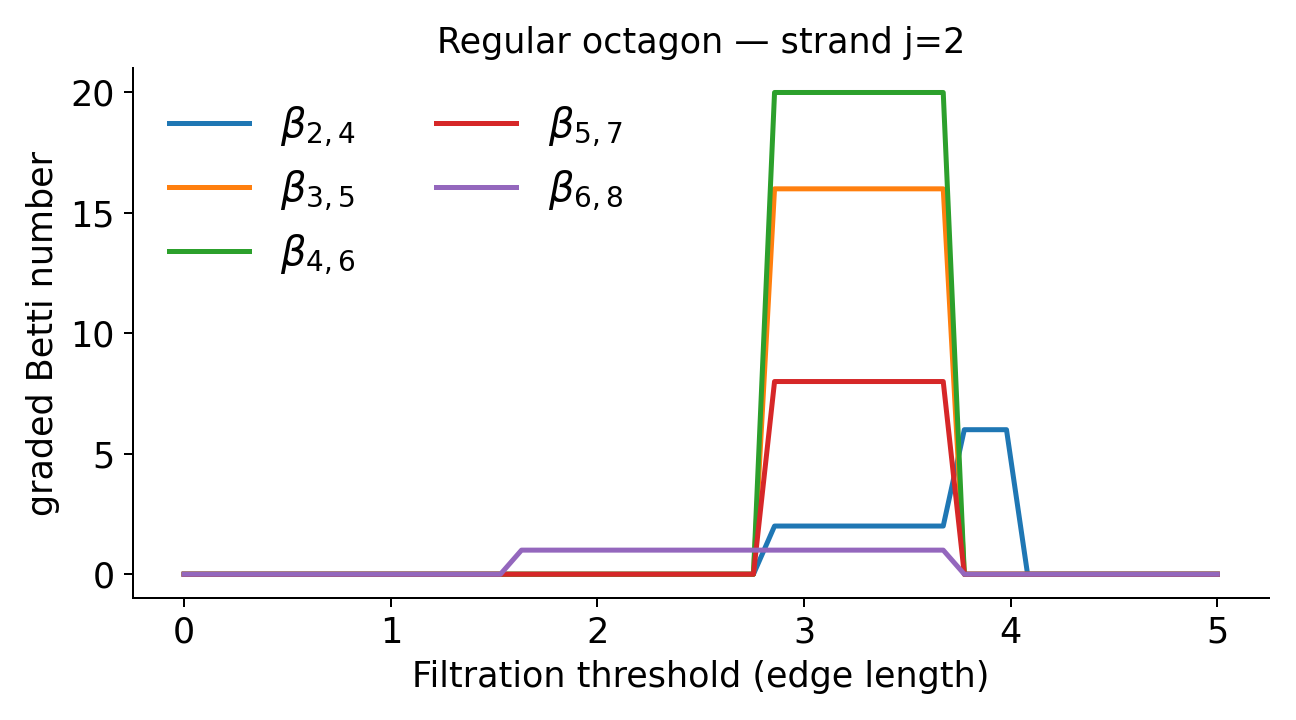}
    \subcaption*{(c) Graded Betti curves, dim$=1$ (strand $j=2$).}
  \end{subfigure}\hfill
  \begin{subfigure}{0.5\linewidth}
    \centering\includegraphics[width=\linewidth]{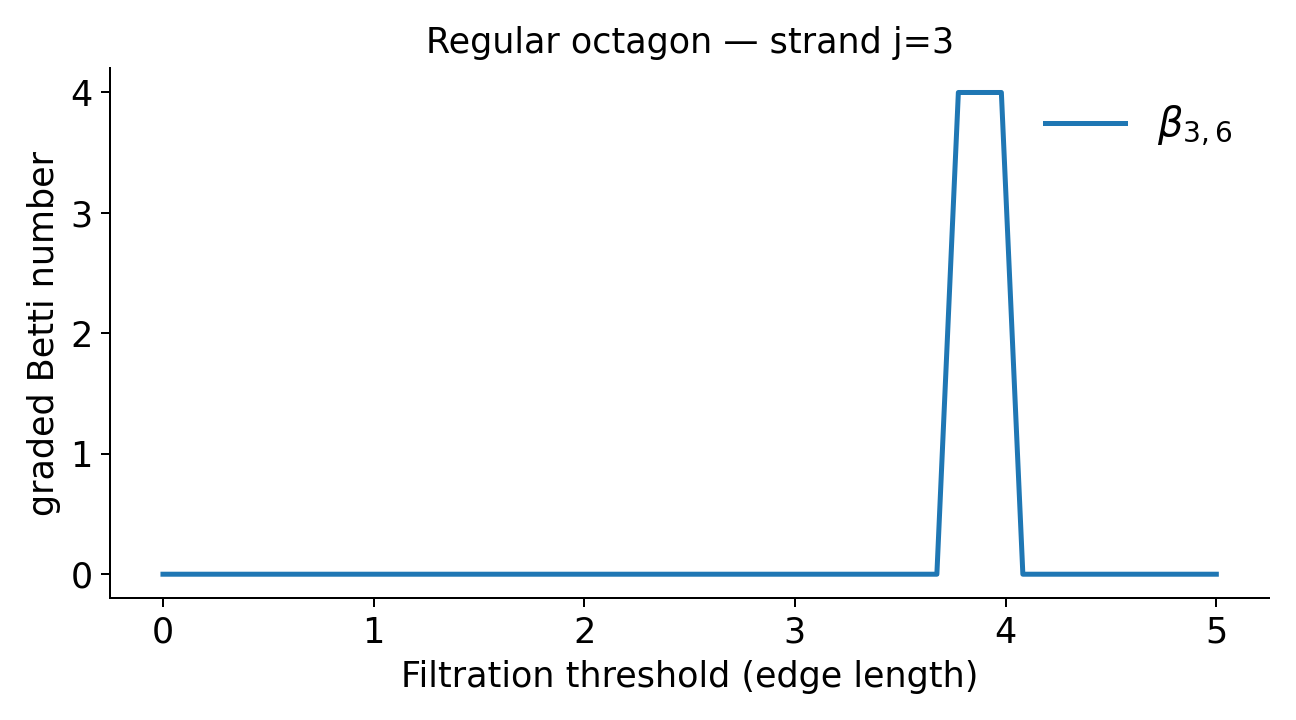}
    \subcaption*{(d) Graded Betti curves, dim$=2$ (strand $j=3$).}
  \end{subfigure}
  \vspace{0.5em}
  \begin{subfigure}{0.48\linewidth}
    \centering\includegraphics[width=\linewidth]{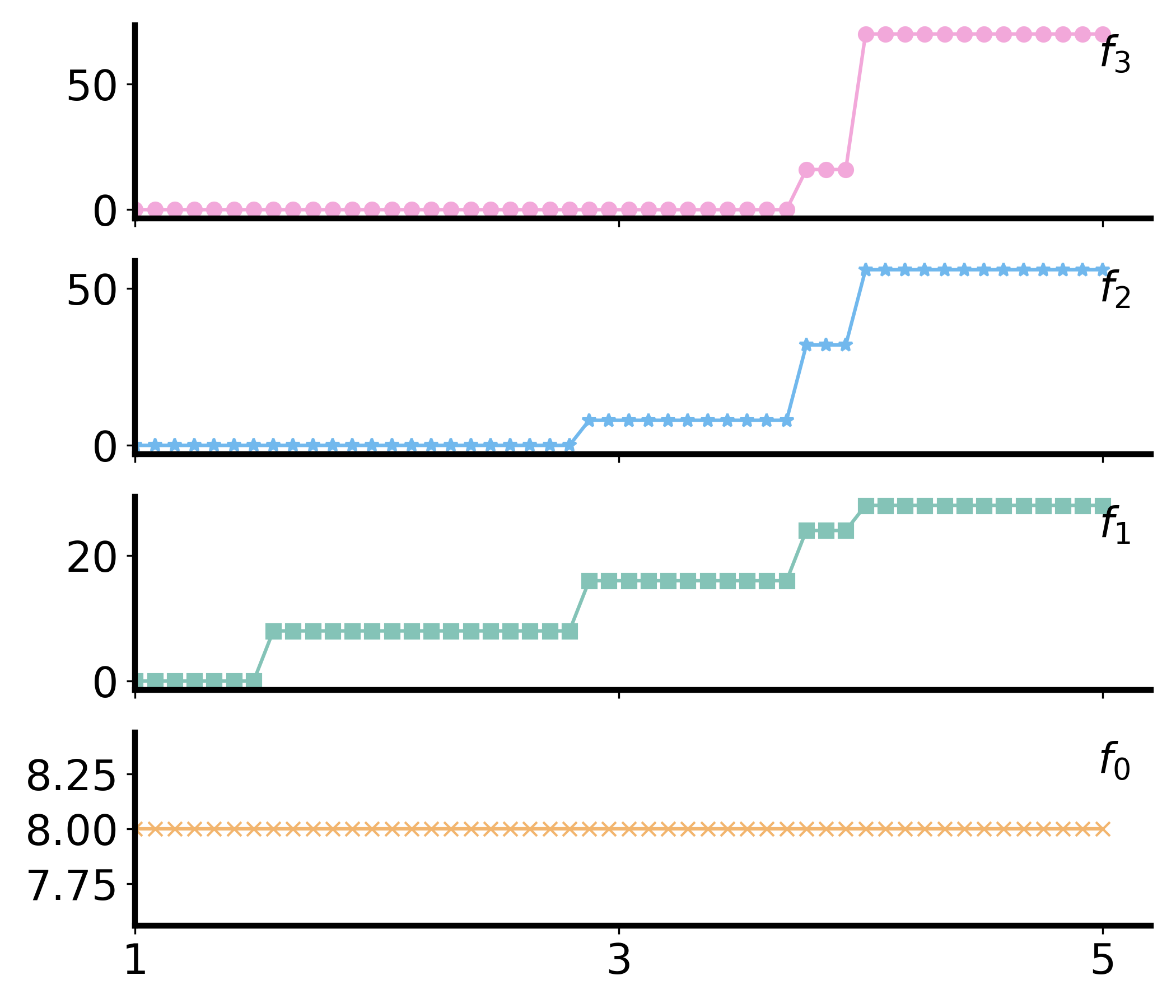}
    \subcaption*{(e) $f$‑vector curves.}
  \end{subfigure}\hfill
  \begin{subfigure}{0.48\linewidth}
    \centering\includegraphics[width=\linewidth]{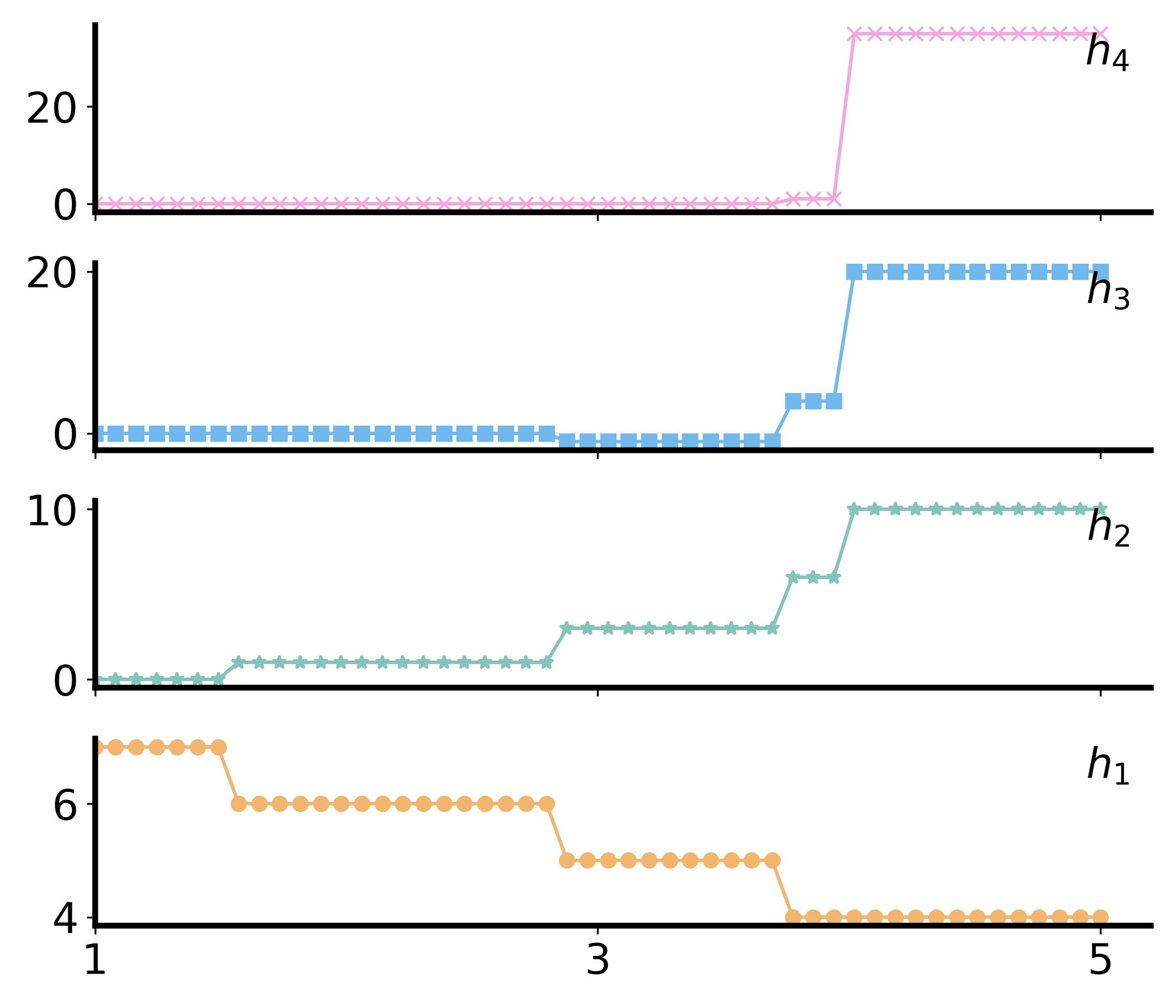}
    \subcaption*{(f) $h$‑vector curves.}
  \end{subfigure}

  \caption{Illustrations of persistent commutative algebra  analysis on the regular octagon using a Rips complex-based filtration process.}
  \label{fig:octagon_CA}
\end{figure}

\section{Interpretability and Representability of PH, PL and PCA}

Building on the constructions in previous sections, we compare three lenses on the same Vietoris--Rips filtration \cite{vietoris1927hoheren} across all datasets, ranging from a benchmark example to a high–dimensional synthetic complex, a fullerene structure, and biomolecular systems.  First, we examine PH through barcodes of the Betti numbers \(\beta_{0}\), \(\beta_{1}\), and \(\beta_{2}\). These intervals track the birth and death of independent connected components, loops, and voids as the filtration radius increases. Second, we study PL through the smallest nonharmonic eigenvalue curves \(\lambda_{0}\), \(\lambda_{1}\), and \(\lambda_{2}\). These spectra quantify how strongly the underlying geometry sustains the observed features and provide sensitivity beyond raw Betti counts. 
Third, within the Stanley--Reisner framework of PCA, we analyze graded Betti numbers, facet persistence, and the \(f\)– and \(h\)–vectors. The graded Betti numbers are algebraic invariants from a minimal free resolution. By Hochster’s formula, they admit a topological interpretation as sums of reduced homology dimensions of induced subcomplexes \((\Delta^{\varepsilon})_{W}\) as in Eq.~\ref{eq:hochster_long}, in which the strand index \(j\) aligns with homological dimension \(j-1\). In this section, we report graded Betti numbers on strands $j=1,2,3$ for dimensions $0,1,2$.  The \(f\)–vector counts simplices by dimension at each filtration scale, and the \(h\)–vector follows from the standard binomial relations applied to the \(f\)–vector as in Eq.~\ref{eq:f_h_transform_long}. The facet persistence barcode records when a simplex is maximal in the complex and when it ceases to be maximal after the inclusion of a strict superset. 


\subsection{Regular Octagon: A Benchmark Example}

We study the Vietoris--Rips filtration on eight planar points positioned at the vertices of a regular octagon with circumradius $R = 2$. 
Four critical distances determine the evolution of the filtration:
\[
r_1 = 2R\sin\!\left(\frac{\pi}{8}\right) \approx 1.53~\text{\AA}, \qquad
r_2 = 2R\sin\!\left(\frac{\pi}{4}\right) \approx 2.83~\text{\AA},
\]
\[
r_3 = 2R\sin\!\left(\frac{3\pi}{8}\right) \approx 3.70~\text{\AA}, \qquad
r_4 = 2R\sin\!\left(\frac{\pi}{2}\right) = 4.00~\text{\AA}.
\]

As the Vietoris--Rips threshold $\varepsilon$ increases, the complex evolves through these four characteristic radii, as illustrated in the top panel of Fig.~\ref{fig:octagon_PL}. 
At $r_1$, connections form between adjacent vertices, outlining the perimeter of the octagon. 
At $r_2$, diagonals connecting every second vertex appear, allowing each set of three consecutive vertices to form a triangle. 
When $\varepsilon$ reaches $r_3$, edges connecting vertices separated by three steps along the octagon are added, bridging across the interior and creating higher-dimensional simplices.
Finally, at $r_4$, the remaining diagonals emerge and all vertices become mutually connected. 
Beyond this point, the complex rapidly fills in to form the complete simplex on eight vertices, where every subset of vertices defines a simplex.

The bottom panel of Fig.~\ref{fig:octagon_PL} summarizes the features of PH and PL through three dimensions for the octagon. The harmonic spectra of PL are isomorphic to PH, and they track the births and deaths of connected components, loops, and voids across the filtration. Eight \(\beta_{0}\) bars, corresponding to eight isolated vertices, appear. Seven die at \(r_{1}\), one bar remains for the global component, and no further zero-dimensional topological changes occur. Nonharmonically, \(\lambda_{0}\) still increases in steps after \(r_{1}\), reflecting stronger overall connectivity as more edges are added. In dimension one, a single loop appears at \(r_{1}\) and persists until \(r_{3}\). Hence, there is one \(\beta_{1}\) bar on \([r_{1},r_{3})\).  Triangles created at \(r_{2}\) reinforce the boundary without yet eliminating this loop. It is filled only when tetrahedra first become available at \(r_{3}\). The nonharmonic spectrum \(\lambda_{1}\) rises near \(r_{2}\) as many triangles stiffen the one–dimensional structure. It continues to change toward \(r_{4}\) even after the \(\beta_{1}\) bar has died, registering geometric consolidation that barcodes alone do not record. No \(\beta_{2}\) bars appear, yet the rise of \(\lambda_{2}\) still reflects strengthening of two-dimensional couplings induced by higher-dimensional simplices. In short, PL encodes geometric consolidation that barcodes alone do not reveal.


The commutative–algebra viewpoint complements PH and PL for the octagon, in Fig.~\ref{fig:octagon_CA}, by examining facet ideals, together with the \(f\)– and \(h\)–vectors, and the graded Betti numbers. In Fig.~\ref{fig:octagon_CA}(a), the eight red facet bars in dimension \(0\) encode the eight isolated vertices and persist until \(r_{1}\). Yellow facet bars in dimension \(1\) terminate until \(r_{2}\), when the first short diagonals enter and every edge is contained in at least one triangle. Blue facet bars in dimension \(2\) begin at \(r_{2}\), from the new triangles, and persist until \(r_{3}\). Longer diagonals, at \(r_{3}\), create tetrahedra, after which triangles are no longer maximal. By \(r_{4}\), all low–dimensional facets are contained in higher simplices, and vanish. The \(f\)–vector curves, in Fig.~\ref{fig:octagon_CA}(e), count simplices by dimension across the filtration. The value \(f_{0}\) remains at eight. The curve \(f_{1}\) rises at \(r_{1}\), when polygon edges appear, then grows again, together with \(f_{2}\), at \(r_{2}\), when short diagonals and triangles enter, and increases further, together with \(f_{2}\) and \(f_{3}\), at \(r_{3}\), as longer diagonals and tetrahedra are added. At \(r_{4}\), the complex becomes the full simplex on eight vertices, and all \(f_{k}\) reach their maximal values. The \(h\)–vector curves, in Fig.~\ref{fig:octagon_CA}(f), follow from the binomial transform of the \(f\)–vector, and respond consistently to these transitions.

The PCA framework captures the same evolution of filtration for the octagon through graded Betti numbers. In Fig.~\ref{fig:octagon_CA}(b–d), the strands \(j=1,2,3\) aggregate contributions from induced subcomplexes \((\Delta^{\varepsilon})_{W}\) that carry reduced homology in dimensions \(0,1,2\), respectively, in accordance with Hochster’s formula (Eq.~\ref{eq:hochster_long}). In dimension \(0\), the graded Betti curves remain at their combinatorial plateaus until \(r_{1}\). At \(r_{1}\), the perimeter edges appear, most induced subcomplexes become connected, and the curves drop sharply. The emergence of many triangles further improves local connectivity and produces another step down at \(r_{2}\). The newly added longer diagonals eliminate the remaining disconnected configurations near \(r_{3}\), and all \(j=1\) curves vanish by \(r_{4}\). In dimension \(1\), A single loop exists on the full octagon from \(r_{1}\) to \(r_{3}\), and in parallel, many induced subcomplexes on \(4\)–\(7\) vertices support short cycles between \(r_{1}\) and \(r_{3}\). This yields rises on \(j=2\) that are most pronounced after \(r_{1}\) and again near \(r_{2}\) as triangles consolidate local structure but do not yet fill every cycle. At \(r_{3}\), the longer edges create enough \(2\)–faces within the relevant induced subcomplexes to cap these loops, and the \(j=2\) curves rapidly decay to zero. A short peak in \(\beta_{2,4}\) appears exactly at \(r_{3}\), when a four–edge shell exists, but its internal diagonals have not fully formed. In dimension \(2\), a brief signal can appear near \(r_{3}\) when certain \(6\)–vertex subsets momentarily support triangular shells before higher–dimensional simplices are available. As \(\varepsilon\) increases past \(r_{3}\), the same subsets admit \(3\)–simplices that fill these shells, so the \(j=3\) signal disappears by \(r_{4}\).


\subsection{Octahedron: Higher-Dimensional Structure}
\begin{figure}
    \centering
     \includegraphics[width=0.9\linewidth]{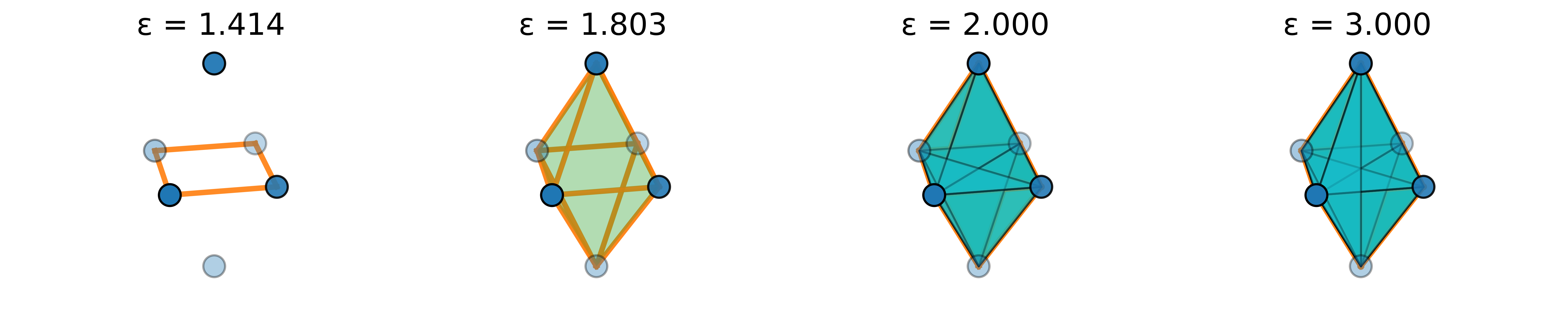}
     \includegraphics[width=0.8\linewidth, height=8cm]{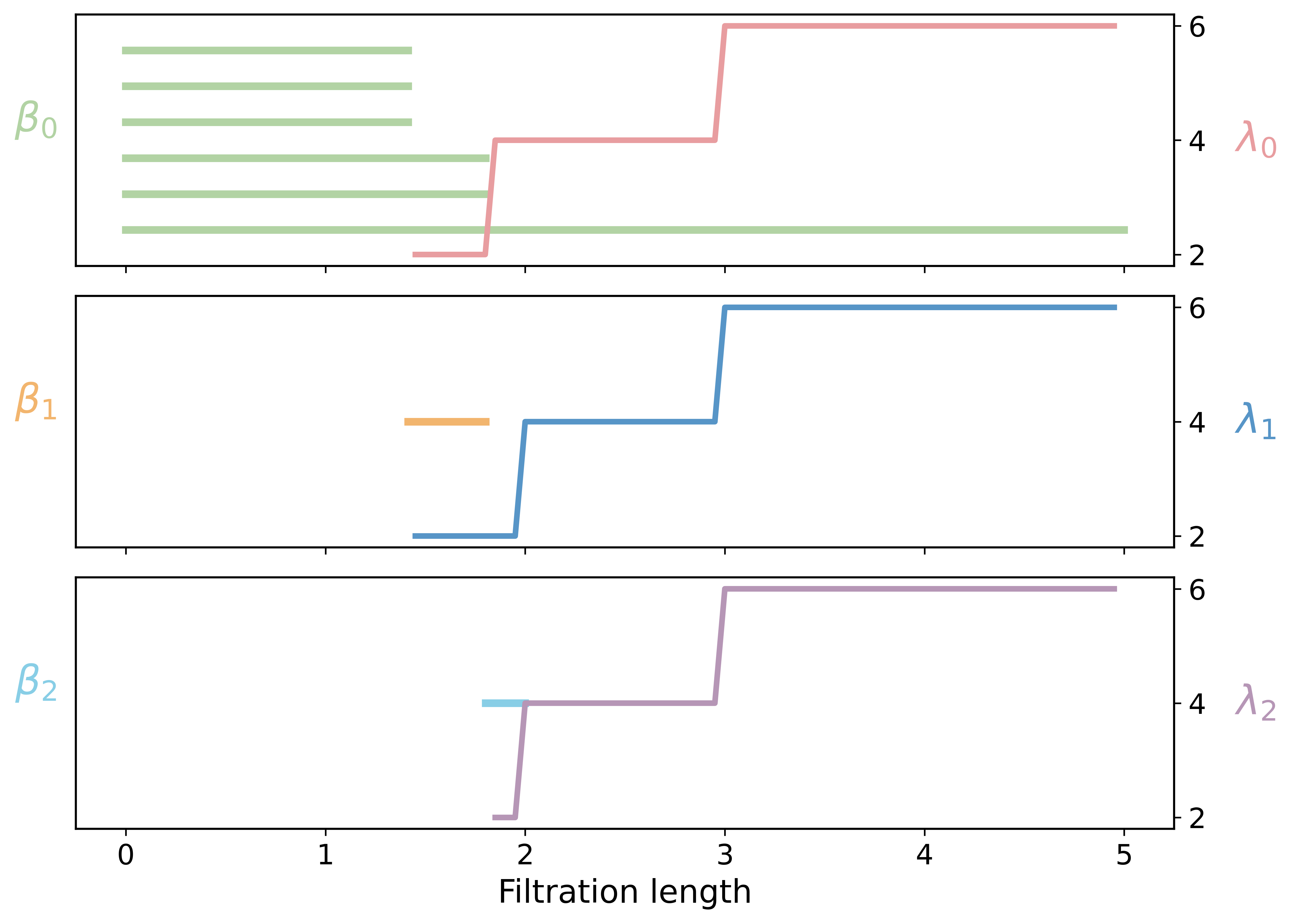}
     \caption{{Illustration of PH and PL on the regular octahedron. Top: Vietoris--Rips filtration process. Bottom: PL and PH for $0$-, $1$-, and $2$-dimensional features.}}
         \label{fig:octa_PL}
\end{figure}

We study the Vietoris--Rips filtration on the stretched octahedron,
\[
P=\{(\pm 1,0,0),\ (0,\pm 1,0),\ (0,0,\pm 1.5)\}.
\]
As shown in the top panel of Fig.~\ref{fig:octa_PL}, four critical radii govern the growth of the complex,
\[
r_{1}=\sqrt{2}\approx 1.41~\text{\AA},\quad
r_{2}=\sqrt{3.25}\approx 1.80~\text{\AA},\quad
r_{3}=2.00~\text{\AA},\quad
r_{4}=3.00~\text{\AA}.
\]
At \(r_{1}\), the four equatorial vertices join to form a square in the \(xy\)-plane. At \(r_{2}\), each equatorial vertex connects to both poles, producing eight triangular faces that form the octahedral surface. At \(r_{3}\), the two long diagonals of the equatorial square appear, creating equatorial triangles and tetrahedra that involve one pole and three equatorial vertices. At \(r_{4}\), the two poles connect directly, all edges among the six vertices are present, and the complex becomes the full \(5\)-simplex.

The bottom panel of Fig.~\ref{fig:octa_PL} compiles the PH barcodes and PL spectra for the octahedron for dimensions \(0,1,2\). At \(r_{1}\), the four equatorial vertices connect to form one component while the two poles remain isolated, leaving three components in total. At \(r_{2}\), each pole joins the equatorial set, so the two remaining \(\beta_{0}\) bars die and a single global \(\beta_{0}\) class persists thereafter.  A short \(\beta_{1}\) bar appears on \([r_{1},r_{2})\), corresponding to the equatorial square, and dies at \(r_{2}\) when triangles cap the loop. A \(\beta_{2}\) bar is born at \(r_{2}\) for the hollow triangular shell of the octahedron and dies at \(r_{3}\) once tetrahedra fill the interior. Beyond \(r_{3}\), no further topological changes occur. The nonharmonic PL spectra reveal additional geometric evolution. At \(r_{1}\), the spectra \(\lambda_{0}\) and \(\lambda_{1}\) turn on, as equatorial edges form. At \(r_{2}\), \(\lambda_{0}\) continues to rise, and the \(\lambda_{2}\) activates, with the appearance of the triangular shell. At \(r_{3}\), both \(\lambda_{1}\) and \(\lambda_{2}\) increase sharply, as many triangles are completed and tetrahedra fill the shell. At \(r_{4}\), the complex stabilizes as the complete simplex, and the spectral curves plateau, reflecting the absence of further geometric or topological change.

\begin{figure}[htbp!]
  \centering
  \begin{subfigure}{0.52\linewidth}
    \centering\includegraphics[width=\linewidth]{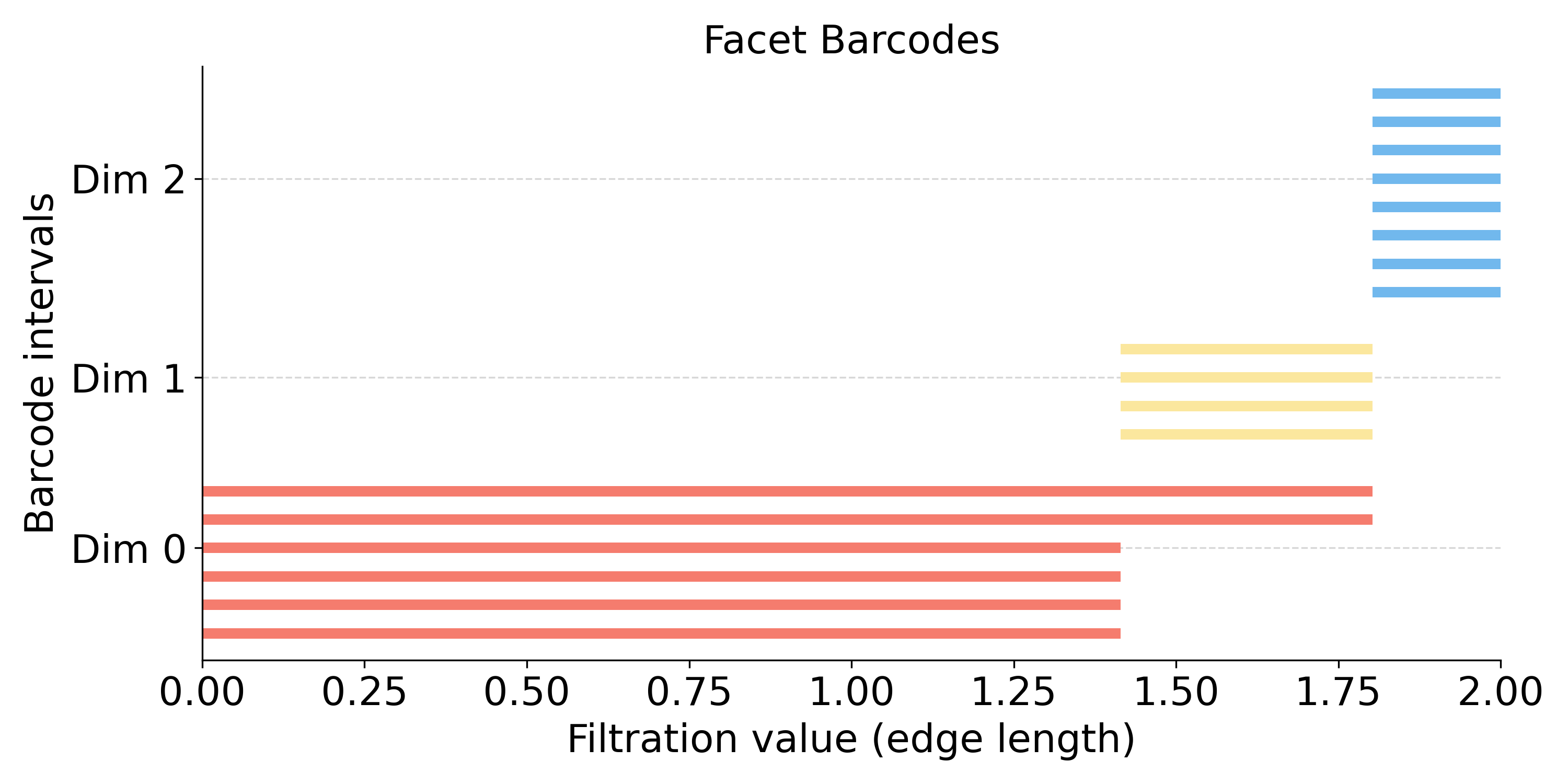}
    \subcaption*{(a) Facet persistence.}
  \end{subfigure}\hfill
  \begin{subfigure}{0.48\linewidth}
    \centering\includegraphics[width=\linewidth]{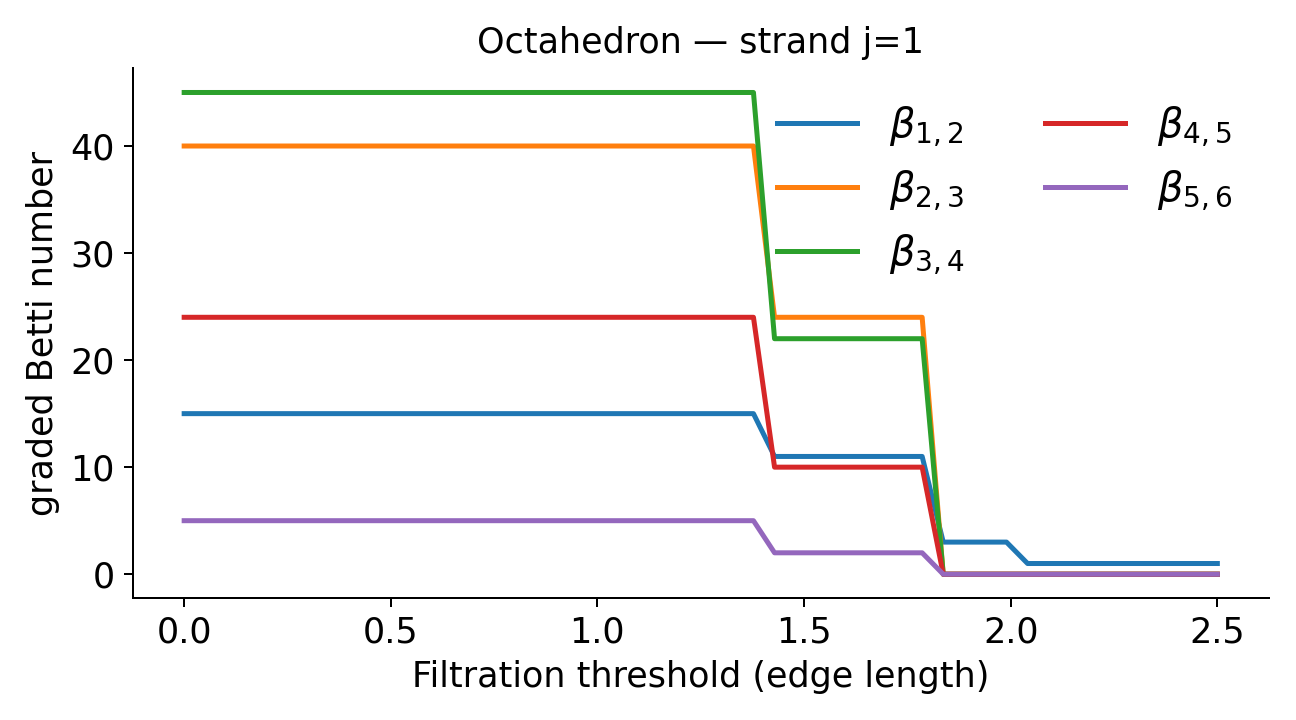}
    \subcaption*{(b) Graded Betti curves, dim$=0$ (strand $j=1$).}
  \end{subfigure}
  \vspace{0.5em}
  \begin{subfigure}{0.5\linewidth}
    \centering\includegraphics[width=\linewidth]{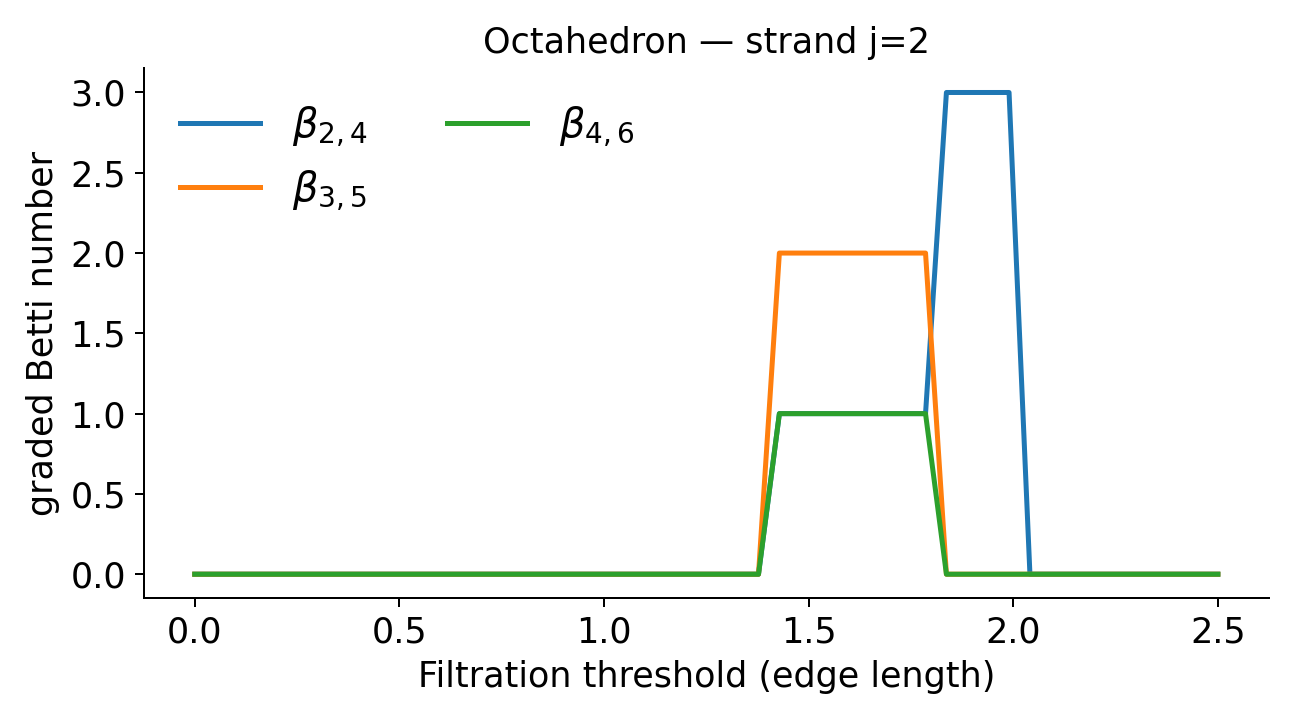}
    \subcaption*{(c) Graded Betti curves, dim$=1$ (strand $j=2$).}
  \end{subfigure}\hfill
  \begin{subfigure}{0.5\linewidth}
    \centering\includegraphics[width=\linewidth]{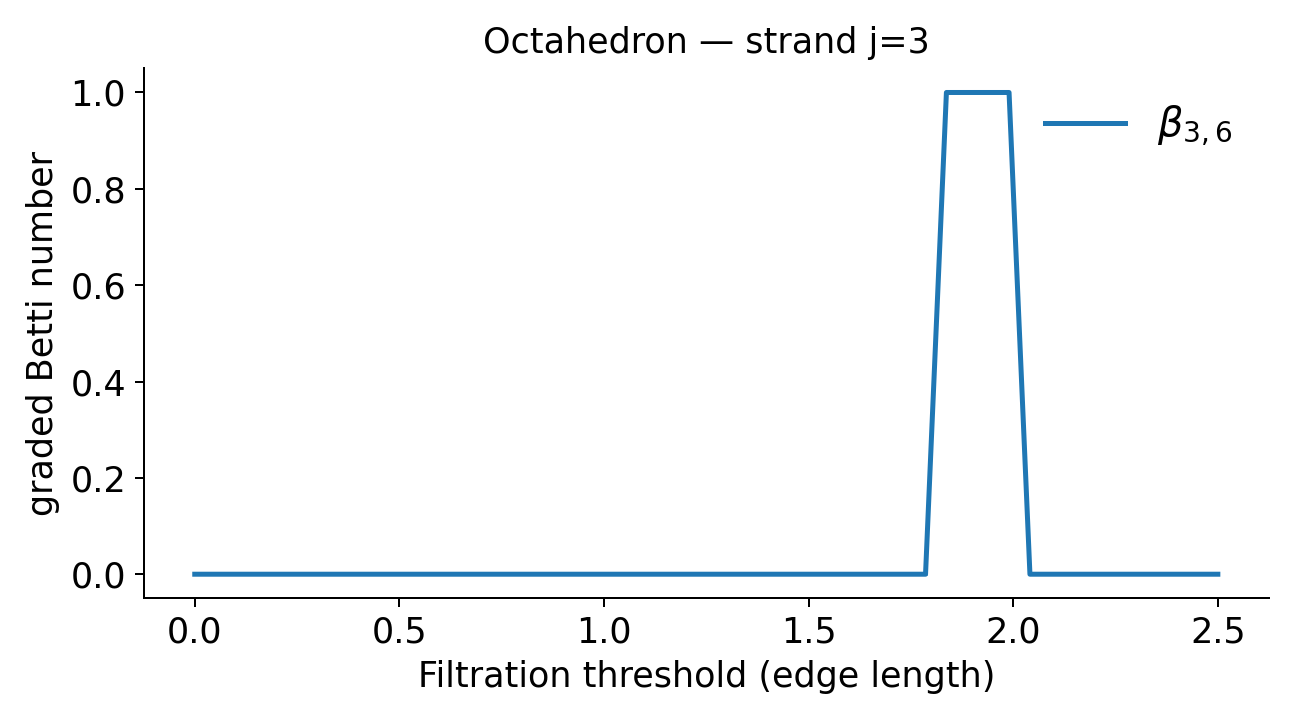}
    \subcaption*{(d) Graded Betti curves, dim$=2$ (strand $j=3$).}
  \end{subfigure}
  \vspace{0.5em}
  \begin{subfigure}{0.48\linewidth}
    \centering\includegraphics[width=\linewidth, height=5.5cm]{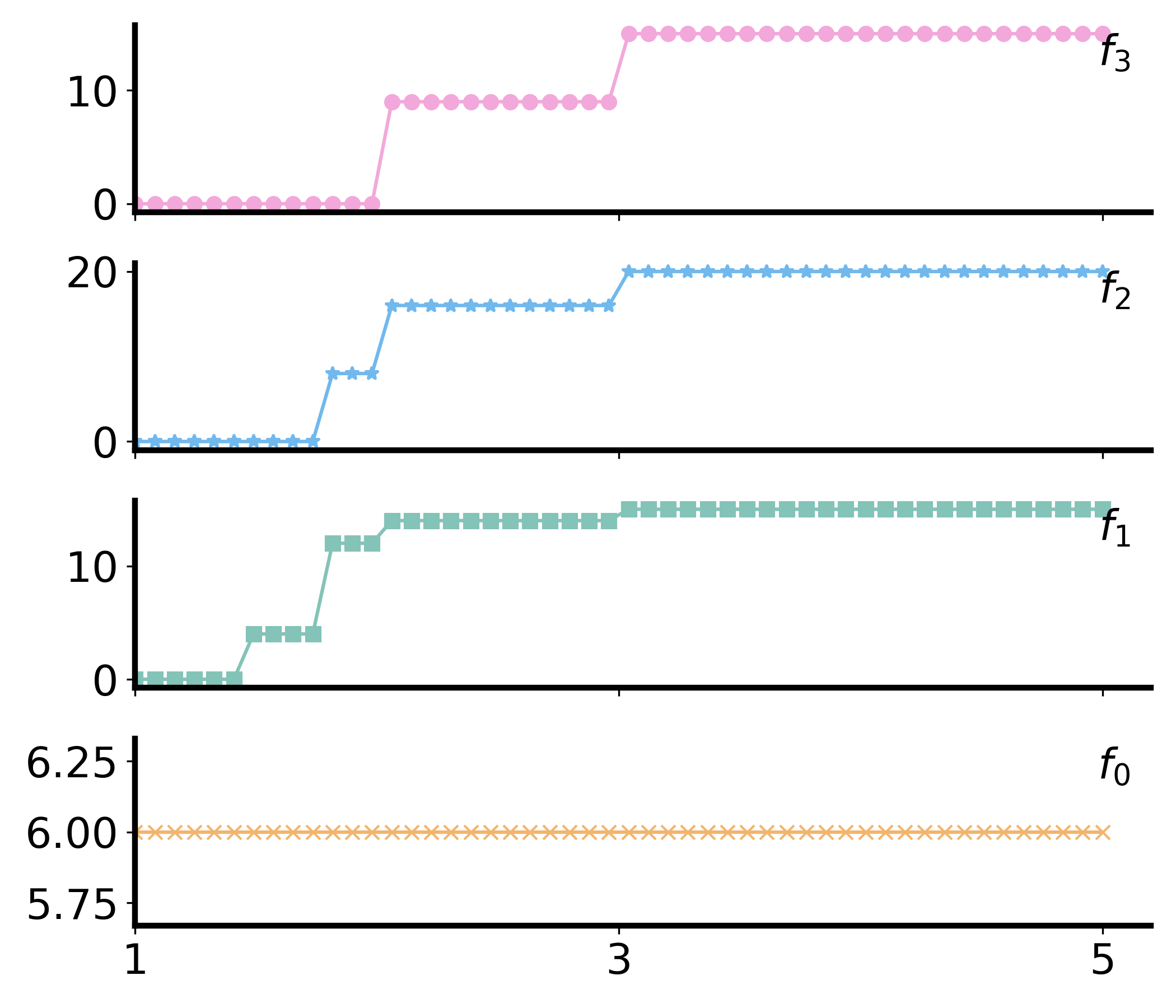}
    \subcaption*{(e) $f$‑vector curves.}
  \end{subfigure}\hfill
  \begin{subfigure}{0.48\linewidth}
    \centering\includegraphics[width=\linewidth, height=5.5cm]{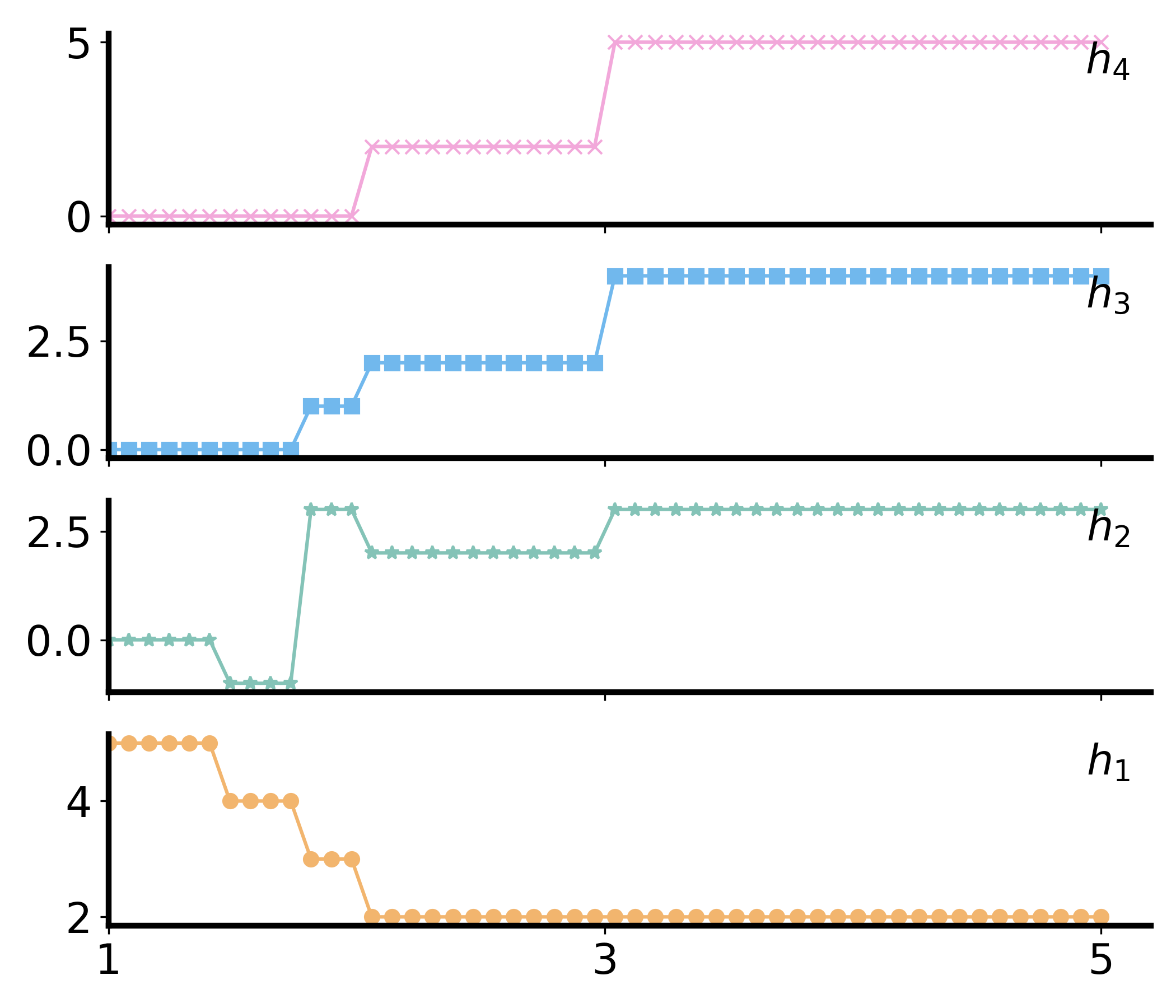}
    \subcaption*{(f) $h$‑vector curves.}
  \end{subfigure}

  \caption{Illustrations of persistent commutative algebra  analysis on the octahedron using a Rips complex-based filtration process.}
  \label{fig:octa_CA}
\end{figure}

Under the PCA framework, four key invariants are tracked across the filtration, and they move in lockstep with the geometry of the stretched octahedron. Facet persistence shows the same sequence of transitions, as in Fig.~\ref{fig:octa_CA}(a), where six long \(0\)–dimensional red bars correspond to six isolated vertices, with the four equatorial vertices persisting until \(r_{1}\), and the two poles disappearing at \(r_{2}\). One–dimensional yellow facet bars vanish at \(r_{2}\), once the equatorial edges participate in pole–equator triangles. These triangles are then absorbed by tetrahedra when the equatorial diagonals are introduced at \(r_{3}\), resulting in the disappearance of eight blue facet bars in dimension 2. The \(f\)– and \(h\)–vector plots capture these combinatorial changes in a stepwise fashion, as in Fig.~\ref{fig:octa_CA}(e) and Fig.~\ref{fig:octa_CA}(f), respectively. The vertex count remains \(f_{0}=6\). The edge count \(f_{1}\) increases by four at \(r_{1}\) from the equatorial square, by eight at \(r_{2}\) from pole–equator connections, by two at \(r_{3}\) from equatorial diagonals, and by one at \(r_{4}\) from the pole–pole edge, reaching the complete value \(15\). The triangle count \(f_{2}\) turns on at \(r_{2}\) with eight pole–equator triangles, rises to \(16\) at \(r_{3}\) with four equatorial and four pole–opposite–equator triangles, and attains the full value \(20\) at \(r_{4}\). Tetrahedra first appear at \(r_{3}\), where nine are present, eight formed by one pole with three equatorial vertices, together with one central equatorial tetrahedra, and they increase to \(15\) at \(r_{4}\). The \(h\)–curves respond accordingly, since the \(h\)–vector is the binomial transform of the \(f\)–vector.

The graded Betti numbers on strands \(j=1,2,3\) offer an algebraic view of the same structural transitions for the octahedron, shown in Fig.~\ref{fig:octa_CA}(b)–(d). They admit a topological interpretation via Hochster’s formula. The strand \(\beta_{i,i+1}\) drops at \(r_{1}\) when equatorial edges form, decreases further at \(r_{2}\) as each equatorial vertex connects to both poles, and falls sharply at \(r_{3}\) when tetrahedra appear, vanishing by \(r_{4}\). The strand \(\beta_{i,i+2}\) is nonzero on \([r_{1},r_{2})\), capturing the equatorial square, then declines at \(r_{2}\) once triangular faces form. Only \(\beta_{2,4}\) shows a brief peak at \(r_{2}\), because three distinct four–vertex cycles exist at that scale. When the two long equatorial diagonals are added at \(r_{3}\), these cycles triangulate immediately, and the signal collapses. The strand \(\beta_{i,i+3}\) remains zero until a closed triangular shell appears, shows a narrow spike just below \(r_{3}\) that reflects the hollow interior, and returns to zero once tetrahedra fill the cavity.

\subsection{Fullerene}

\begin{figure}
    \centering
      \includegraphics[width=0.86\linewidth]{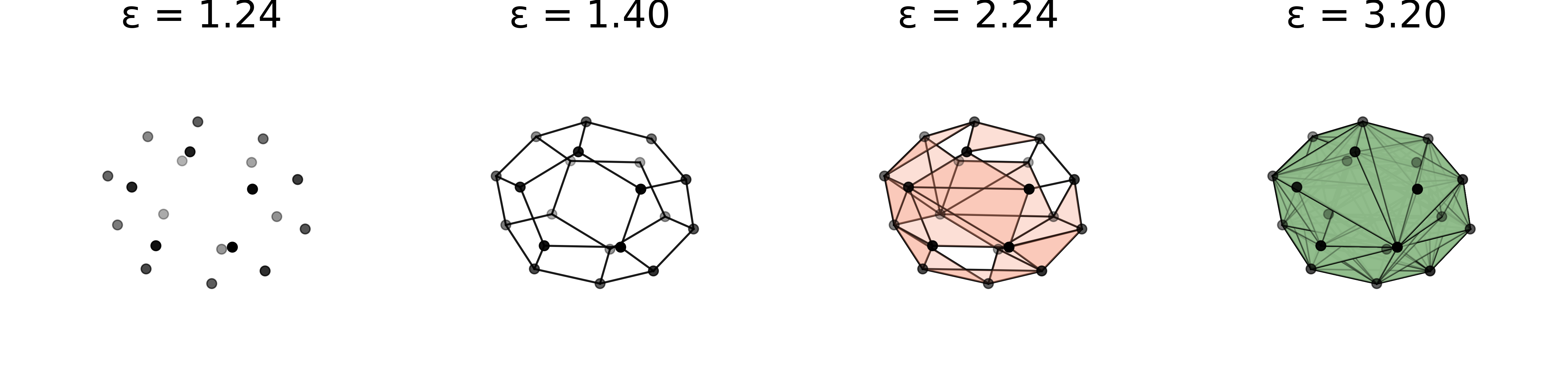}
     \includegraphics[width=0.8\linewidth, height=8cm]{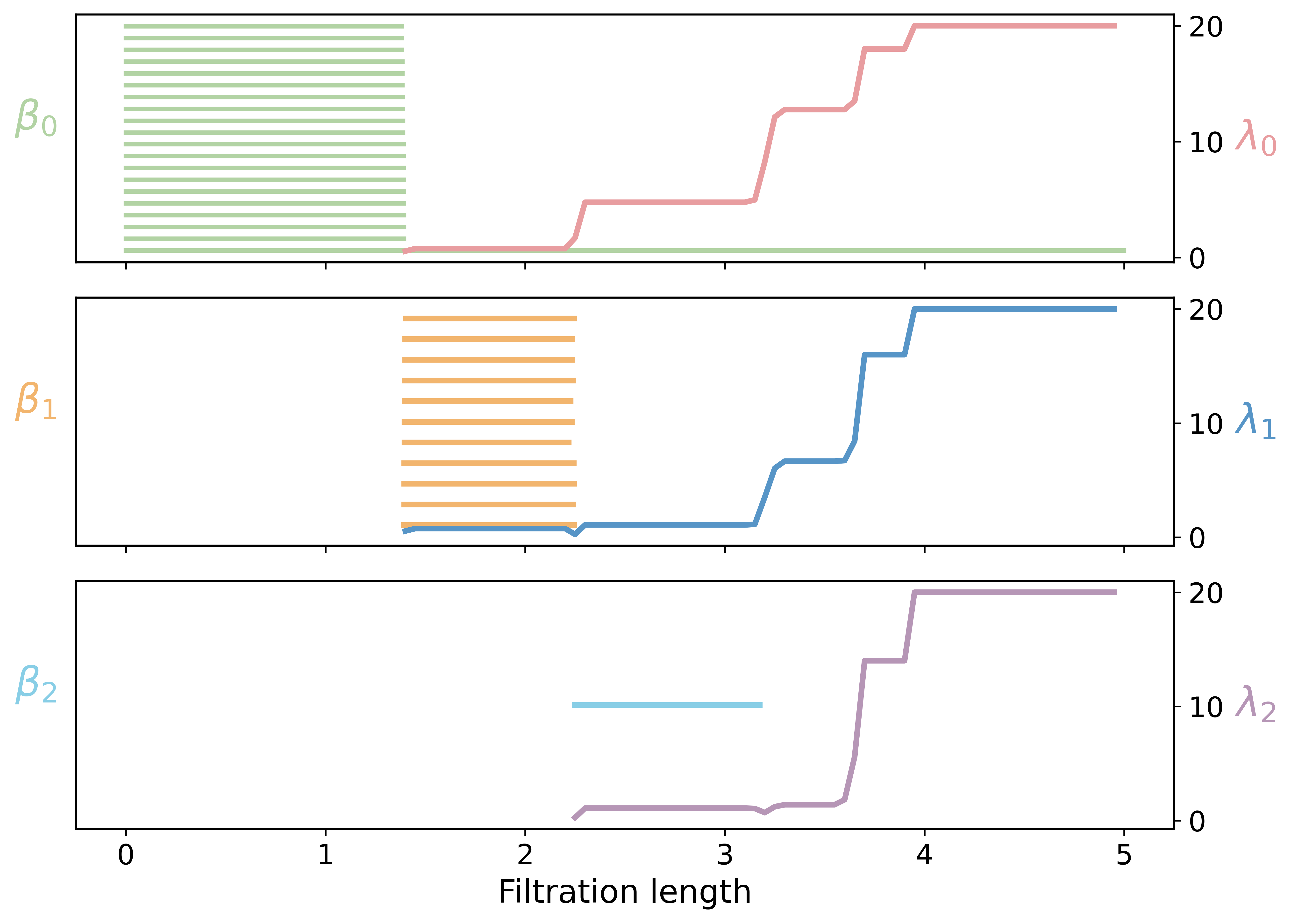}
\caption{Illustration of PH and PL on the fullerene structure C$_{20}$. Top: Vietoris--Rips filtration process. Bottom: PL and PH for $0$-, $1$-, and $2$-dimensional features. }
    \label{fig:phl_c20}
\end{figure}

\begin{figure}[htbp!]
  \centering

  \begin{subfigure}{0.52\linewidth}
    \centering
    \includegraphics[width=\linewidth,height=5cm]{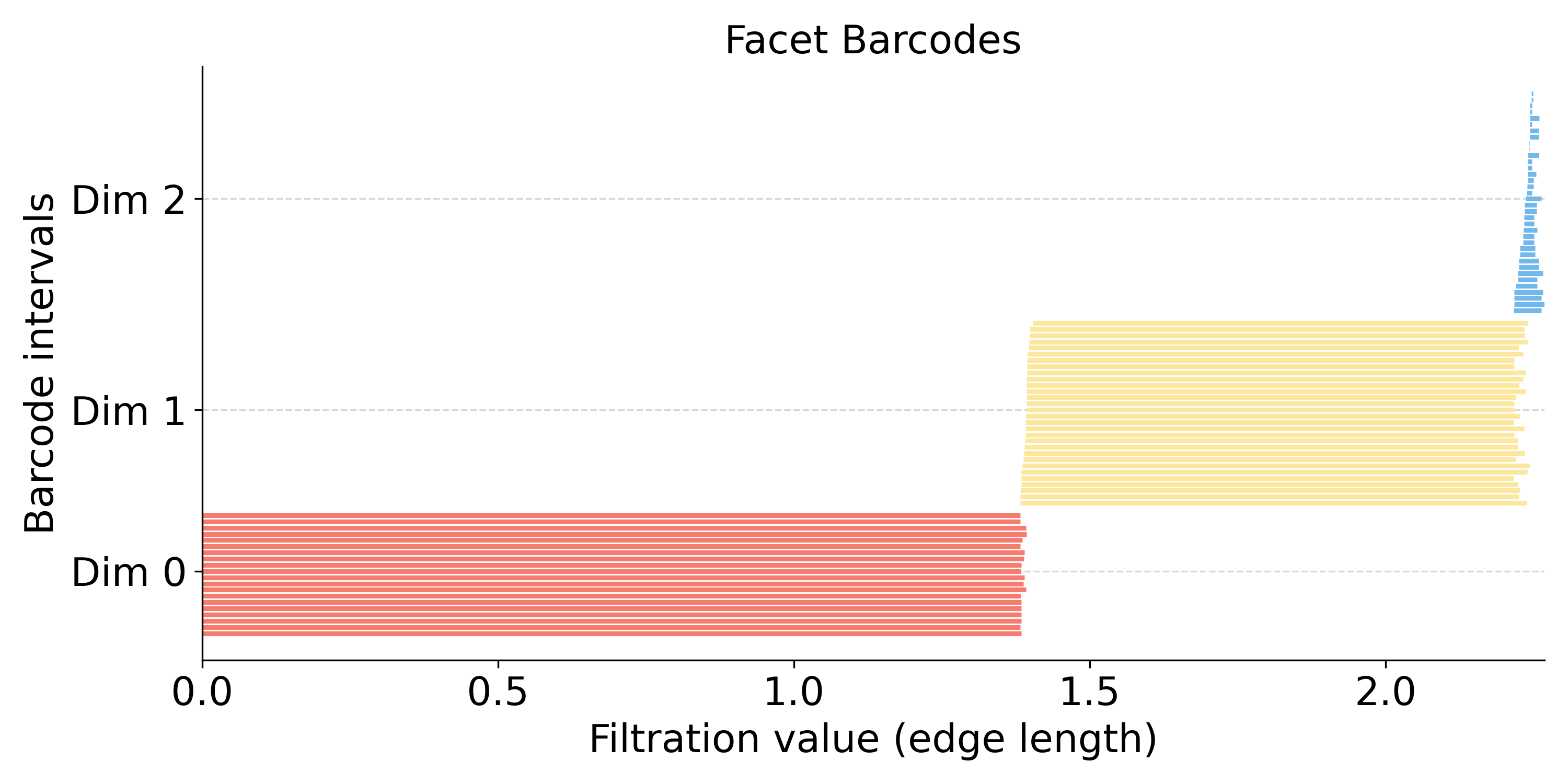}
    \subcaption*{(a) Facet persistence.}
  \end{subfigure}\hfill
  \begin{subfigure}{0.48\linewidth}
    \centering
    \includegraphics[width=\linewidth,height=5cm]{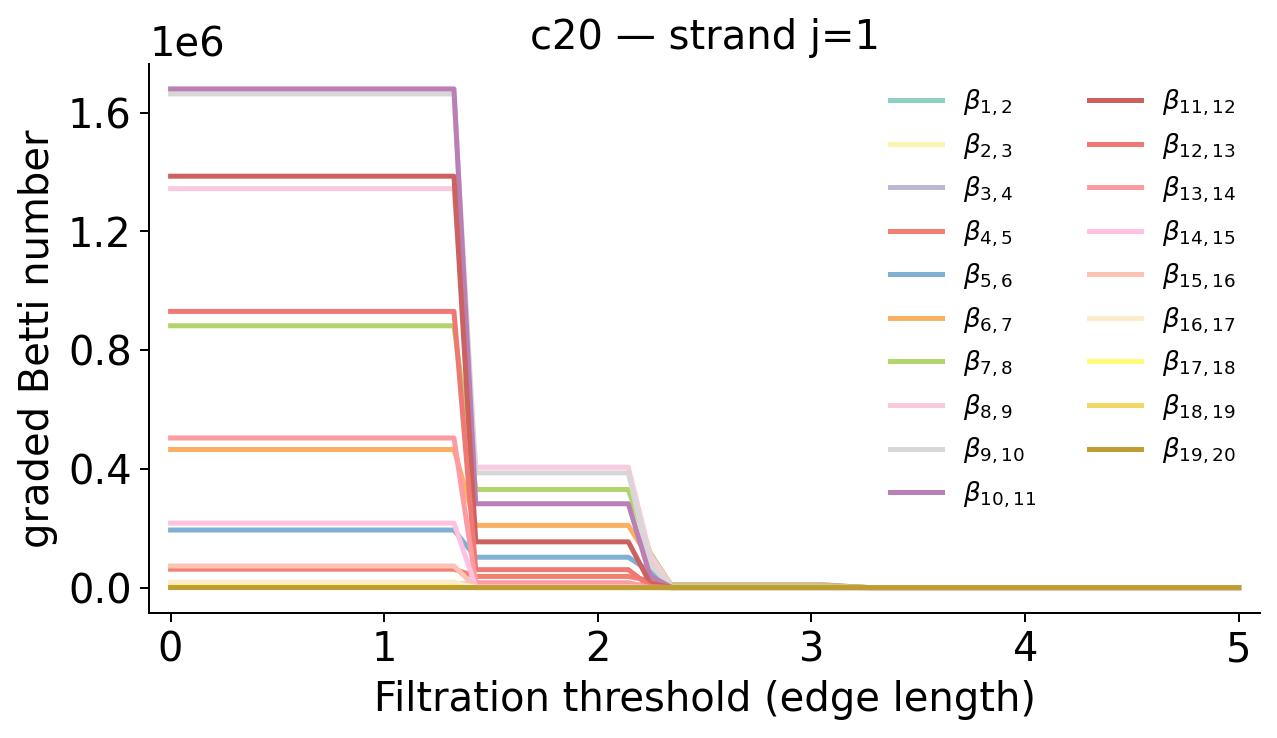}
    \subcaption*{(b) Graded Betti curves, dim$=0$ (strand $j=1$).}
  \end{subfigure}

  \vspace{0.5em}

  \begin{subfigure}{0.49\linewidth}
    \centering
    \includegraphics[width=\linewidth,height=5cm]{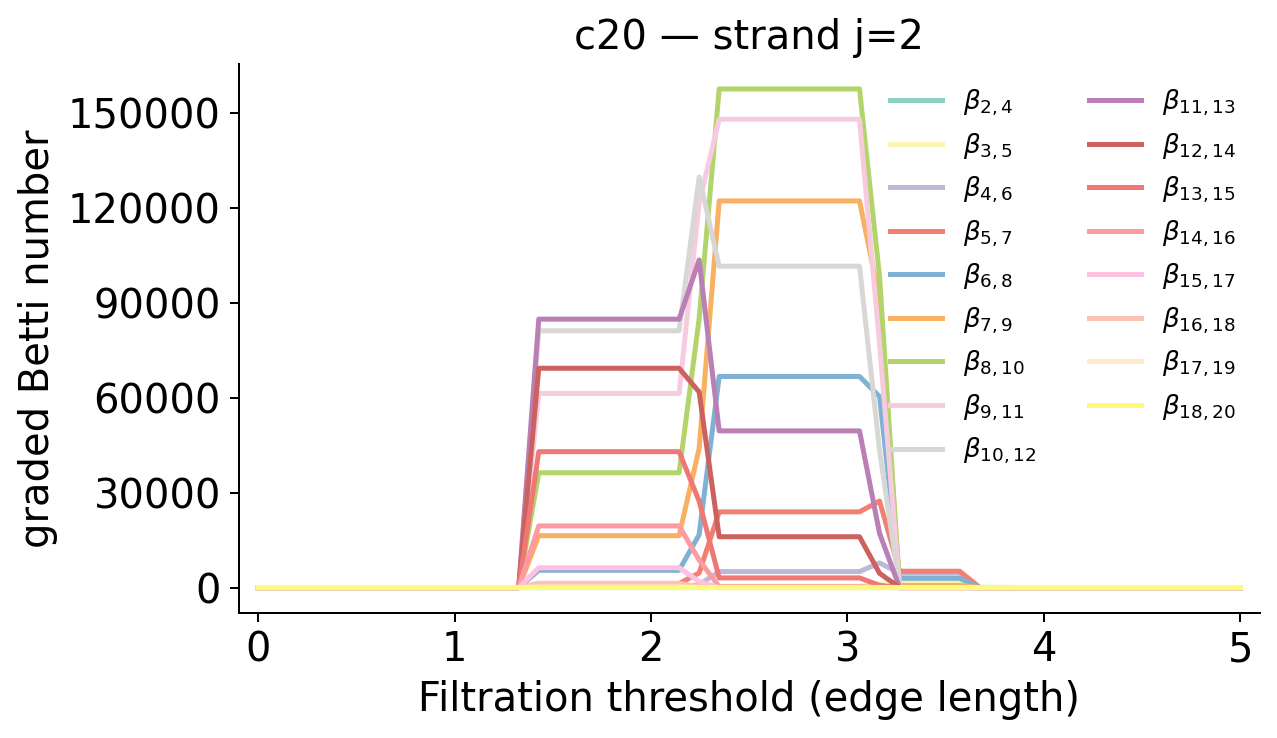}
    \subcaption*{(c) Graded Betti curves, dim$=1$ (strand $j=2$).}
  \end{subfigure}\hfill
  \begin{subfigure}{0.49\linewidth}
    \centering
    \includegraphics[width=\linewidth,height=5cm]{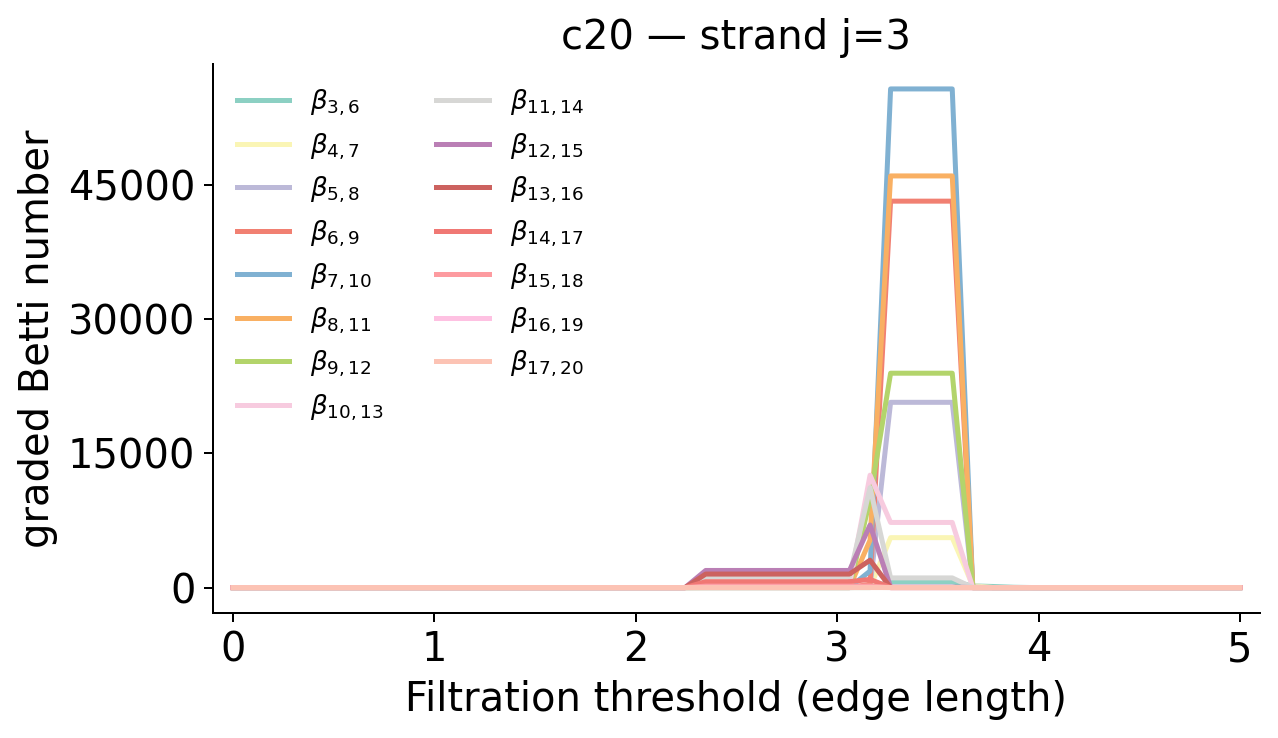}
    \subcaption*{(d) Graded Betti curves, dim$=2$ (strand $j=3$).}
  \end{subfigure}

  \vspace{0.5em}

  \begin{subfigure}{0.48\linewidth}
    \centering
    \includegraphics[width=\linewidth, height=5.5cm]{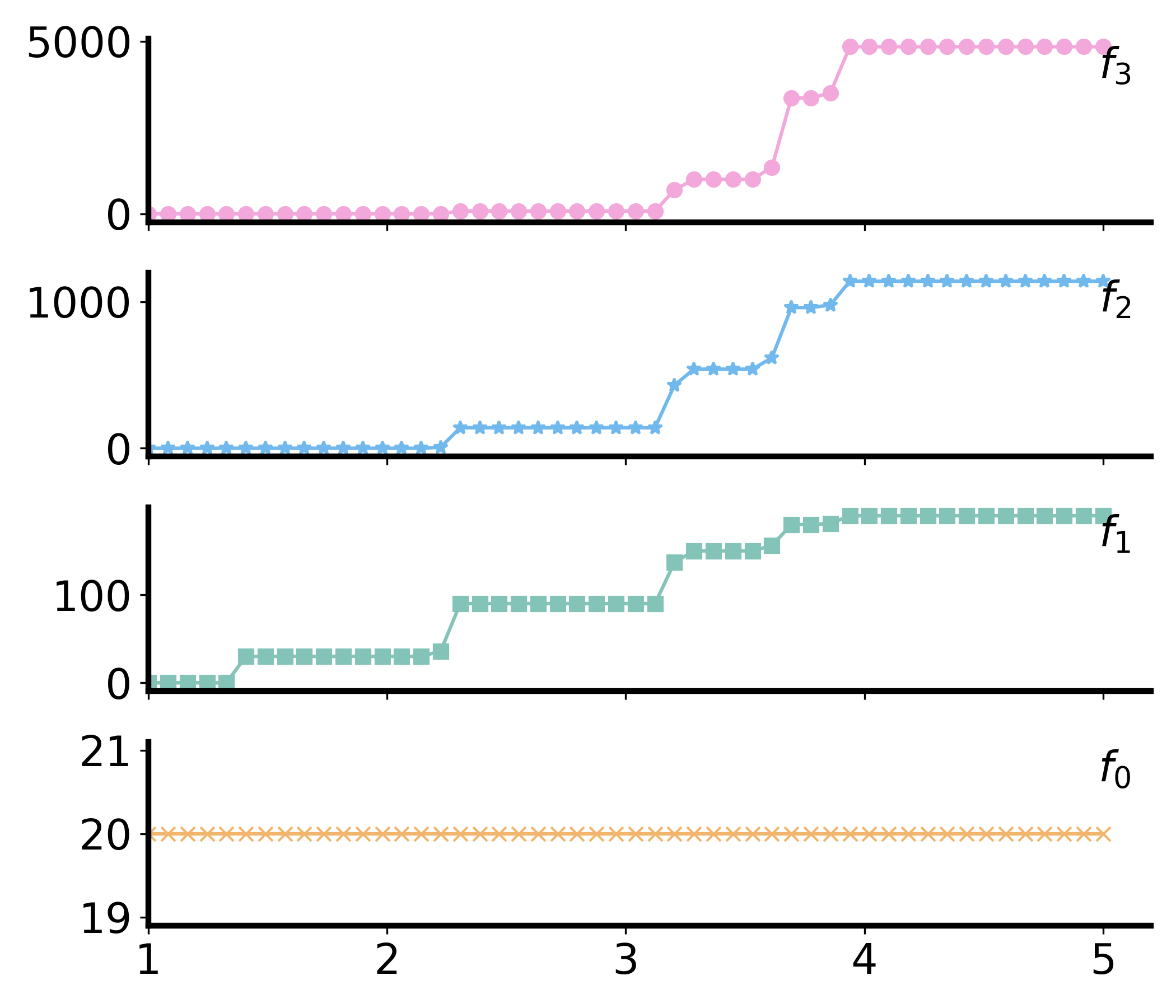}
    \subcaption*{(e) $f$-vector curves.}
  \end{subfigure}\hfill
  \begin{subfigure}{0.48\linewidth}
    \centering
    \includegraphics[width=\linewidth, height=5.5cm]{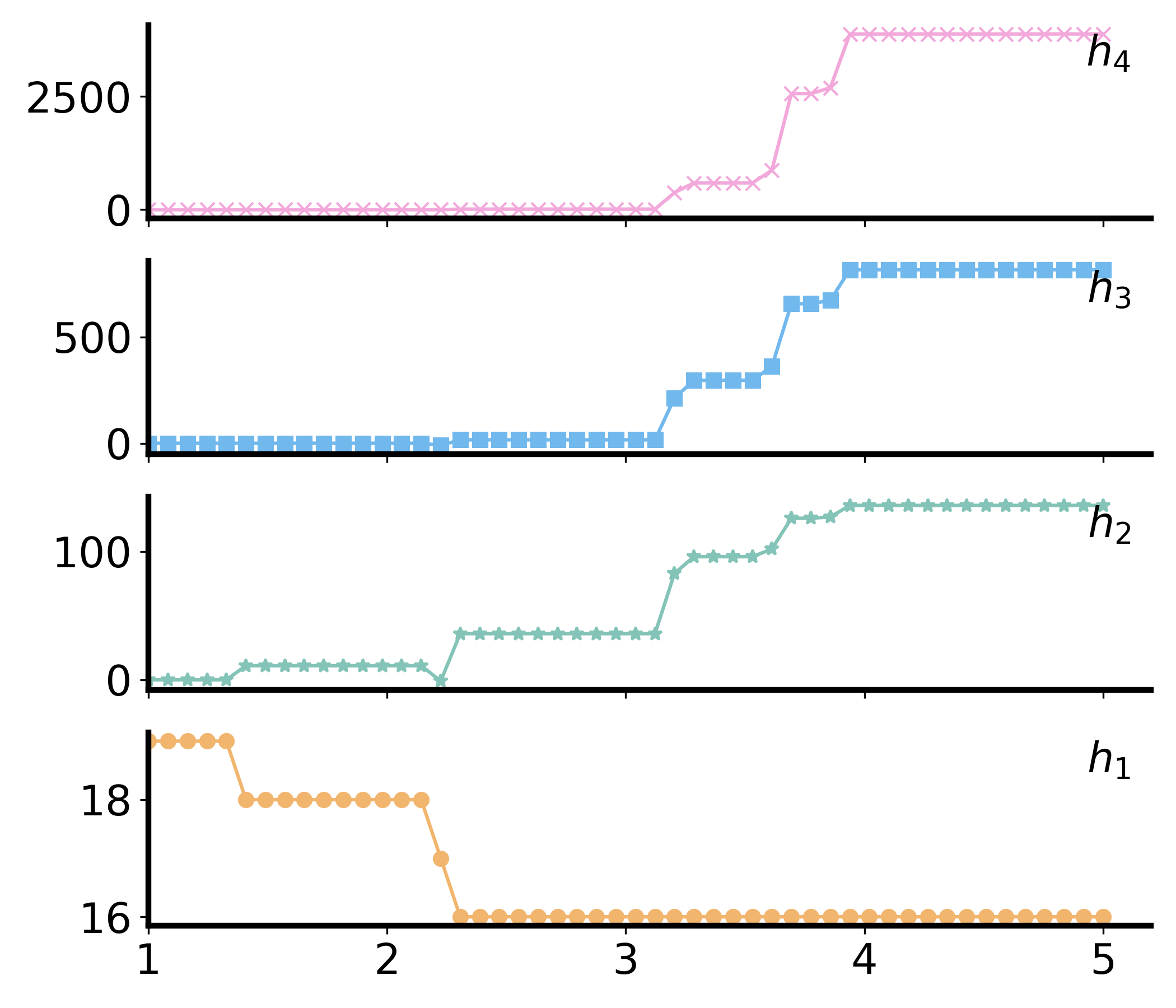}
    \subcaption*{(f) $h$-vector curves.}
  \end{subfigure}

  \caption{Persistent commutative algebra  analysis of the C\textsubscript{20} molecule using a Rips complex-based filtration process.}
  \label{fig:pca_c20}
\end{figure}
In this subsection, we use the fullerene \(\mathrm{C}_{20}\) to illustrate the interpretability and representational power of PH, PL, and PCA. The \(\mathrm{C}_{20}\) cage comprises twenty carbon atoms arranged into twelve pentagonal faces, each formed by five atoms. The Vietoris–Rips filtration on the twenty vertices is governed by a sequence of critical scales, shown in the top panel of Fig.~\ref{fig:phl_c20}, which determine when specific simplices enter. At the nearest–neighbor bond length, \(\varepsilon \approx 1.4~\text{\AA}\), edges along pentagonal rims first appear. At the next scale, \(\varepsilon \approx 2.0\)–\(2.3~\text{\AA}\), diagonals within pentagons create filled triangles, so each pentagonal face becomes triangulated. As \(\varepsilon\) increases further, \(\varepsilon \approx 3.2~\text{\AA}\), large batches of edges are added simultaneously, introducing many new triangles and higher–dimensional shells. In the range \(\varepsilon \approx 3.7\)–\(3.9~\text{\AA}\), tetrahedra begin to fill these triangular shells. Finally, at \(\varepsilon \approx 4.0~\text{\AA}\), which matches the longest atomic separation, the complex becomes the full simplex on twenty vertices.

For the C\(_{20}\) fullerene, the harmonic spectra of PL reproduce PH in the bottom panel of Fig.~\ref{fig:phl_c20}, as seen in the Betti–barcode summary. The \(\beta_{0}\) bars record twenty initial components, one per atom, which merge to a single component. The \(\beta_{1}\) bars capture eleven independent loops, reflecting the twelve pentagons with one topological dependence. The \(\beta_{2}\) barcode shows a single class, representing the global void enclosed by the \(\mathrm{C}_{20}\) cage. The nonharmonic spectra of PL then enrich this picture through the smallest positive eigenvalues. At \(\varepsilon \approx 1.4~\text{\AA}\), \(\lambda_{0}\) and \(\lambda_{1}\) turn on with the emergence of bonds, while \(\lambda_{2}\) is absent since no triangles are present. At \(\varepsilon \approx 2.3~\text{\AA}\), triangles appear, \(\lambda_{2}\) turns on, and both \(\lambda_{0}\) and \(\lambda_{1}\) increase slightly. Notably, one–dimensional topological activity ceases by \(\varepsilon \approx 2.3~\text{\AA}\) in PH, yet the nonharmonic spectra of PL continue to evolve in several steps, reflecting geometric reinforcement beyond topological change. At \(\varepsilon \approx 3.2~\text{\AA}\), numerous edges are added, producing a jump in \(\lambda_{0}\) and \(\lambda_{1}\), and a slight dip in \(\lambda_{2}\), because the number of triangles grows faster than the number of tetrahedra. Subsequent rises of all three curves at \(\varepsilon \approx 3.7\) and \(3.9~\text{\AA}\) reflect the addition of further edges and faces, while all three spectra flatten at the constant value \(20\) at \(4.0~\text{\AA}\), consistent with the complex becoming the complete simplex.

We adopt a  PCA perspective on the \(\mathrm{C}_{20}\) cage, interpreting its features through facet persistence, \(f\)– and \(h\)–vectors, and graded Betti numbers. Figure~\ref{fig:pca_c20}(a) shows the persistent facet ideals in dimensions \(0\), \(1\), and \(2\). In dimension \(0\), facet bars persist up to the bond length, \(\varepsilon \approx 1.4\,\text{\AA}\), where every vertex acquires an incident edge. These edges remain until the diagonals within the pentagons shorten enough to create many triangles, causing the yellow facet bars of dimension \(1\) to disappear and the blue facet bars of dimension \(2\) to appear around \(\varepsilon \approx 2.2\)–\(2.3\,\text{\AA}\). In dimension \(2\), triangle facet bars are short-lived, since groups of four vertices quickly come into range as \(\varepsilon\) grows, tetrahedra appear, and their triangular faces are absorbed. The \(f\)– and \(h\)–curves in Fig.~\ref{fig:pca_c20}(e) and Fig.~\ref{fig:pca_c20}(f) count simplices present at each \(\varepsilon\). One has \(f_{0}(\varepsilon)= 20\). The edge count \(f_{1}\) makes its first large jump at \(\varepsilon \approx 1.4\,\text{\AA}\), when bond length edges form on the pentagonal faces, rises again near \(\varepsilon \approx 2.2\)–\(2.3\,\text{\AA}\), as many new edges arrive together with the first triangles, and increases around \(\varepsilon \approx 3.2\,\text{\AA}\), when a dense batch of longer edges appears. The triangle count \(f_{2}\) turns on at \(\varepsilon \approx 2.2\)–\(2.3\,\text{\AA}\), and grows through \(3.2\,\text{\AA}\). The tetrahedron count \(f_{3}\) begins near \(3.2\,\text{\AA}\), increases sharply in the range \(\varepsilon \approx 3.7\)–\(3.9\,\text{\AA}\), and reaches its final value at \(\varepsilon = 4.0\,\text{\AA}\). At this last threshold, the complex becomes the full simplex on twenty vertices, and all counts stabilize.

We organize the graded Betti numbers \(\beta_{i,i+j}\) for \(\mathrm{C}_{20}\) by strands \(j=1,2,3\), as shown in Fig.~\ref{fig:pca_c20}(b)–(d). On strand \(j=1\), the curves start at large combinatorial plateaus and drop sharply at \(\varepsilon \approx 1.4\,\text{\AA}\), when many \((i+j)\)–vertex subcomplexes first acquire the short edges that connect them. Additional step–downs occur near \(\varepsilon \approx 2.2\)–\(2.3\,\text{\AA}\) and around \(3.2\,\text{\AA}\), reflecting successive gains in connectivity as diagonals within pentagons and then longer edges across the cage appear.  On strand \(j=2\), new edges create one–dimensional cycles in many induced subcomplexes at \(\varepsilon \approx 2.0\)–\(2.3\,\text{\AA}\), so graded Betti numbers with mid–size vertex such as \(\beta_{8,10}\) and \(\beta_{9,11}\) rise. At the same radii, other pairs, for example \(\beta_{12,14}\), \(\beta_{13,15}\), and \(\beta_{14,16}\), drop because triangles already form inside those larger \((i+j)\)–subcomplexes, filling the loops as soon as they appear. After the dense edge layer near \(3.2\,\text{\AA}\), most \(j=2\) curves collapse toward zero. A small tail persists until \(\varepsilon \approx 3.4\)–\(3.5\,\text{\AA}\), which implies that a minority of induced subcomplexes require a slightly larger threshold before the final diagonals fall in range and triangles cap the remaining loops. On strand \(j=3\), weak signals emerge once triangles are common, near \(2.3\,\text{\AA}\), followed by a pronounced surge across moderate \((i+j)\) at \(\varepsilon \approx 3.1\)–\(3.3\,\text{\AA}\) when many triangular shells form. Two bands of "death times” are visible. Some curves drop quickly at \(\varepsilon \approx 3.2\)–\(3.4\,\text{\AA}\), which corresponds to compact subcomplexes whose pairwise edges already lie below the threshold, so tetrahedra are immediately available and the shells become boundaries at once. Other curves persist and drop later, in the \(\varepsilon \approx 3.7\)–\(3.9\,\text{\AA}\) range, because they involve larger or more spread–out \((i+j)\)–subsets and require a larger radius before the last closing edges appear and tetrahedra can form. 
At \(\varepsilon = 4.0\,\text{\AA}\) the complex is the full simplex on twenty vertices, and all graded Betti numbers vanish. We note that the raw magnitudes of these graded Betti numbers are large on \(\mathrm{C}_{20}\) data. The graded Betti number on the strand \(j=1\) reaches about \( 10^{6}\), \(j=2\) reaches about \( 10^{5}\), and \(j=3\) reaches about \(10^{4}\), which is a consequence of combinatorial growth in the number of induced \((i+j)\)–subcomplexes.

\subsection{Protein Structures: Molecular Fingerprints}
\begin{figure}
    \centering
     \includegraphics[width=0.8\linewidth, height=8cm]{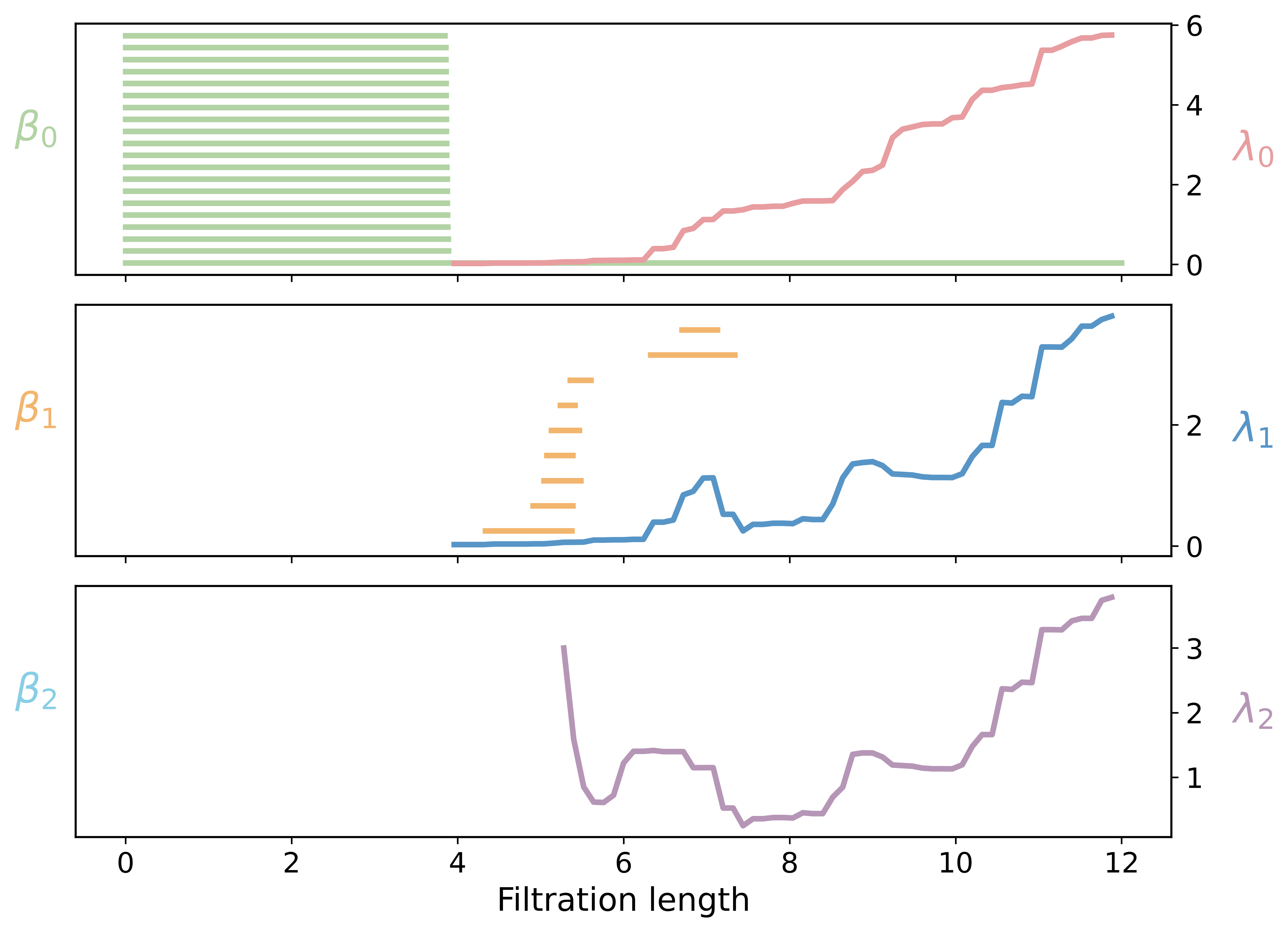}
     \caption{{Illustration of PH and PL on Protein structure (PDBID:1L2Y). Barcodes and spectra summarize the evolution of $0$-, $1$-, and $2$-dimensional features across the filtration.}}
         \label{fig:1L2Y}
\end{figure}

\begin{figure}[htbp!]
  \centering

  \begin{subfigure}{0.42\linewidth}
    \centering
    \includegraphics[width=\linewidth,height=5cm]{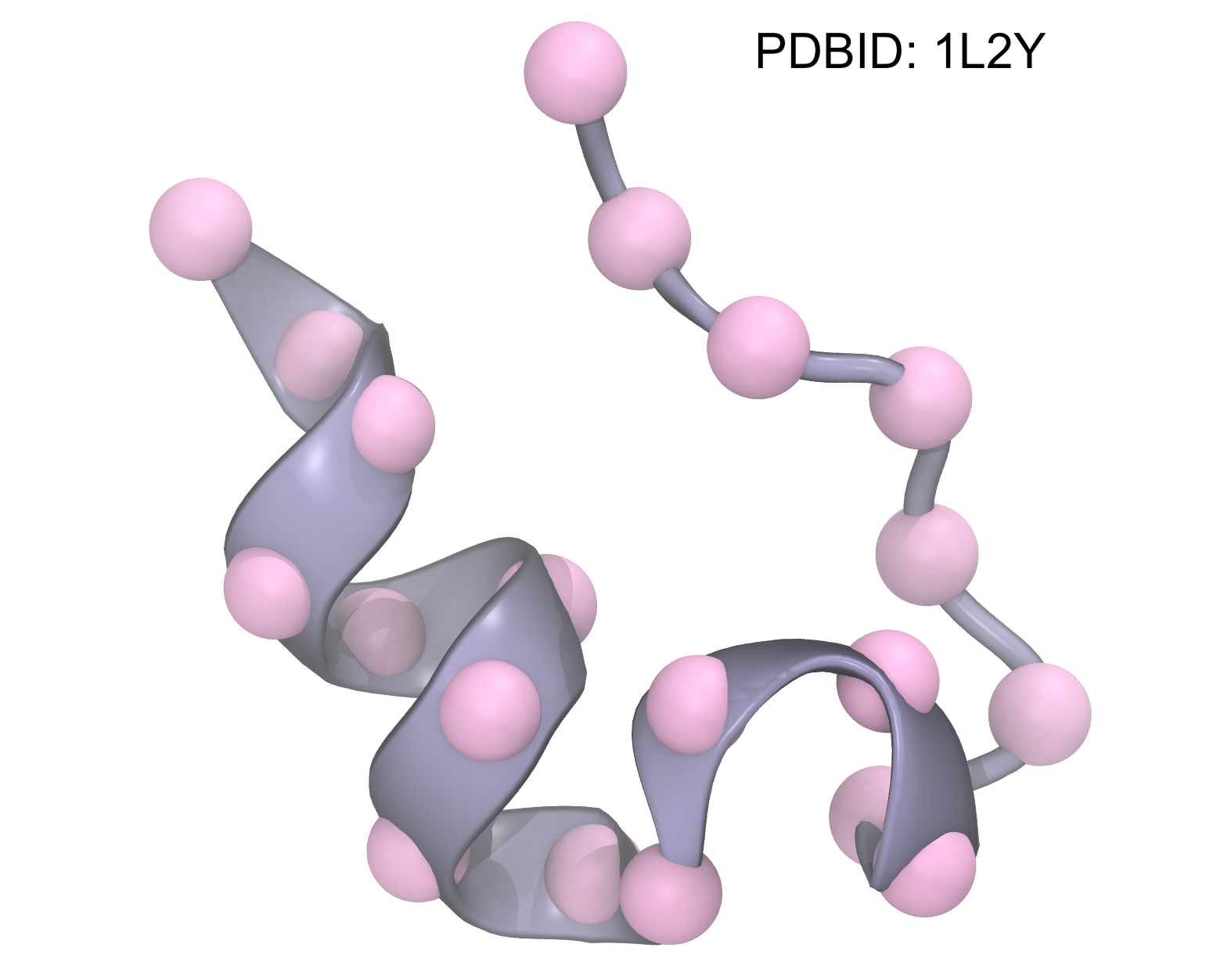}
    \subcaption*{(a) Protein structure.}
  \end{subfigure}\hfill
  \begin{subfigure}{0.54\linewidth}
    \centering
    \includegraphics[width=\linewidth,height=5cm]{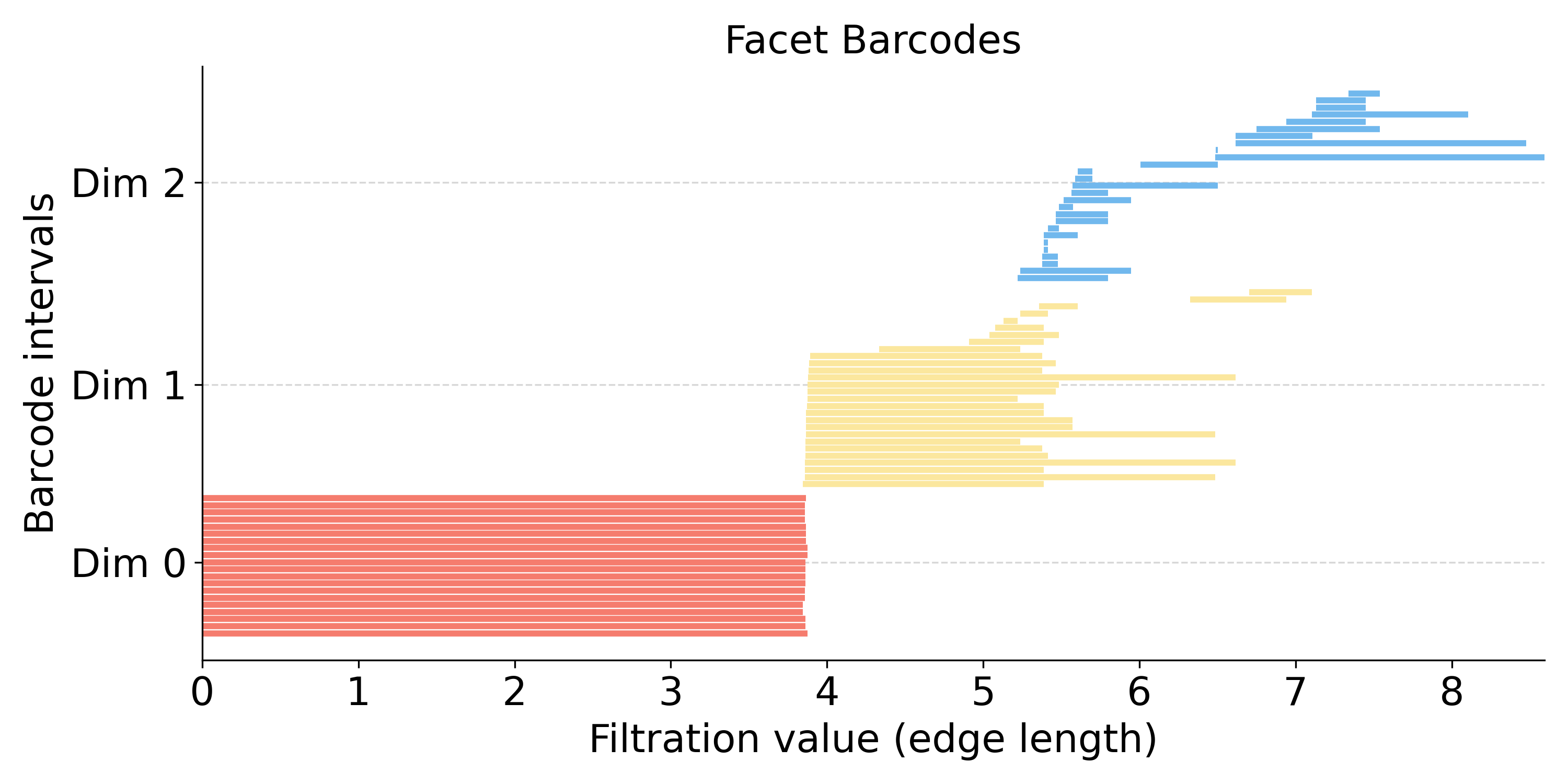}
    \subcaption*{(b) Facet persistence}
  \end{subfigure}

  \begin{subfigure}{0.49\linewidth}
    \centering
    \includegraphics[width=\linewidth, height=5cm]{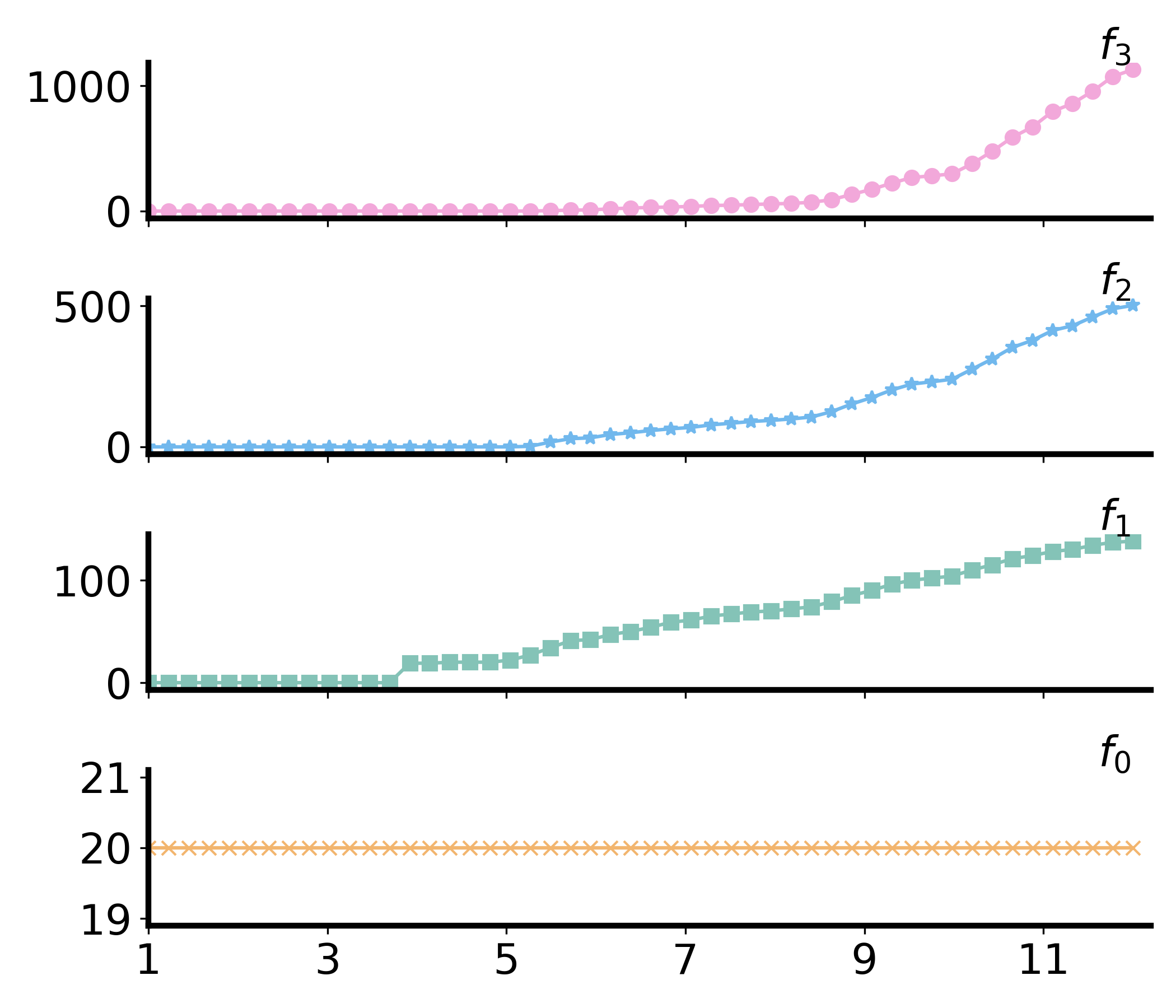}
    \subcaption*{(c) $f$-vector curves}
  \end{subfigure}\hfill
  \begin{subfigure}{0.49\linewidth}
    \centering
    \includegraphics[width=\linewidth, height=5cm]{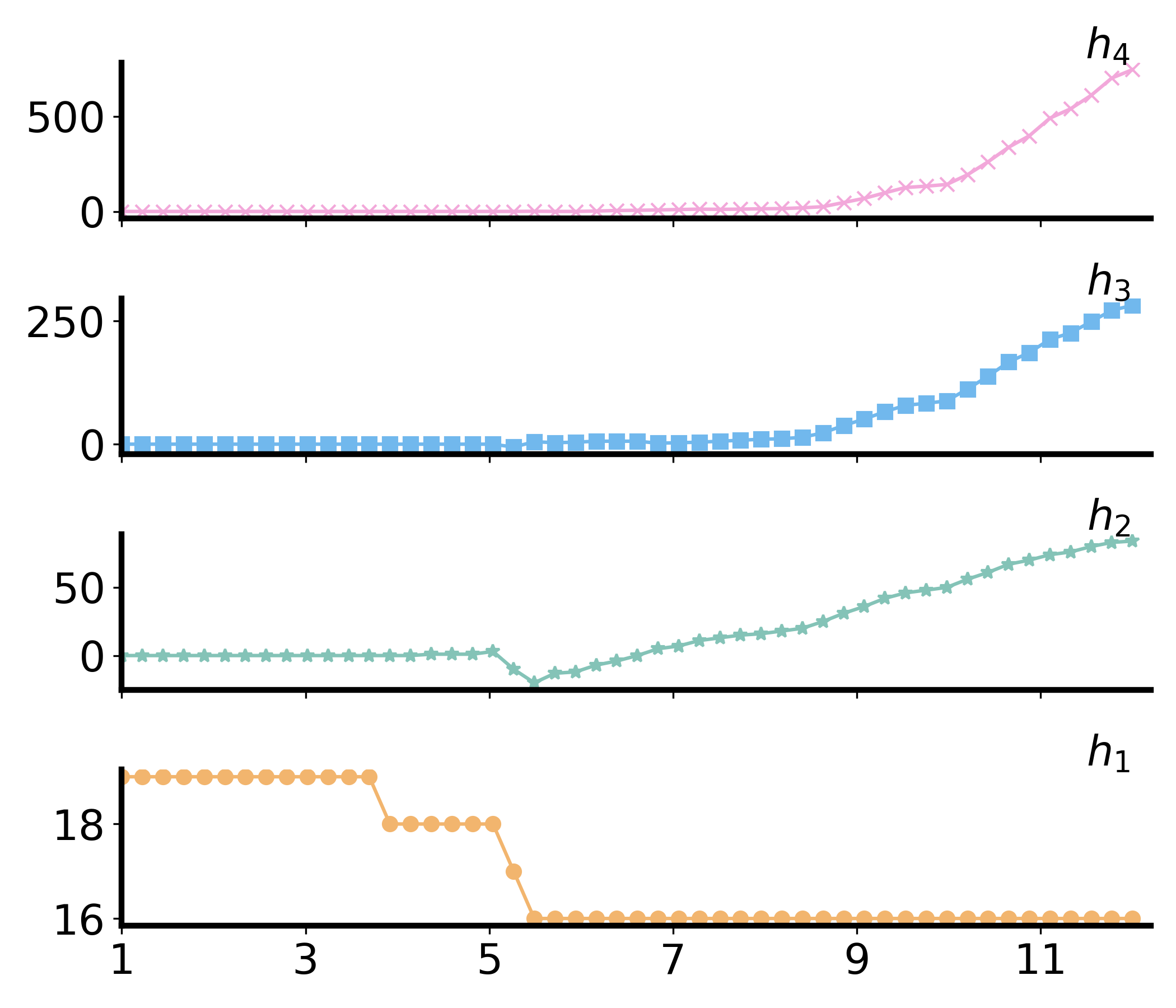}
    \subcaption*{(d) $h$-vector curves}
  \end{subfigure}

\begin{subfigure}{0.49\linewidth}
    \centering
    \includegraphics[width=\linewidth,height=5cm]{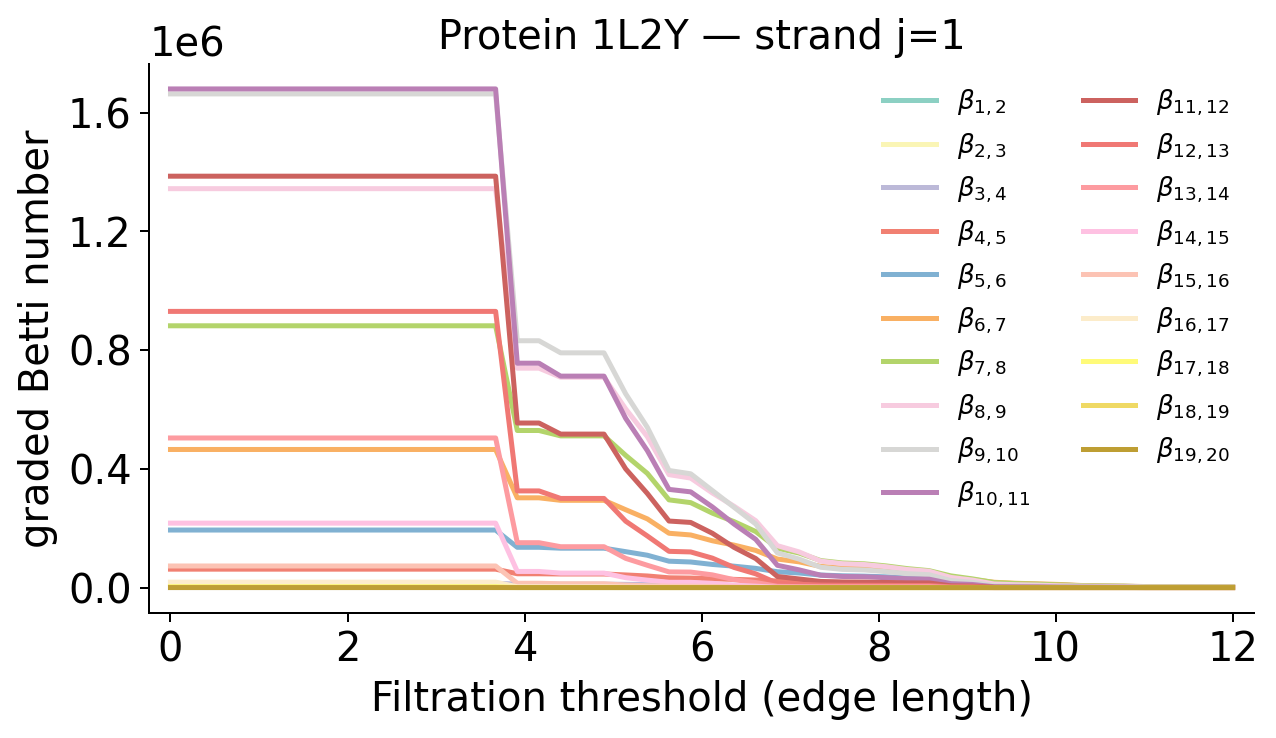}
    \subcaption*{(e) Graded Betti, dim$=0$ ($j=1$)}
  \end{subfigure}\hfill
  \begin{subfigure}{0.49\linewidth}
    \centering
    \includegraphics[width=\linewidth,height=5cm]{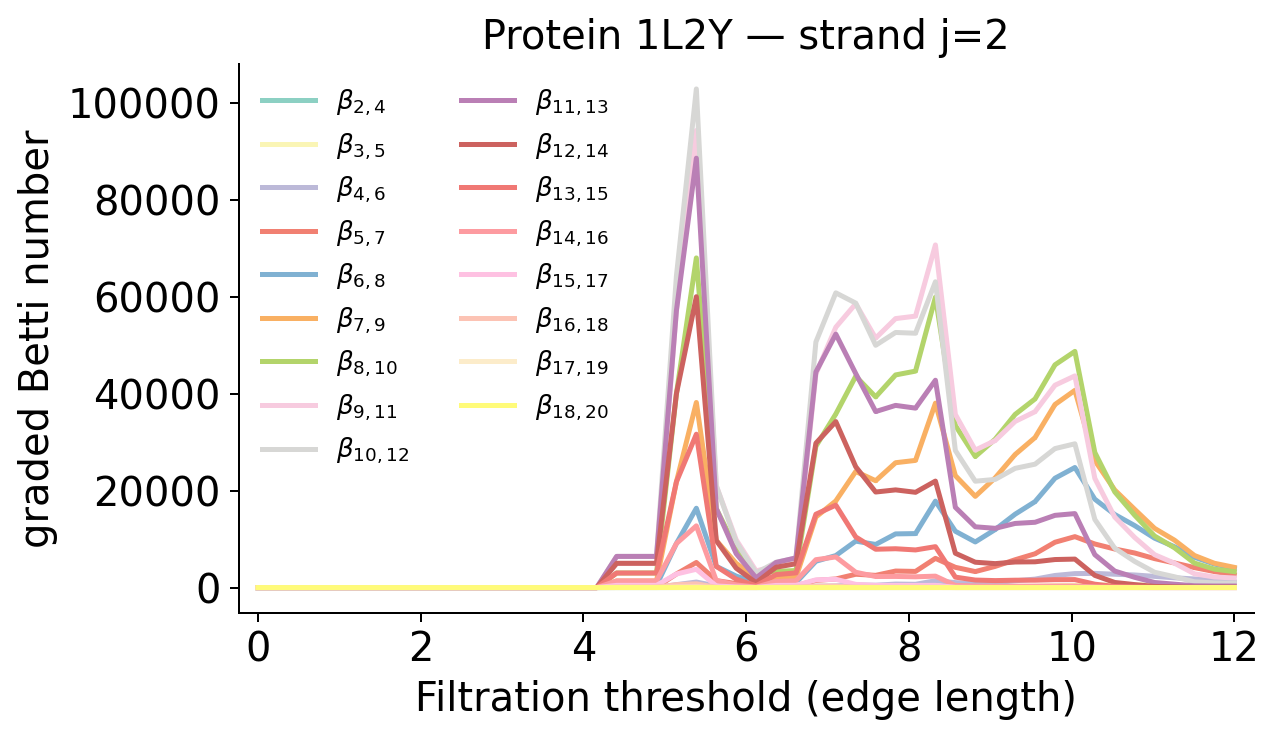}
    \subcaption*{(f) Graded Betti, dim$=1$ ($j=2$)}
  \end{subfigure}


  \begin{subfigure}{0.49\linewidth}
    \centering
    \includegraphics[width=\linewidth,height=5cm]{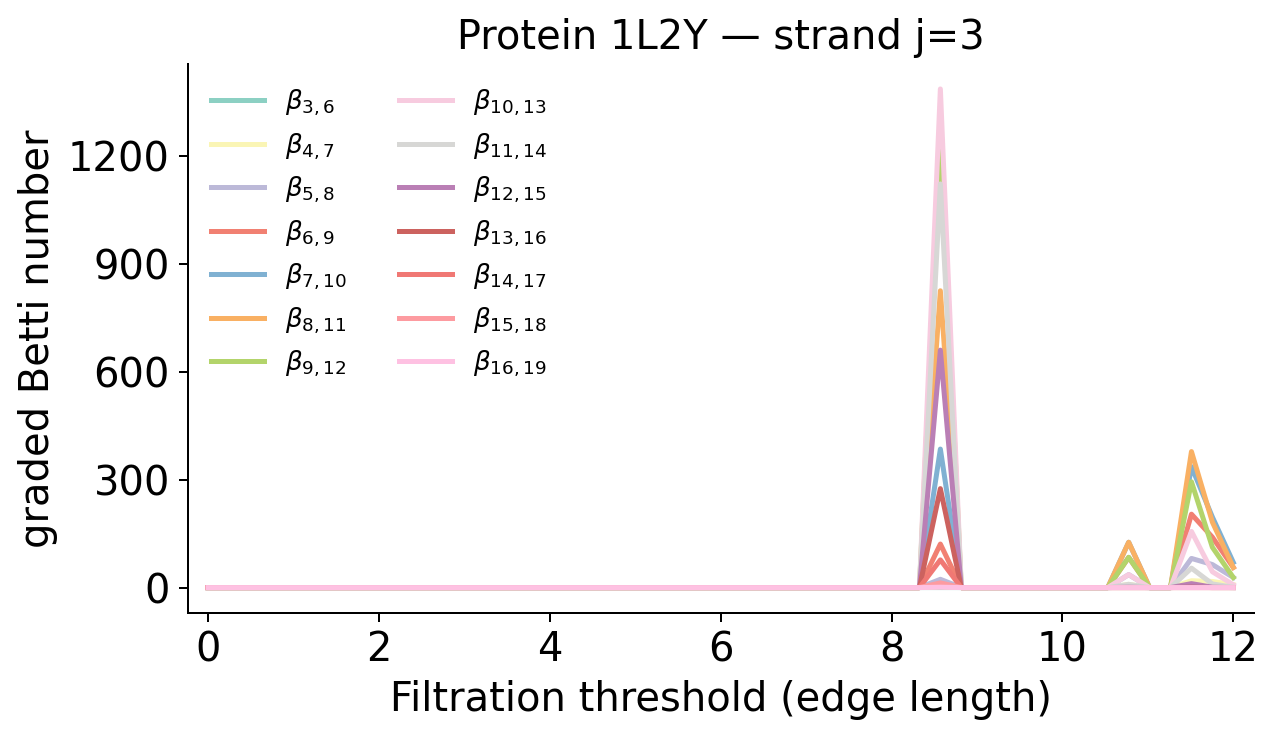}
    \subcaption*{(g) Graded Betti, dim$=2$ ($j=3$)}
  \end{subfigure}

  \caption{Persistent commutative algebra  analysis of protein 1L2Y using a Vietoris--Rips filtration}
  \label{fig:1L2Y2}
\end{figure}

\begin{figure}
    \centering
     \includegraphics[width=0.8\linewidth, height=8cm]{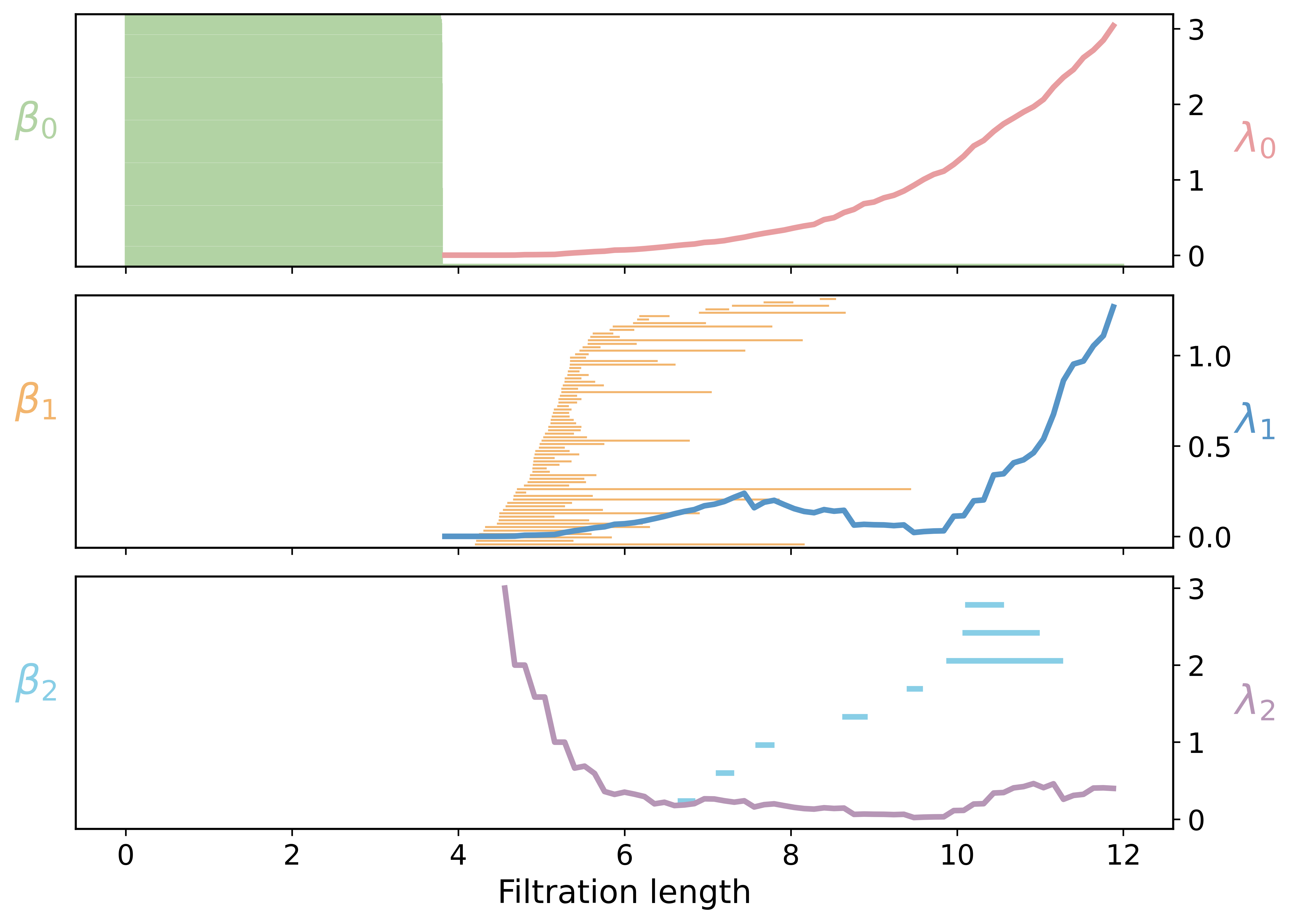}
     \caption{{Illustration of PH and PL on Protein structure (PDBID: 2LYZ). Barcodes and spectra summarize the evolution of $0$-, $1$-, and $2$-dimensional features across the filtration.}}
         \label{fig:2LYZ}
\end{figure}

Proteins are composed of one or more coiled peptide chains that fold into highly complex three-dimensional structures, with each protein adopting a unique conformation. Understanding these structures requires mathematical tools capable of capturing both geometry and topology. Our objective is to employ three types of molecular fingerprints—PH, PL, and PCA—to reveal the intrinsic topological and geometric characteristics of proteins and to identify distinctive structural fingerprints.

\begin{figure}[htbp!]
  \centering

  \begin{subfigure}{0.42\linewidth}
    \centering
    \includegraphics[width=\linewidth,height=5.5cm]{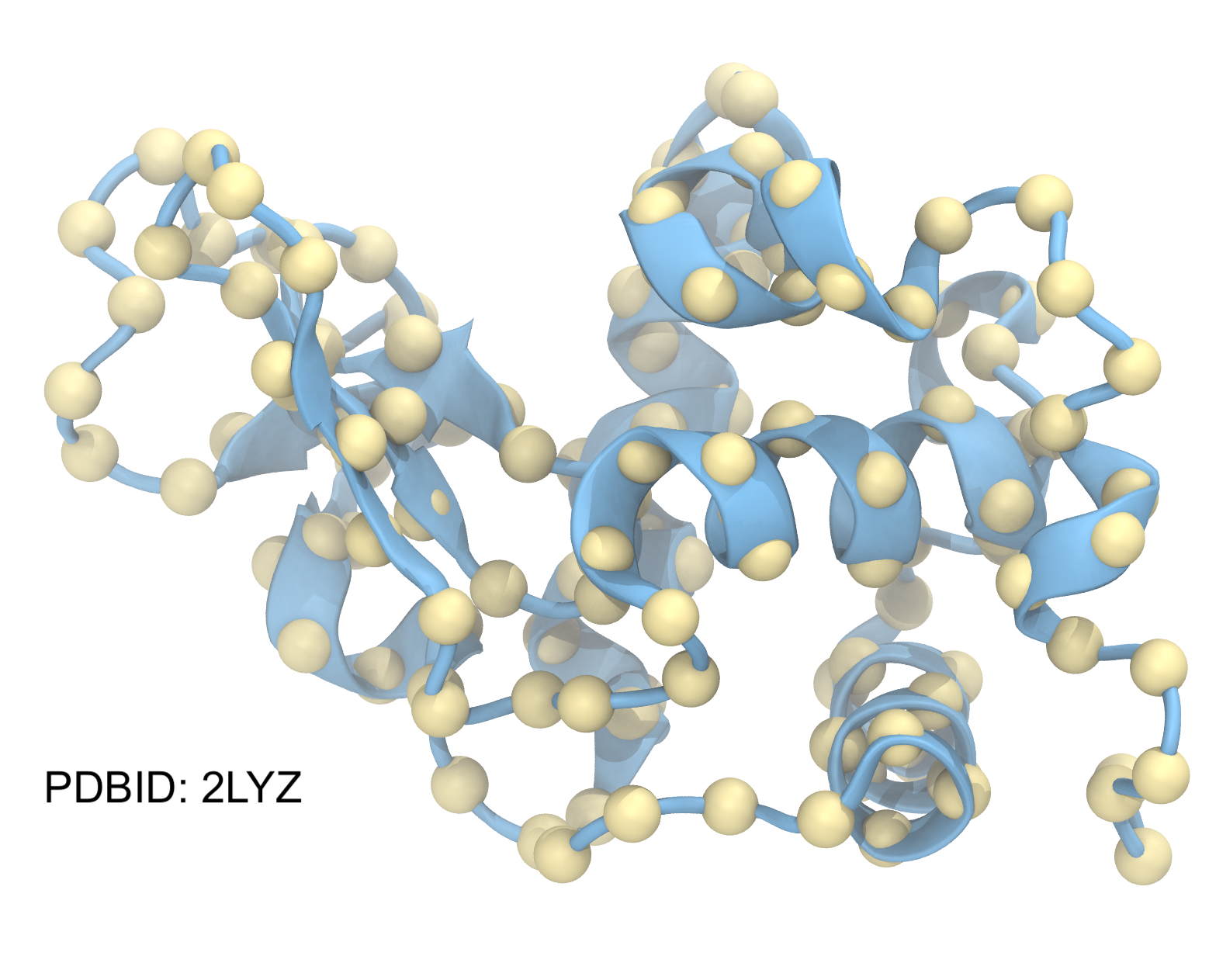}
    \subcaption*{(a) Protein 2LYZ}
  \end{subfigure}\hfill
  \begin{subfigure}{0.54\linewidth}
    \centering
    \includegraphics[width=\linewidth,height=5cm]{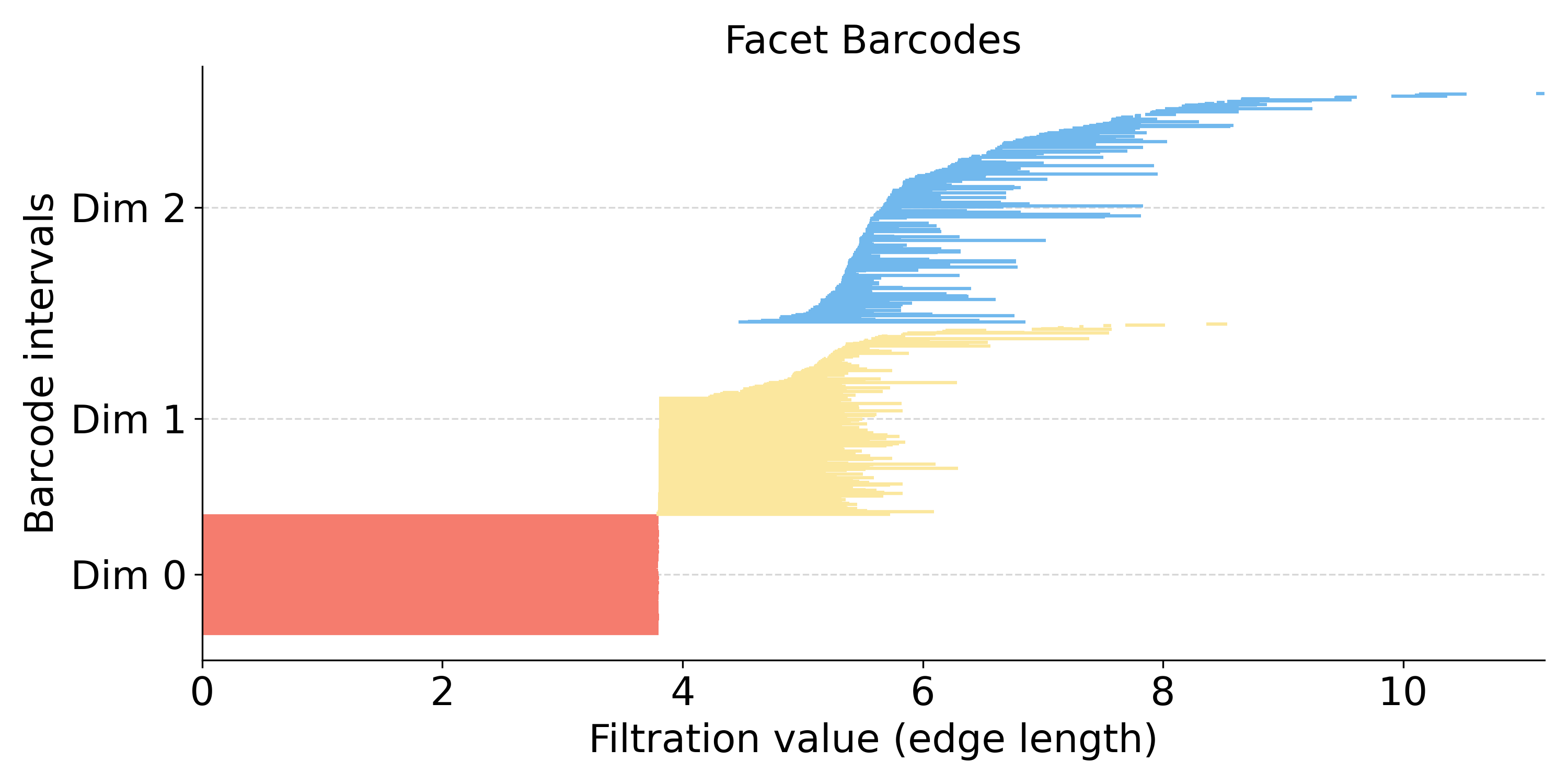}
    \subcaption*{(b) Facet persistence}
  \end{subfigure}

  \vspace{0.6em}

  \begin{subfigure}{0.49\linewidth}
    \centering
    \includegraphics[width=\linewidth]{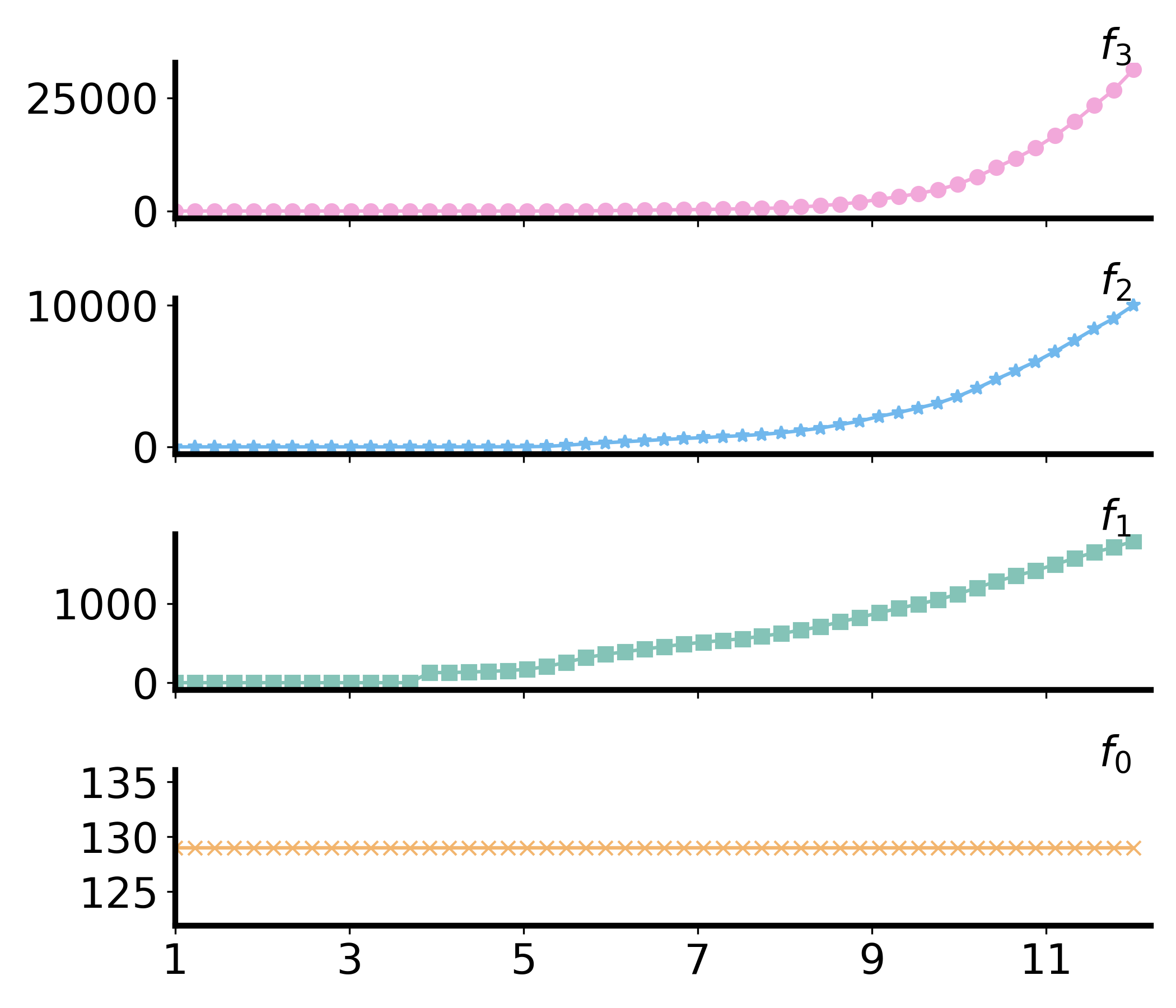}
    \subcaption*{(c) $f$-vector curves}
  \end{subfigure}\hfill
  \begin{subfigure}{0.49\linewidth}
    \centering
    \includegraphics[width=\linewidth]{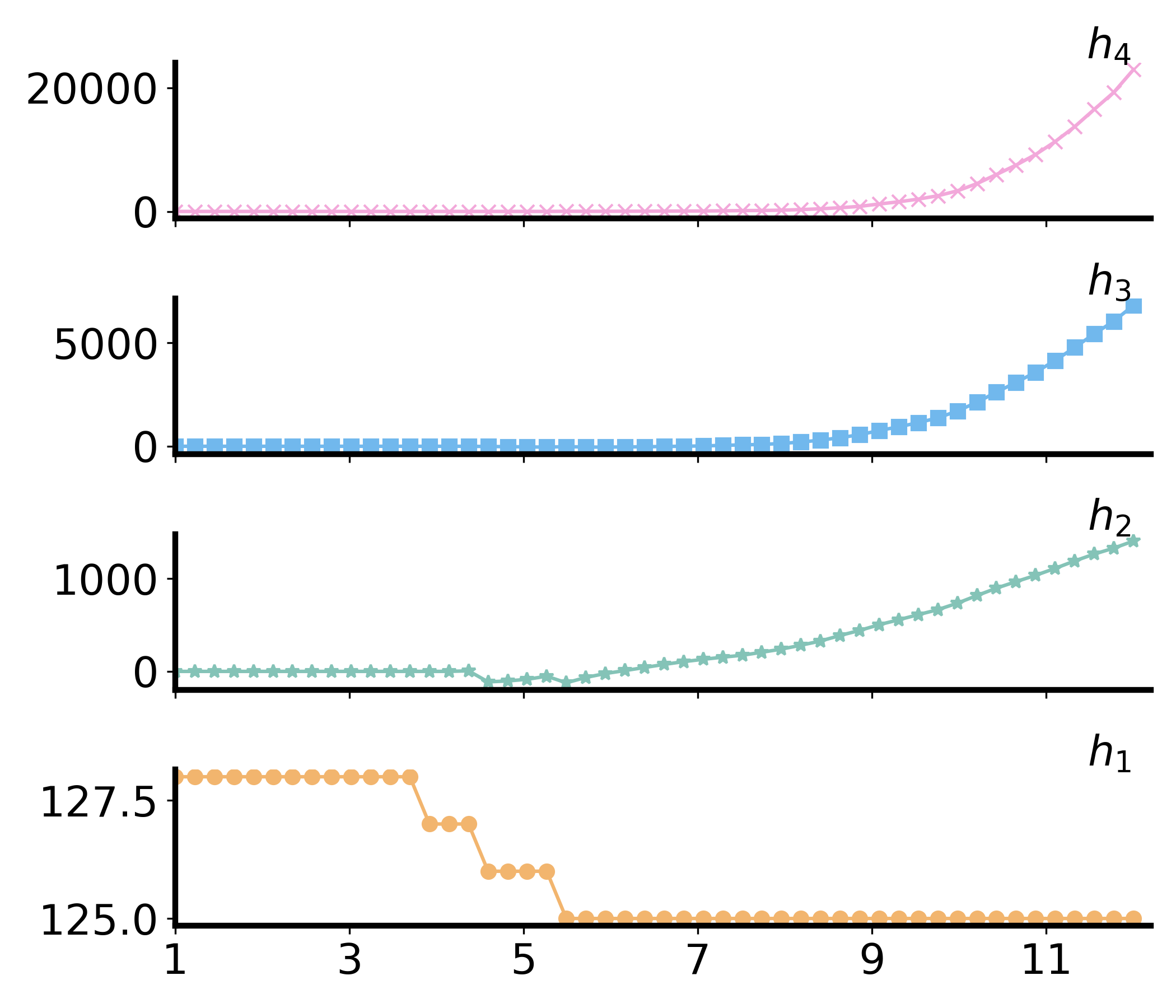}
    \subcaption*{(d) $h$-vector curves}
  \end{subfigure}

  \caption{Persistent commutative algebra  analysis of protein 2LYZ using a Rips complex-based filtration process.}
  \label{fig:2LYZ2}
\end{figure}

Protein structural data are publicly available through the Protein Data Bank (PDB), providing a rich source for such analyses. To illustrate, we focus on an $\alpha$-helical protein with PDB ID: 1L2Y, a de novo design consisting of a single chain of 20 residues. For clarity in extracting structural patterns, the protein is simplified by mapping each amino acid to its C$\alpha$ atom, which serves as a backbone representation of the helix. As shown in Fig.~\ref{fig:2LYZ2}(a), the pink spheres mark the C\(\alpha\) positions of Protein 1L2Y, and the translucent ribbon traces the helical segment.

We construct a Vietoris–Rips filtration on the set of 20 C\(\alpha\) atoms of protein 1L2Y, and summarize PH and PL features for three dimensions in Fig.~\ref{fig:2LYZ}. In dimension \(0\), nineteen of \(\beta_{0}\) bars terminate by \(3.8\,\text{\AA}\), consistent with the average C\(\alpha\)–C\(\alpha\) bond length. At this scale, all local connections have formed, and the chain merges into a single connected component, leaving only one \(\beta_{0}\) bar that represents the global connectivity of the structure. The spectrum \(\lambda_{0}\) turns on at this merge, then rises in steps that reflect the progressive addition of cross–turn edges within the helix. In dimension \(1\), \(\beta_{1}\) bars are born near \(5\,\text{\AA}\), when nonsequential helical neighbors begin to connect and close short ring–like loops around the helix axis. As triangles appear and fill these loops, the \(\beta_{1}\) bars die by \(8\,\text{\AA}\). The spectrum \(\lambda_{1}\) turns on at \(4\,\text{\AA}\), shows a small bump around \(6\,\text{\AA}\) when early triangles stabilize a subset of loops, dips near \(7\,\text{\AA}\) as longer edges spanning multiple turns occasionally including residues from the tail transiently appear. At \(8\,\text{\AA}\), no additional one–dimensional topological events occur, yet the spectrum \(\lambda_{1}\) rises again as triangles and higher–dimensional faces accumulate, reflecting increased geometric stiffening of the complex.
In dimension \(2\), no \(\beta_{2}\) bars appear. Nevertheless, \(\lambda_{2}\) records geometric stabilization of higher–order structure. An early shoulder of \(\lambda_{2}\) at \(5\)–\(6\,\text{\AA}\) corresponds to the first triangles among cross–turn neighbors within the helix. A broad minimum of \(\lambda_{2}\) near \(7\)–\(8\,\text{\AA}\) arises when triangles accumulate faster than tetrahedra. A subsequent rise of \(\lambda_{2}\) around \(10\,\text{\AA}\) indicates that tetrahedra become abundant across multiple helical turns and between the helix and the tail.

The PCA analysis of protein 1L2Y is shown in Fig.~\ref{fig:2LYZ2}. Facet persistence of protein 1L2Y is examined across three dimensions in Fig.~\ref{fig:2LYZ2}(b). All \(0\)–dimensional red bars terminate near \(3.8\,\text{\AA}\), marking the scale at which the C\(\alpha\) network becomes connected at the canonical nearest–neighbor spacing of protein backbones. The \(1\)–dimensional yellow facet bars emerge as short edges link adjacent C\(\alpha\) atoms along the \(\alpha\)–helical segment. Most of these edges disappear by \(5.5\,\text{\AA}\) when neighboring triangles form and absorb them, while a few persist up to \(7.0\,\text{\AA}\). These longer–lived bars correspond to edges spanning multiple helical turns, or involving terminal residues. The \(2\)–dimensional blue facet bars appear around \(5.0\,\text{\AA}\) as three–atom configurations within and across turns begin to form. Many triangles are quickly capped by tetrahedra, yet a subset persists to larger filtration values. The longer blue facet bars correspond to triangular faces that connect the compact helical core with the extended coil or terminal regions. The $f$-vector and $h$-vector are well-suited for structural data analysis by using C$\alpha$ atoms of protein 1L2Y in Fig.~\ref{fig:1L2Y2}(c) and Fig.~\ref{fig:1L2Y2}(d). The curve \(f_{0}\) is constant at \(20\), reflecting the fixed number of C\(\alpha\) vertices. The edge count \(f_{1}\) shows a distinct jump at \(3.8\,\text{\AA}\), when nearly all backbone neighbors connect simultaneously, then grows steadily as cross–turn edges accumulate within the helix, and as longer multi–turn and helix–tail contacts appear. Triangle counts \(f_{2}\) begin to rise rapidly around \(5\,\text{\AA}\) as nearby residues from adjacent turns close many triangles, and continue to grow as edge density becomes sufficient for widespread triangular closure. Tetrahedra counts \(f_{3}\) appear only after the complex contains many triangles, then increase sharply near \(8\,\text{\AA}\) as four–cliques close across multiple helical turns, with occasional participation from the tail, indicating dense local packing of residues.

Figure~\ref{fig:1L2Y2}(e)–(g) shows the evolution of the graded Betti numbers \(\beta_{i,i+j}\) for the \(\alpha\)–helical protein 1L2Y on three strands.  For strand \(j=1\), a broad plateau extends to \(3.8\,\text{\AA}\). Immediately beyond this threshold, the curves drop sharply, as nearly all C\(\alpha\)–C\(\alpha\) backbone edges appear and connect consecutive residues along the helix. The subsequent, gradual step–downs reflect the introduction of short cross–turn contacts between nearby helical turns, which progressively merge local substructures and reduce the number of independent connectivity patterns in the induced \((i+j)\)–subcomplexes. For strand \(j=2\), nonzero graded Betti numbers first appear near \( 4.0\,\text{\AA}\) and rise steeply to a dominant peak around \(5\,\text{\AA}\). This peak corresponds to numerous short loops created by local cross–turn diagonals inside the helix. As filtration length \(\varepsilon\) increases toward \(6\,\text{\AA}\), additional cross–turn and turn–to–turn edges act as diagonals that fill many of these loops, producing a rapid decline after the first peak. A second, broader rise appears between \(7\) and \(9\,\text{\AA}\), consistent with larger loops that span multiple helical turns or bridge the helix and terminal regions. Then, the curves decrease more gradually near \(9\,\text{\AA}\), as triangles and tetrahedra proliferate and progressively transform extended loops into faces of higher–dimensional simplices. Beyond \(10\,\text{\AA}\), the strand \(j=2\) curves decay toward zero as the complex becomes saturated and most one–dimensional cycles are filled. For strand \(j=3\), the signal summarizes the formation and capping of triangular shells. A pronounced, narrow peak typically appears near \(8\,\text{\AA}\), when cross–turn contacts are sufficiently common to complete many local four–vertex configurations. These transient shells register strongly in \(\beta_{i,i+3}\). With further growth of filtration length, additional edges supply the missing faces, the shells become boundaries of tetrahedra, and the signal decays quickly. Smaller, broader bumps may occur at \(11\)–\(12\,\text{\AA}\), arising from longer–range four–vertex groupings, including bridges between the helix and more flexible segments. These features are likewise short–lived and vanish once the complex is sufficiently saturated. Similar to \(\mathrm{C}_{20}\), the raw graded–Betti magnitudes are excessively large even for a 20–vertex system. 

We now examine hen egg–white lysozyme (PDB ID: 2LYZ), a 129-residue enzyme that adopts a compact, globular fold dominated by tightly packed $\alpha$-helices interconnected by short loops and turns. As shown in Fig.~\ref{fig:2LYZ2}(a), yellow spheres mark the C\(\alpha\) positions and the blue ribbon highlights the helical cores, while flexible connectors route between helices and shape the overall curvature. This C\(\alpha\) representation serves as the input for the Vietoris–Rips filtration and subsequent PH, PL, and PCA analyses.


PH and PL analyses of protein 2LYZ are shown in Fig.~\ref{fig:2LYZ} for three dimensions. 
In dimension \(0\), \(\beta_{0}\) bars corresponds to 129 C\(\alpha\) atoms in protein 2LY. At \(3.8\,\text{\AA}\), local clusters coalesce into a single connected component, leaving one surviving \(\beta_{0}\) class. Consistently, \(\lambda_{0}\) turns on and then increases, reflecting the steady enrichment of cross–turn and inter–helix edges that strengthen global connectivity even after \(\beta_{0}\) stabilizes. In dimension \(1\), \(\beta_{1}\) bars emerge between \(4\) and \(6\,\text{\AA}\). These loops arise when nonsequential residues across neighboring turns or adjacent \(\alpha\)–helices form cross–links that close short cycles around the helical scaffolds. As \(\varepsilon\) grows, triangles progressively fill these cycles, and deaths cluster between \(5\) and \(10\,\text{\AA}\). During this loop–rich interval, the nonharmonic spectrum \(\lambda_{1}\) remains small, indicating that many \(1\)–chains are weakly constrained. Beyond \( 10\,\text{\AA}\), as triangles and higher-dimensional simplices become widespread, \(\lambda_{1}\) rises, signaling increased geometric rigidity of the \(1\)–dimensional structure. Higher–dimensional topological features are sparse and short–lived. Intermittent \(\beta_{2}\) bars appear between \(7\) and \(12\,\text{\AA}\), indicating transient cavities within the compact lysozyme fold. In parallel, Once triangles become frequent around \(5\,\text{\AA}\), \(\lambda_{2}\) is well defined and then decreases until \(6\)–\(10\,\text{\AA}\), reflecting the formation of many triangular shells with few tetrahedral caps. Beyond \(9\)–\(10\,\text{\AA}\), \(\lambda_{2}\) turns upward slightly as tetrahedra proliferate, and it gradually levels off as the complex approaches saturation.

The large C\(\alpha\) atom count in lysozyme makes a full graded Betti table impractical to compute, since the number of triangles and tetrahedra grows combinatorially across the filtration. Even so, the structural evolution is faithfully reflected by facet persistence and by the progression of the \(f\)– and \(h\)–vectors. In Fig.~\ref{fig:2LYZ2}(b), all 0-dimensional red facet bars persist until \(3.8\,\text{\AA}\), consistent with the C\(\alpha\)–C\(\alpha\) spacing along the backbone. At this threshold, a dense cluster of one–dimensional yellow bars appears as backbone edges connect neighboring C\(\alpha\) atoms, then vanishes near \(6\,\text{\AA}\) when short–range triangles close within individual helical turns. A modest tail of longer one–dimensional bars persists until \(7\)–\(9\,\text{\AA}\), primarily in terminal or coil regions where closure requires longer edges. Two–dimensional blue facet bars correspond to persistent triangular faces. Births begin near \(5-10\,\text{\AA}\), and deaths extend through \( 6\)–\(11\,\text{\AA}\) as helices pack and tetrahedra gradually emerge through multi–turn and inter–helix contacts. The \(f\)–vector provides a more scalable aggregate view in Fig.~\ref{fig:2LYZ2}(c), though at the cost of losing detailed persistence information. The vertex count \(f_{0}\) is approximately \(129\), matching the number of C\(\alpha\) atoms. The edge count \(f_{1}\) begins to rise at \(3.8\,\text{\AA}\) with backbone connections, grows steadily as short cross–turn contacts appear, and accelerates beyond \(8\,\text{\AA}\) when inter–turn and inter–helix contacts proliferate. The triangle count \(f_{2}\) starts to increase near \(5\,\text{\AA}\) as close residue triplets within turns close, then steepens beyond \(9\,\text{\AA}\) as triangular faces become ubiquitous. The tetrahedra count \(f_{3}\) remains low until triangles are dense, then surges after \(9\,\text{\AA}\), reflecting numerous four–cliques spanning turns and helices. The \(h\)–vector complements this picture by tracking relative combinatorial growth derived from \(f\), shown in Fig.~\ref{fig:2LYZ2}(d). In this system, the number of triangles and tetrahedra can reach \(10^{4}\), which underscores the combinatorial richness of the fold.

\subsection{Protein-Nucleic Acid Complexes}

\begin{figure}
    \centering
     \includegraphics[width=0.8\linewidth, height=8cm]{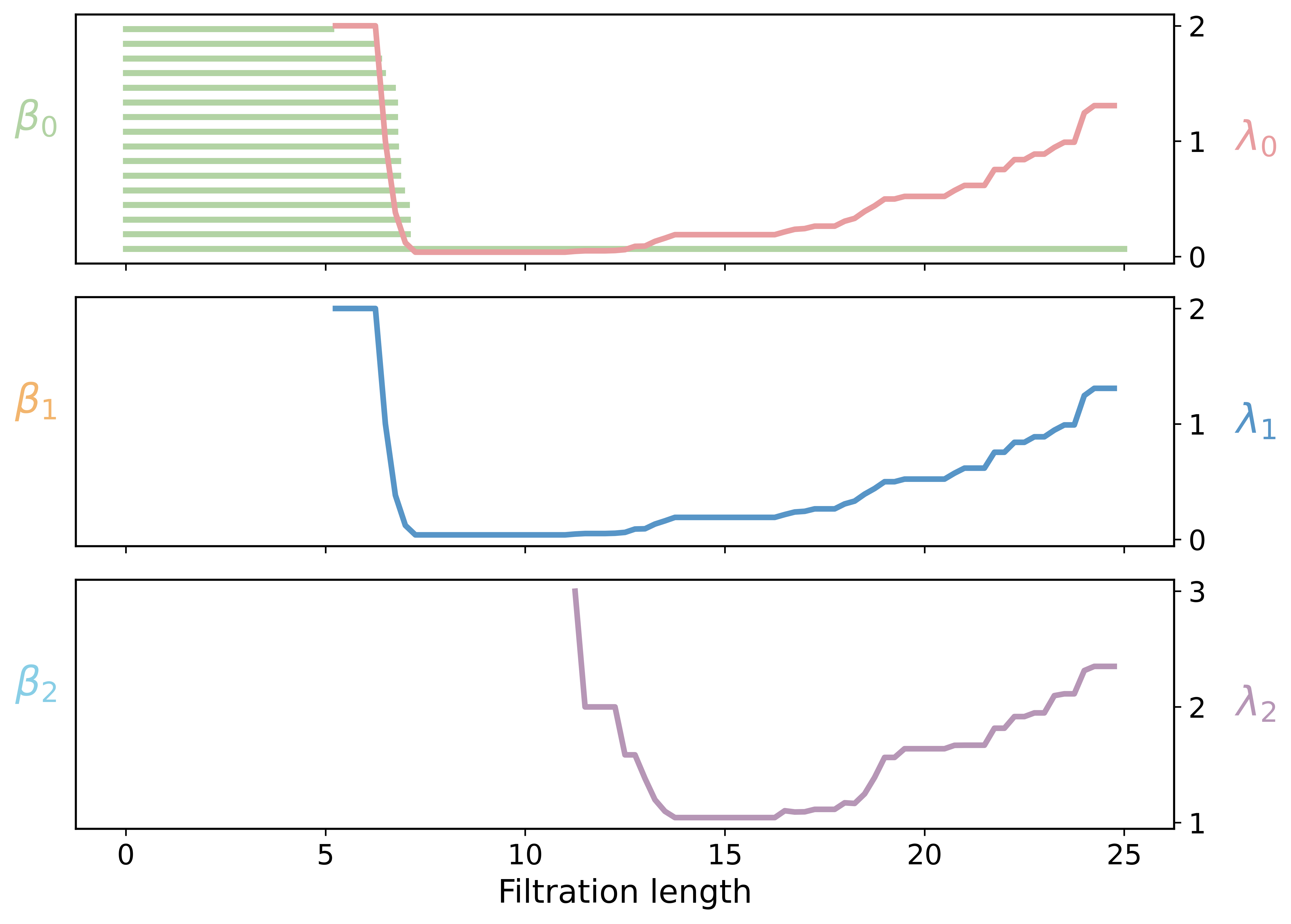}
     \caption{{Illustration of PH and PL on Protein-DNA complex (PDBID:1WET). Barcodes and spectra summarize the evolution of $0$-, $1$-, and $2$-dimensional features across the filtration.}}
         \label{fig:phl_1wet}
\end{figure}

\begin{figure}[htbp!]
  \centering

  \begin{subfigure}{0.42\linewidth}
    \centering
    \includegraphics[width=\linewidth,height=4.8cm]{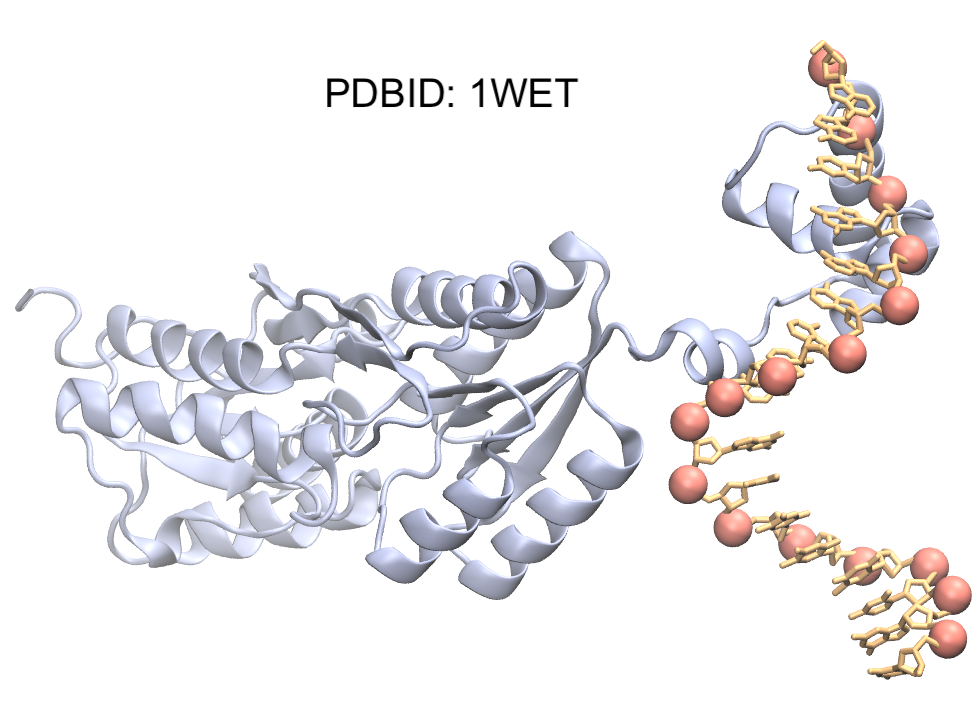}
    \subcaption*{(a) Protein--DNA complex 1WET}
  \end{subfigure}\hfill
  \begin{subfigure}{0.54\linewidth}
    \centering
    \includegraphics[width=\linewidth]{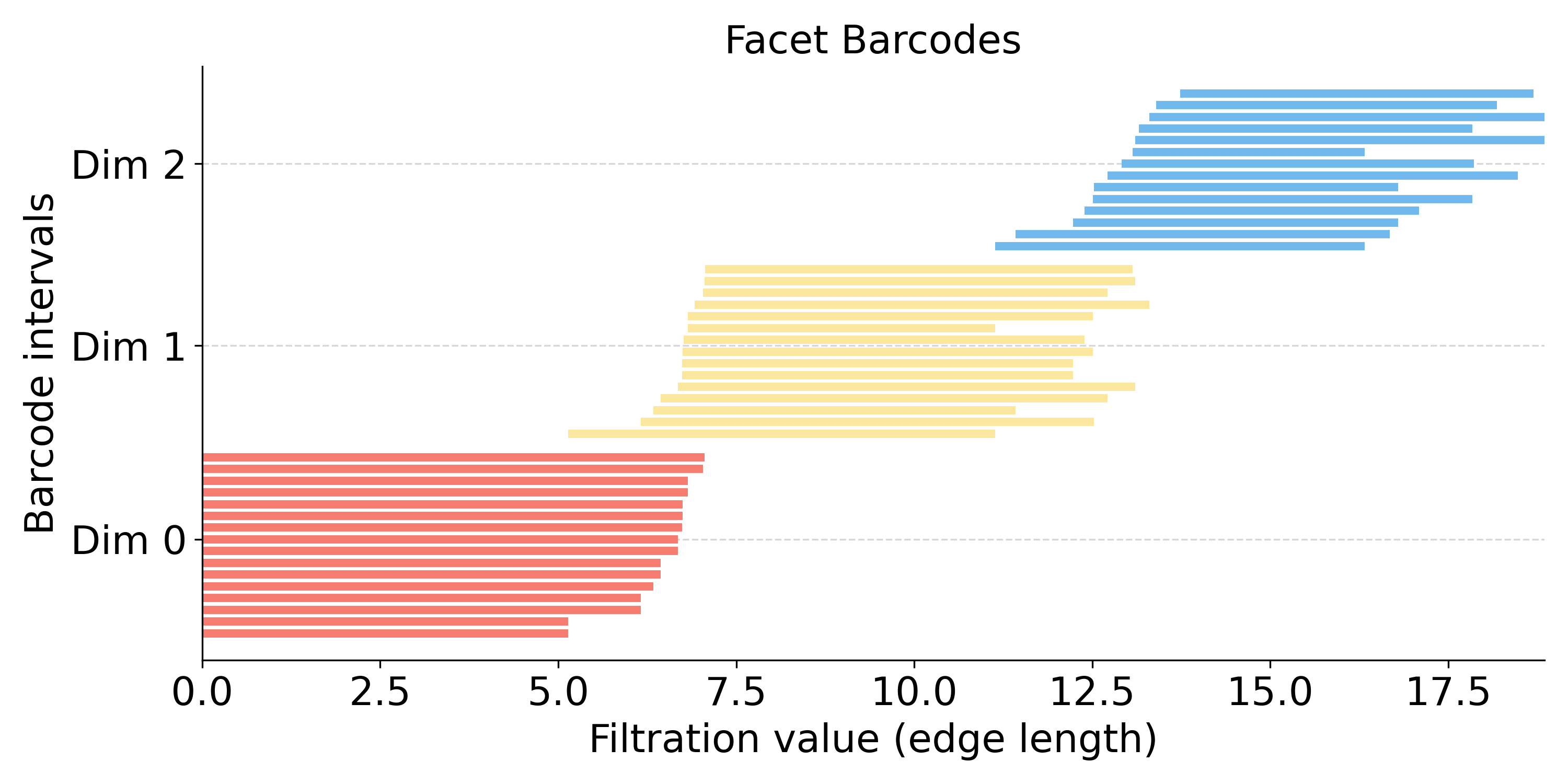}
    \subcaption*{(b) Facet persistence}
  \end{subfigure}

  \vspace{0.6em}

  \begin{subfigure}{0.49\linewidth}
    \centering
    \includegraphics[width=\linewidth, height=5.5cm]{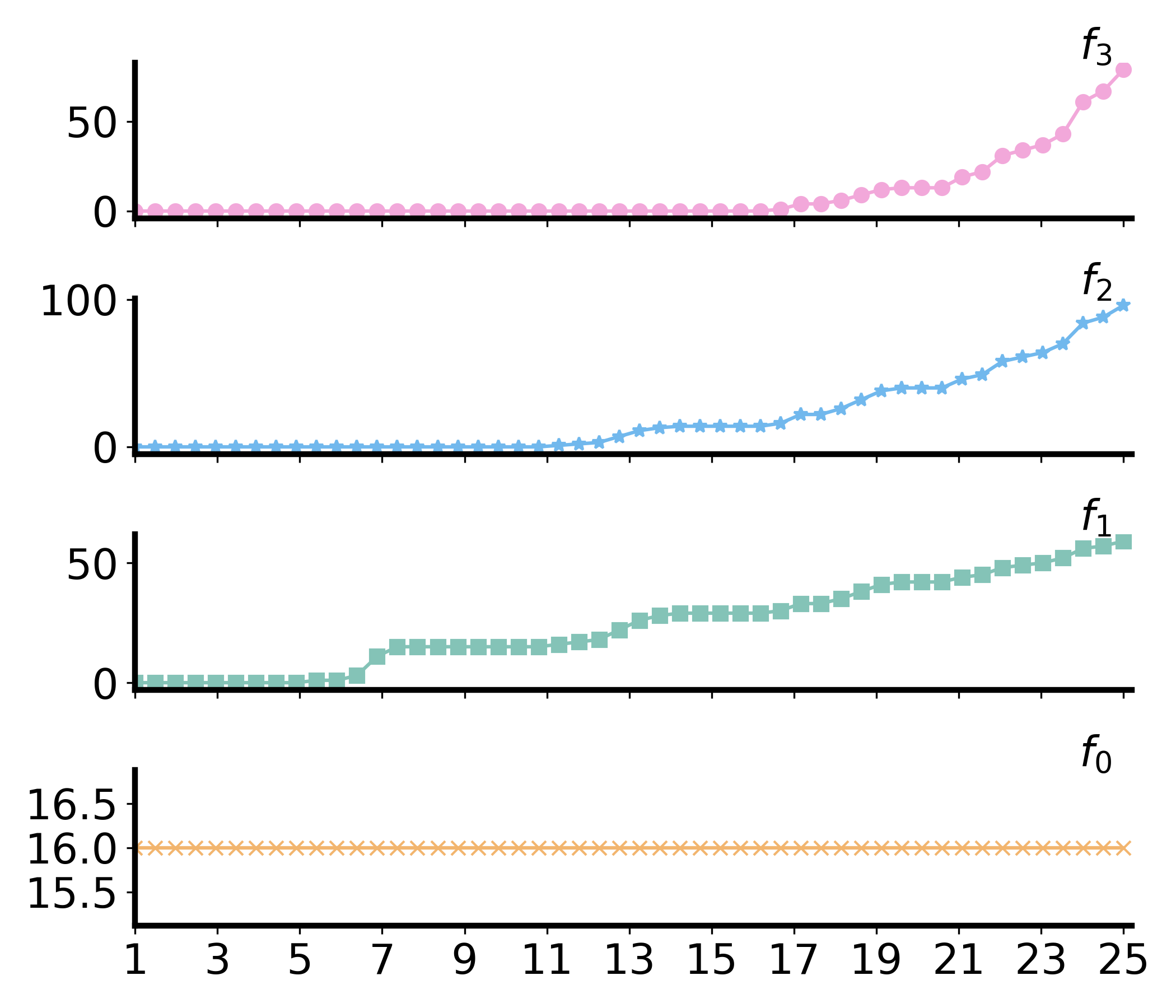}
    \subcaption*{(c) $f$-vector curves}
  \end{subfigure}\hfill
  \begin{subfigure}{0.49\linewidth}
    \centering
    \includegraphics[width=\linewidth, height=5.5cm]{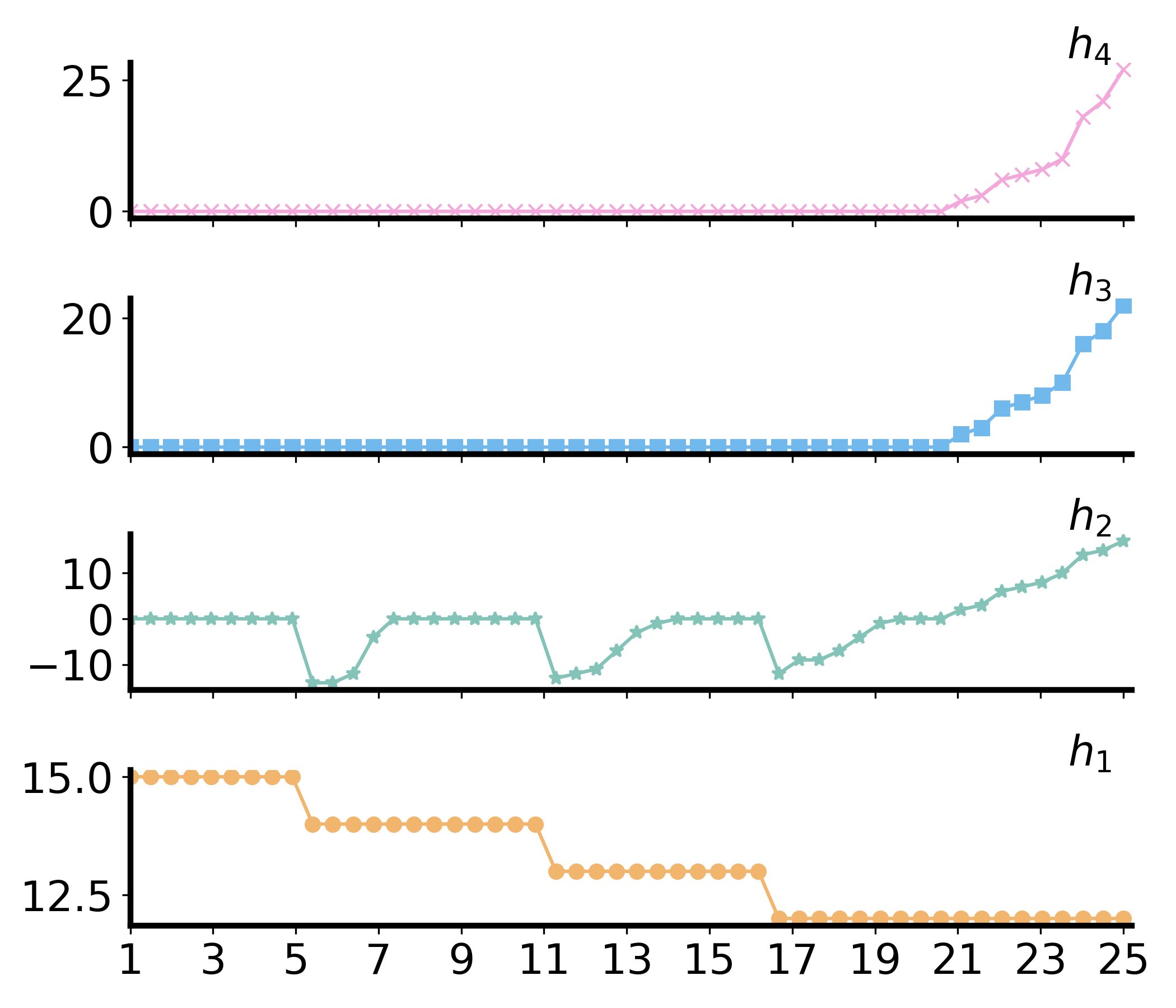}
    \subcaption*{(d) $h$-vector curves}
  \end{subfigure}

  \vspace{0.6em}

  \begin{subfigure}{0.5\linewidth}
    \centering
    \includegraphics[width=\linewidth,height=5cm]{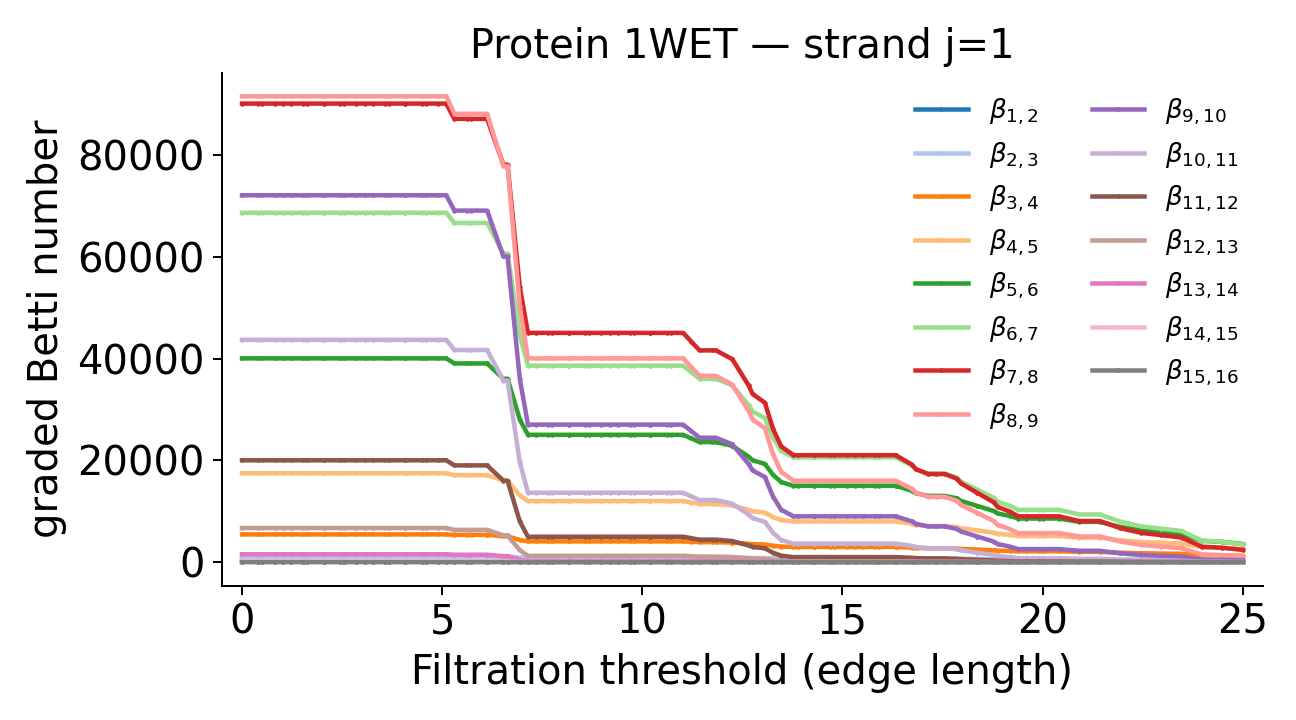}
    \subcaption*{(e) Graded Betti, dim$=0$ ($j=1$)}
  \end{subfigure}

  \caption{Persistent commutative algebra  analysis of the Protein--DNA complex (PDBID: 1WET) using a Rips complex-based filtration process.}
  \label{fig:pca_1wet}
\end{figure}

We begin with a representative protein–DNA system to demonstrate how PH, PL, and PCA characterize and represent structural features from the same filtration. Fig.~\ref{fig:pca_1wet}(a) shows the complex (PDB ID: 1WET). The protein ribbon embraces a short DNA fragment whose sixteen phosphate (P) atoms, displayed as pink spheres, furnish a sequence–ordered point cloud that traces the duplex backbone. Binding bends and locally untwists the fragment, so inter–phosphate distances vary along the strands and within the contact region. A Vietoris–Rips filtration on the P atoms therefore proceeds from along–strand links that reconstruct each DNA backbone, to cross–strand links that join paired phosphates across the duplex, and finally to longer–range links that capture inter-strand and intra-strand connections within the bent fragment.


The PH barcode and PL spectra for the phosphate cloud in the protein-DNA complex 1WET are shown in Fig.~\ref{fig:phl_1wet}. \(\beta_{0}\) bars correspond to sixteen phosphate atoms. No \(\beta_{1}\) or \(\beta_{2}\) bars appear, since the VR complex starts as two path graphs (no loops or cavities), and immediately fills any nascent edge cycles with triangles and any triangular shells with tetrahedra. The spectrum \(\lambda_{0}\) turns on when nearest–neighbor edges form, signaling the onset of local connectivity. As these edges extend into long chain–like segments, \(\lambda_{0}\) decreases toward a small plateau, which is characteristic of weak algebraic connectivity in nearly linear paths. Beyond \(12\,\text{\AA}\), \(\lambda_{0}\) rises again, as cross–strand and local intra-strand edges appear and strengthen the global network. Although \(\beta_{1}\) and \(\beta_{2}\) are absent, their nonharmonic spectra still record geometric transitions. The curve \(\lambda_{1}\) becomes defined once edges exist, then stays near zero across \( 7\)–\(12\,\text{\AA}\) while triangles remain scarce and extended paths dominate. When the filtration reaches the cross–strand and local intra-strand range, \( 11\)–\(13\,\text{\AA}\), triangles become common, leading to a rise of \(\lambda_{1}\). At the same time, \(\lambda_{2}\) first appears, then declines to a broad minimum near \(14\,\text{\AA}\), indicating that triangles accumulate faster than tetrahedra. At larger scales, \(\lambda_{2}\) increases again as tetrahedra form between successive turns and across the two strands, marking the approach to a geometrically saturated scaffold.

The PCA analysis of the DNA phosphate cloud in 1WET is shown in Fig.~\ref{fig:pca_1wet}. The facet evolution of 1WET is displayed in Fig.~\ref{fig:pca_1wet}(b). Red \(0\)-dimensional facets persist until \( 7\,\text{\AA}\), consistent with the typical phosphate–phosphate spacing along each strand. A dense block of \(1\)-dimensional yellow facets appears between \( 6\) and \(8\,\text{\AA}\), corresponding to intra–strand phosphate edges. Most of these \(1\)-dimensional bars disappear by \(13\,\text{\AA}\) as triangles close short loops formed by cross-strand connections and local intra-strand interactions within the slightly curved DNA fragment. Blue \(2\)-dimensional facets, which represent triangular faces, are born primarily between \(12\) and \(14\,\text{\AA}\) and persist through \( 15\)–\(18\,\text{\AA}\) until tetrahedra form. The \(f\)– and \(h\)– vector offers an aggregate view in Fig.~\ref{fig:pca_1wet}(c)(d). The vertex count \(f_{0}\) is constant at \(16\), matching the number of phosphate atoms. The edge count \(f_{1}\) rises at \( 6\,\text{\AA}\) with along–strand connections, then increases gradually as cross–strand and longer–range edges appear, reaching about \(55\) at the largest scale. The triangle count \(f_{2}\) remains negligible until \( 11\)–\(13\,\text{\AA}\), when triplets of phosphates become mutually close enough to form faces, after which \(f_{2}\) grows toward \(100\). The tetrahedra count \(f_{3}\) is absent until \(18\)–\(20\,\text{\AA}\), then increases, indicating the onset of dense four–vertex contacts spanning nearby phosphates and strands. However, the graded Betti numbers of the complex 1WET are nonzero only on the strand \(j=1\), shown in Fig.~\ref{fig:pca_1wet}(e). Induced subcomplexes of the phosphate cloud in the complex 1WET never sustain \(1\)– or \(2\)–dimensional homology, hence the graded Betti numbers on strand \(j=2\) and strand \(j=3\) vanish.

\begin{figure}
    \centering
     \includegraphics[width=0.8\linewidth, height=8cm]{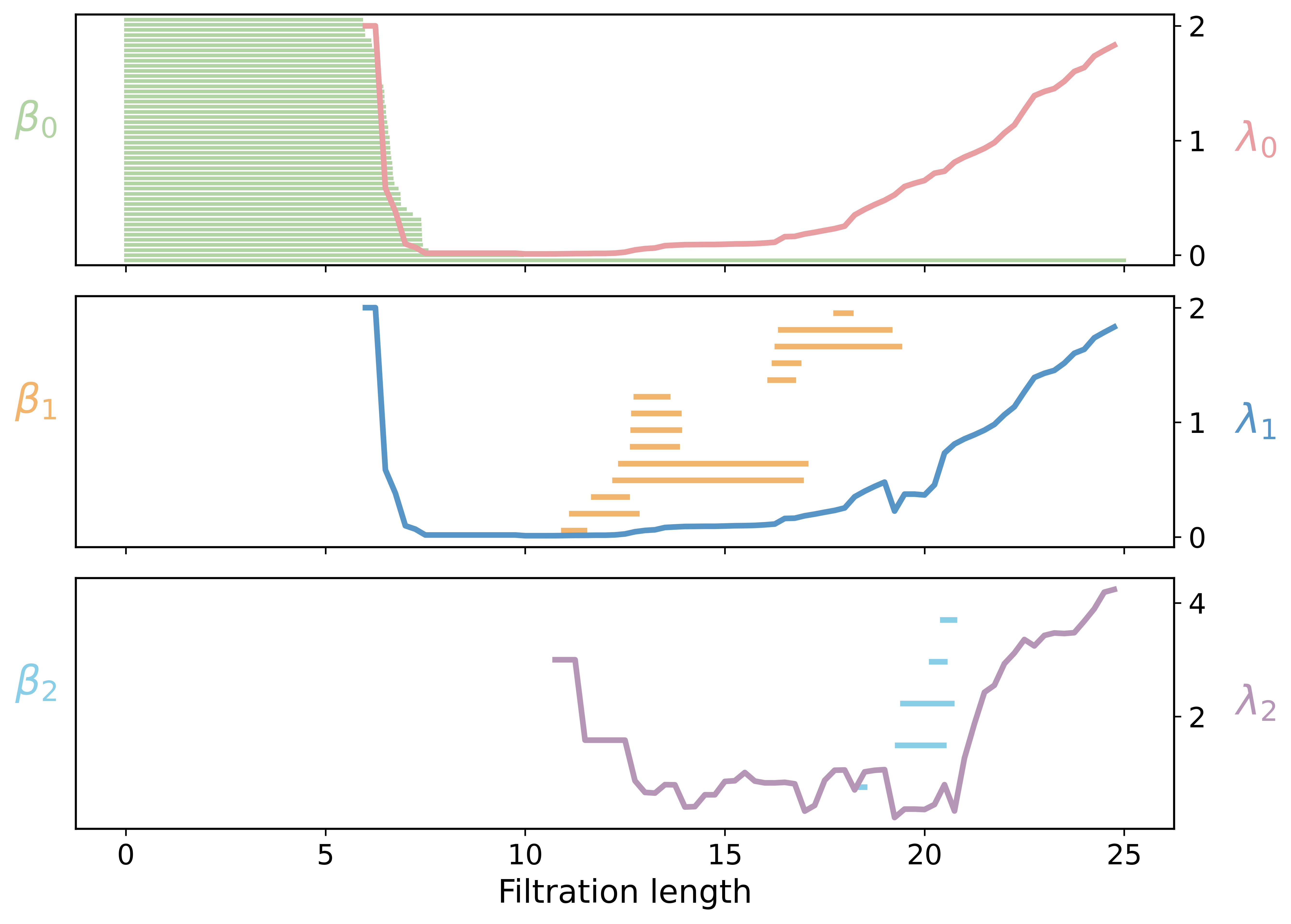}
     \caption{{Illustration of PH and PL on Protein-DNA complex (PDBID:1TW8). Barcodes and spectra summarize the evolution of $0$-, $1$-, and $2$-dimensional features across the filtration.}}
         \label{fig:phl_itw8}
\end{figure}

\begin{figure}[htbp!]
  \centering

  \begin{subfigure}{0.42\linewidth}
    \centering
    \includegraphics[width=\linewidth,height=5.5cm]{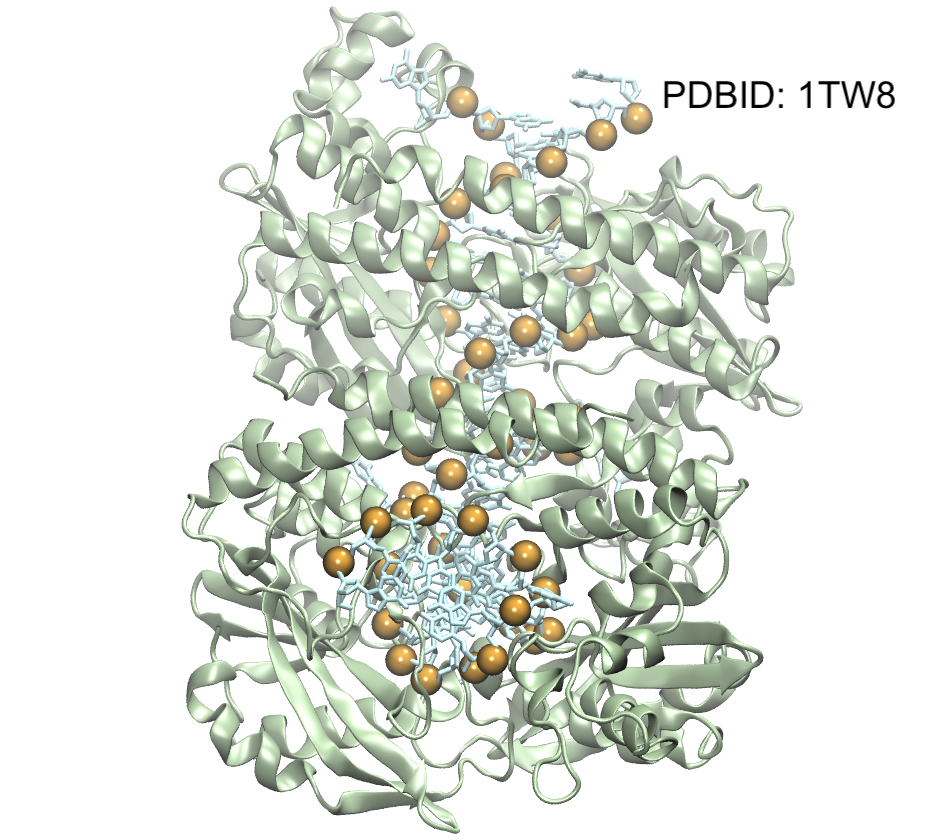}
    \subcaption*{(a) Protein--DNA complex 1TW8}
  \end{subfigure}\hfill
  \begin{subfigure}{0.54\linewidth}
    \centering
    \includegraphics[width=\linewidth,height=5.2cm]{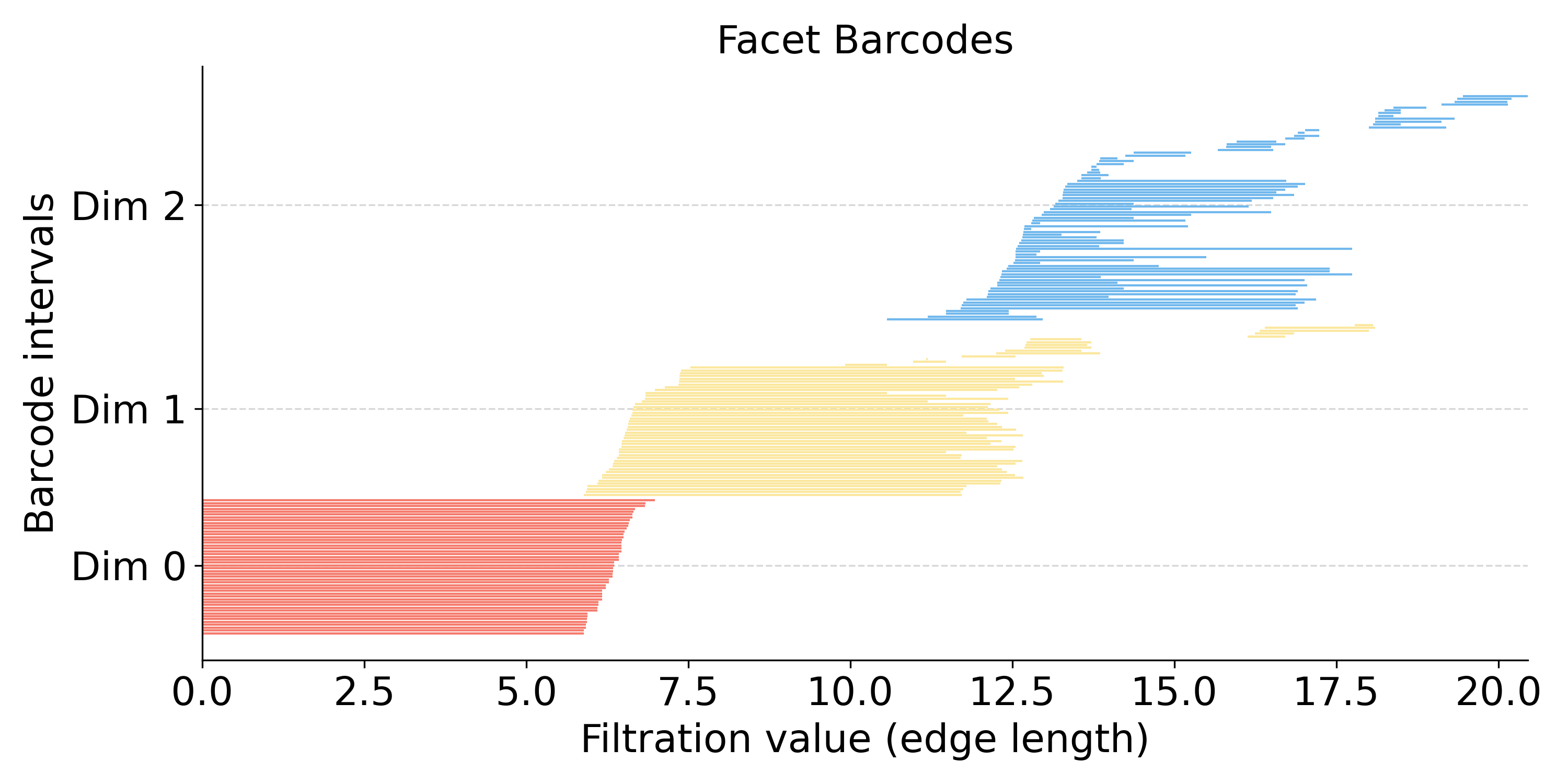}
    \subcaption*{(b) Facet persistence}
  \end{subfigure}

  \vspace{0.5em}

  \begin{subfigure}{0.49\linewidth}
    \centering
    \includegraphics[width=\linewidth]{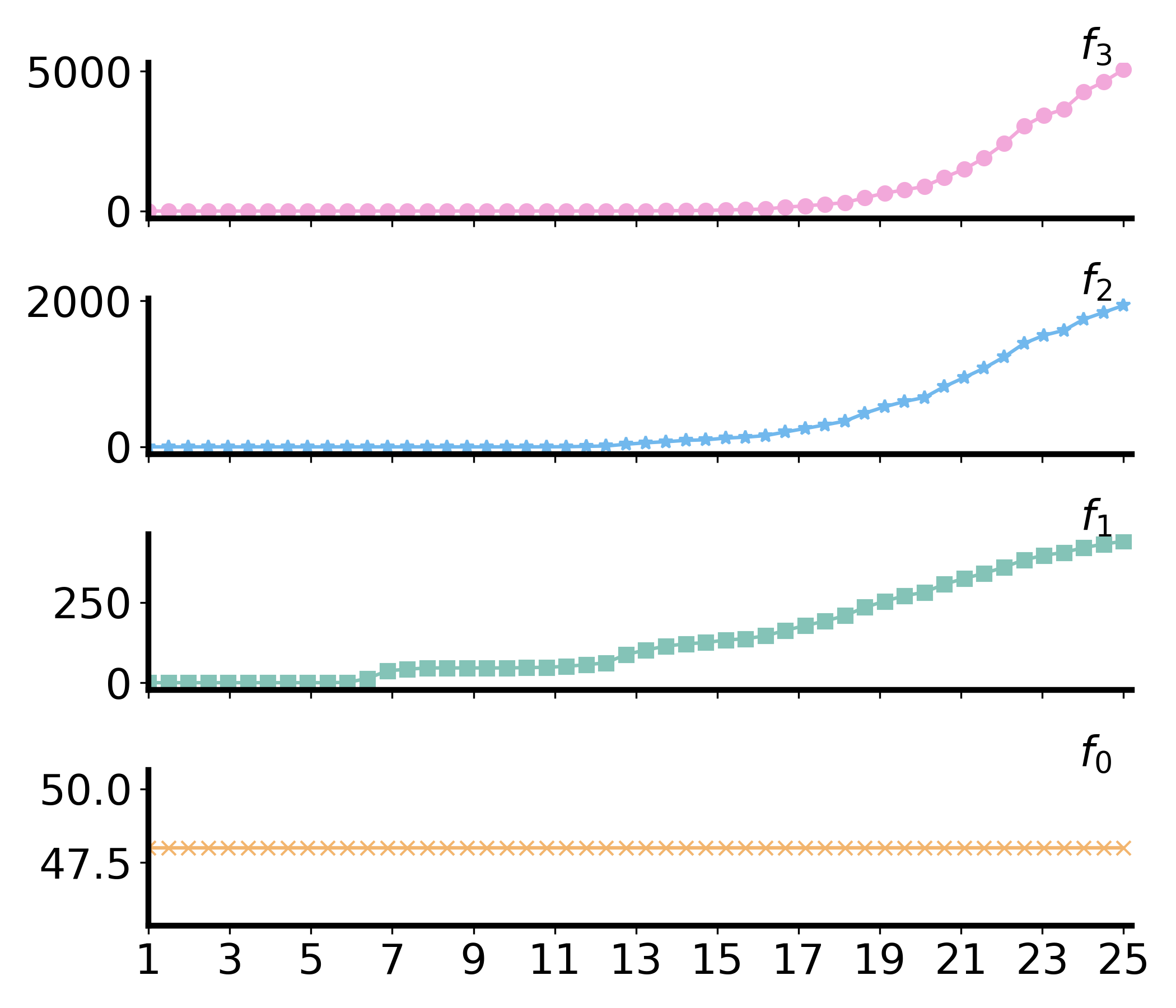}
    \subcaption*{(c) $f$-vector curves}
  \end{subfigure}\hfill
  \begin{subfigure}{0.49\linewidth}
    \centering
    \includegraphics[width=\linewidth]{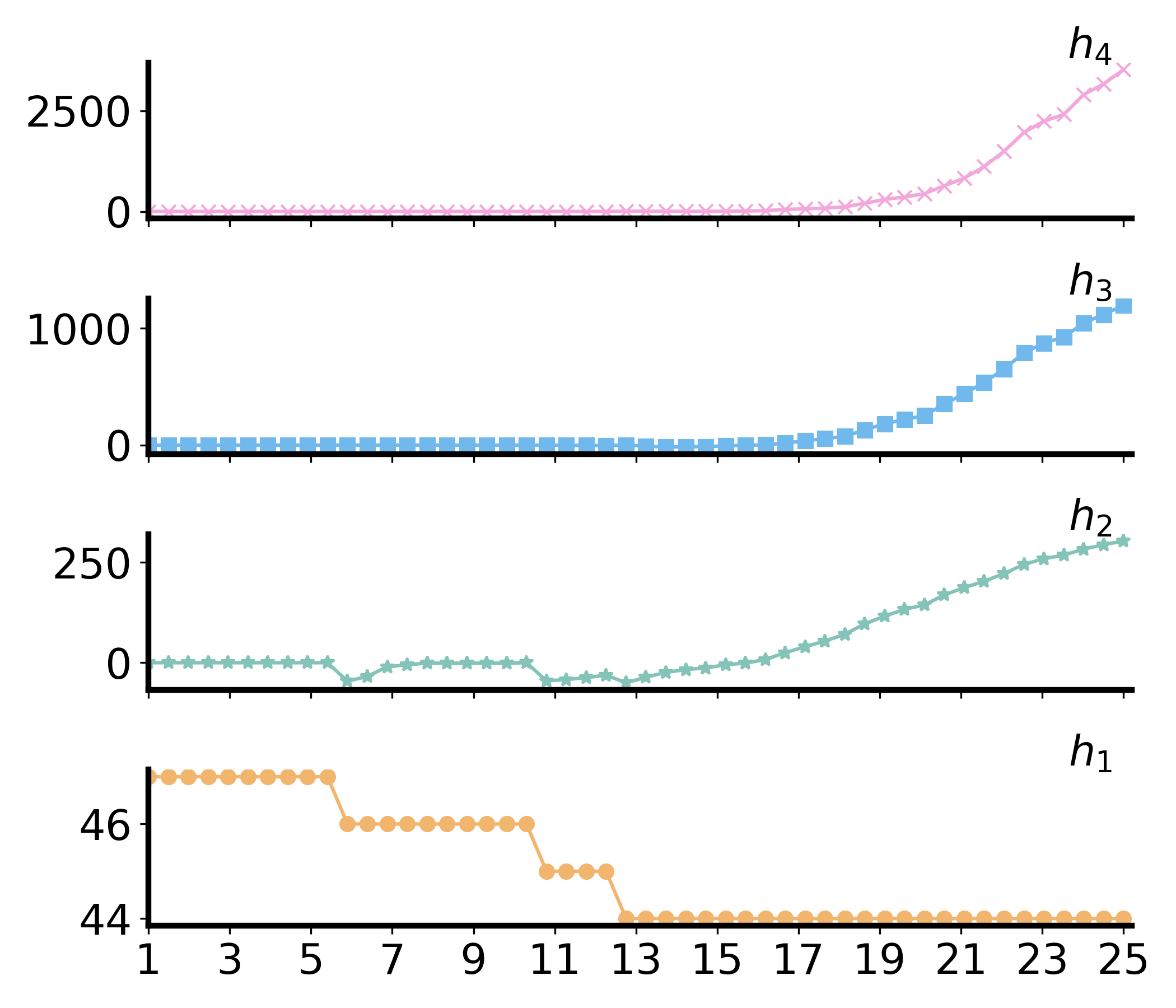}
    \subcaption*{(d) $h$-vector curves}
  \end{subfigure}

  \caption{Persistent commutative algebra  analysis of the Protein--DNA complex (PDBID: 1TW8) using a Rips complex-based filtration process.}
  \label{fig:pca_itw8}
\end{figure}

We next interpret and represent PH, PL, and PCA on the DNA–protein complex (PDB ID: 1TW8), shown in Fig.~\ref{fig:pca_itw8}(a). The protein forms a clamp around a DNA duplex, and the forty–eight phosphate atoms (orange spheres) provide a natural, sequence–ordered scaffold for analysis. Using the phosphates as vertices, the PH and PL results are summarized in Fig.~\ref{fig:phl_itw8}. In dimension \(0\), \(\beta_{0}\) bars represent forty–eight phosphate atoms. Beyond \( 6\,\text{\AA}\), the duplex backbone is globally connected. At the same scale, \(\lambda_{0}\) remains low, reflecting path–like connectivity along the strands, then increases at \( 13\,\text{\AA}\) as numerous cross–strand and inter–segment edges strengthen global connectivity. In dimension \(1\), two main groups of \(\beta_{1}\) bars appear, one for \(11\)–\(16\,\text{\AA}\) and another for \( 16\)–\(20\,\text{\AA}\). These bands correspond to the formation of short loops from inter–strand and inter–segment edges within the protein-bound DNA fragment. As \(\varepsilon\) increases beyond the along–strand P–P spacing, numerous short edge segments merge into long, nearly linear backbones on each DNA strand over the range \(6\)–\(18\,\text{\AA}\). Because triangles remain sparse in this interval, the network is primarily chain-like, leading to a decrease in \(\lambda_{1}\). Once \(\varepsilon\) exceeds \(\sim 19\,\text{\AA}\), triangular faces proliferate, resulting in a steady rise of \(\lambda_{1}\).
In dimension \(2\), a modest number of \(\beta_{2}\) bars are born near \( 18\)–\(20\,\text{\AA}\), marking transient triangular shells in locally crowded regions of the duplex. The spectrum \(\lambda_{2}\) becomes meaningful once triangles are present at \(11\)–\(13\,\text{\AA}\), then decreases toward a broad minimum across \( 15\)–\(20\,\text{\AA}\) due to triangles accumulating faster than tetrahedra. Beyond \(21\,\text{\AA}\), tetrahedra form across nearby phosphates and strands, resulting in a sharp increase of \(\lambda_{2}\).

Under the PCA framework, the phosphate count in 1TW8 renders a full graded–Betti computation infeasible, so we summarize structure via facet persistence and the \(f\)– and \(h\)–vectors in Fig.~\ref{fig:pca_itw8}. The facet persistence analysis is displayed in Fig.~\ref{fig:pca_itw8}(b). All \(0\)–dimensional bars persist until \( 6\,\text{\AA}\), consistent with the typical along–strand P–P spacing. One–dimensional (yellow) bars arise for \( 6\)–\(7\,\text{\AA}\), marking edge formation along each strand. Most of these edges die for \( 9\)–\(12\,\text{\AA}\) as triangles close local loops, while a subset persists to \( 16\)–\(18\,\text{\AA}\), capturing longer–range contacts between inter–strand and inter–segment edges. The broad range of deaths reflects geometric variability in spacing and curvature within the DNA fragment. In dimension \(2\), blue facet bars appear near \(12\,\text{\AA}\) when triplets of neighboring phosphates become mutually close and form triangles. These facets persist over a broad range, with deaths extending to \(20\,\text{\AA}\). The evolution of the \(f\)– and \(h\)–vectors in Fig.~\ref{fig:pca_itw8}(c,d) quantifies this progression, documenting the sequential enrichment from vertices to edges to triangular faces and, eventually, to higher–order simplices as the complex becomes more connected and geometrically saturated.

\begin{figure}
    \centering
     \includegraphics[width=0.8\linewidth, height=8cm]{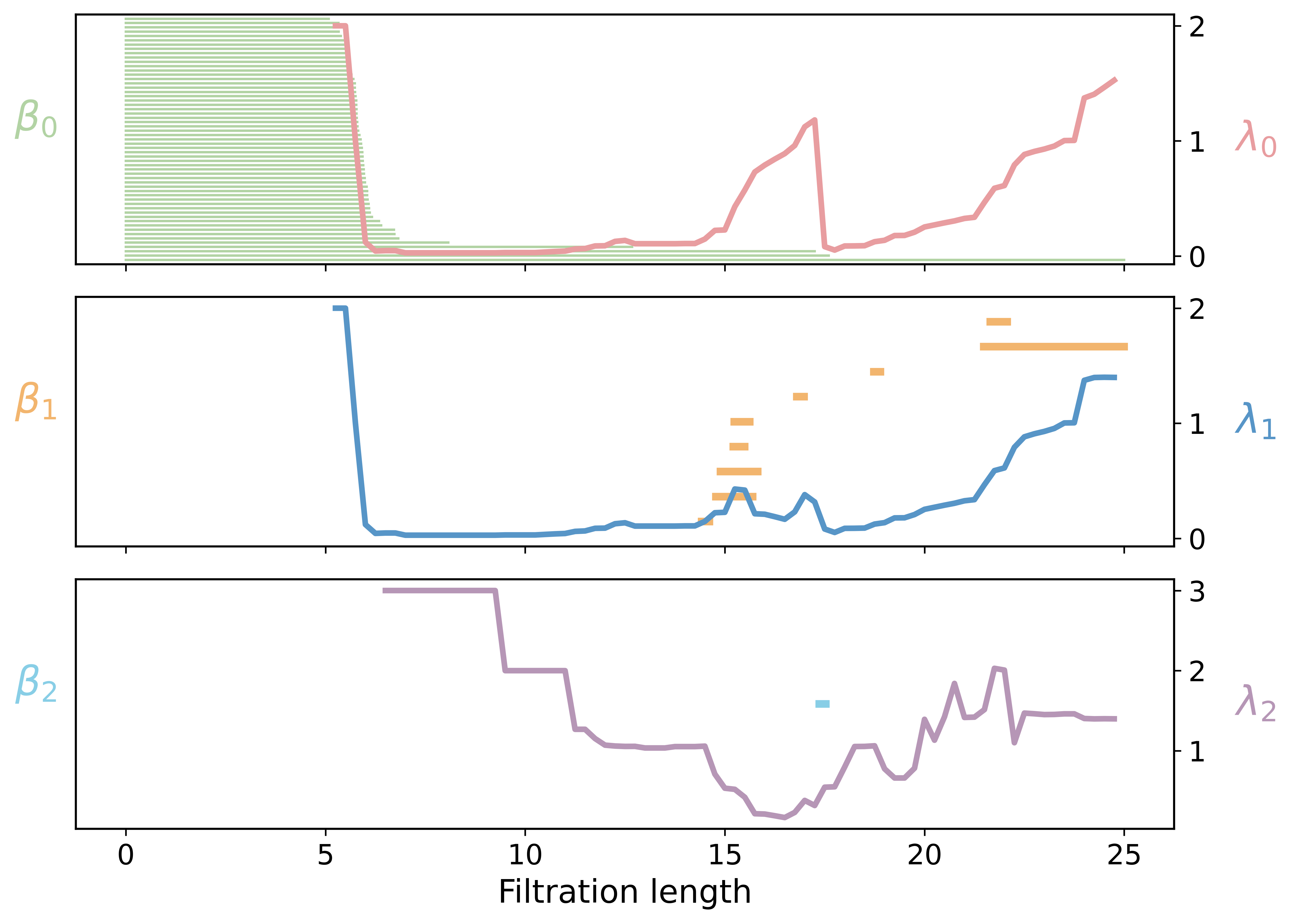}
     \caption{{Illustration of PH and PL on Protein-RNA complex (PDBID:1URN). Barcodes and spectra summarize the evolution of $0$-, $1$-, and $2$-dimensional features across the filtration. }}
         \label{fig:phl_1urn}
\end{figure}

\begin{figure}[htbp!]
  \centering

  \begin{subfigure}{0.42\linewidth}
    \centering
    \includegraphics[width=\linewidth,height=5.5cm]{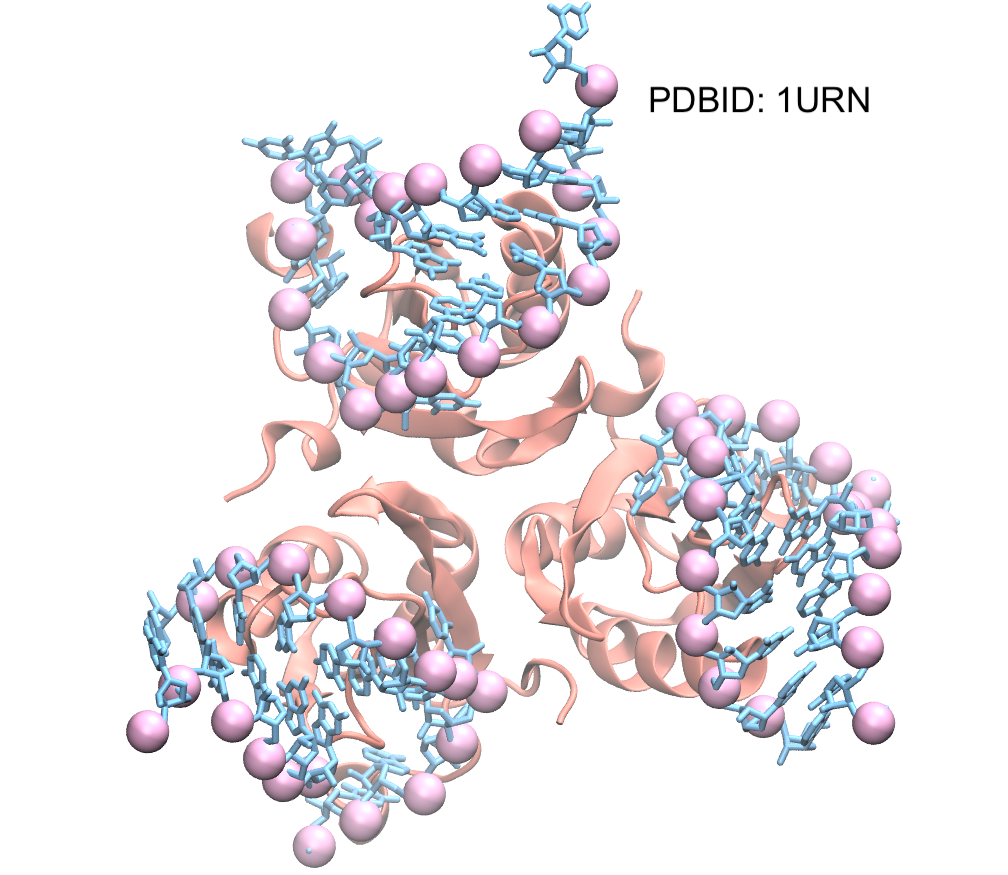}
    \subcaption*{(a) Protein--RNA complex 1URN}
  \end{subfigure}\hfill
  \begin{subfigure}{0.54\linewidth}
    \centering
    \includegraphics[width=\linewidth,height=5.2cm]{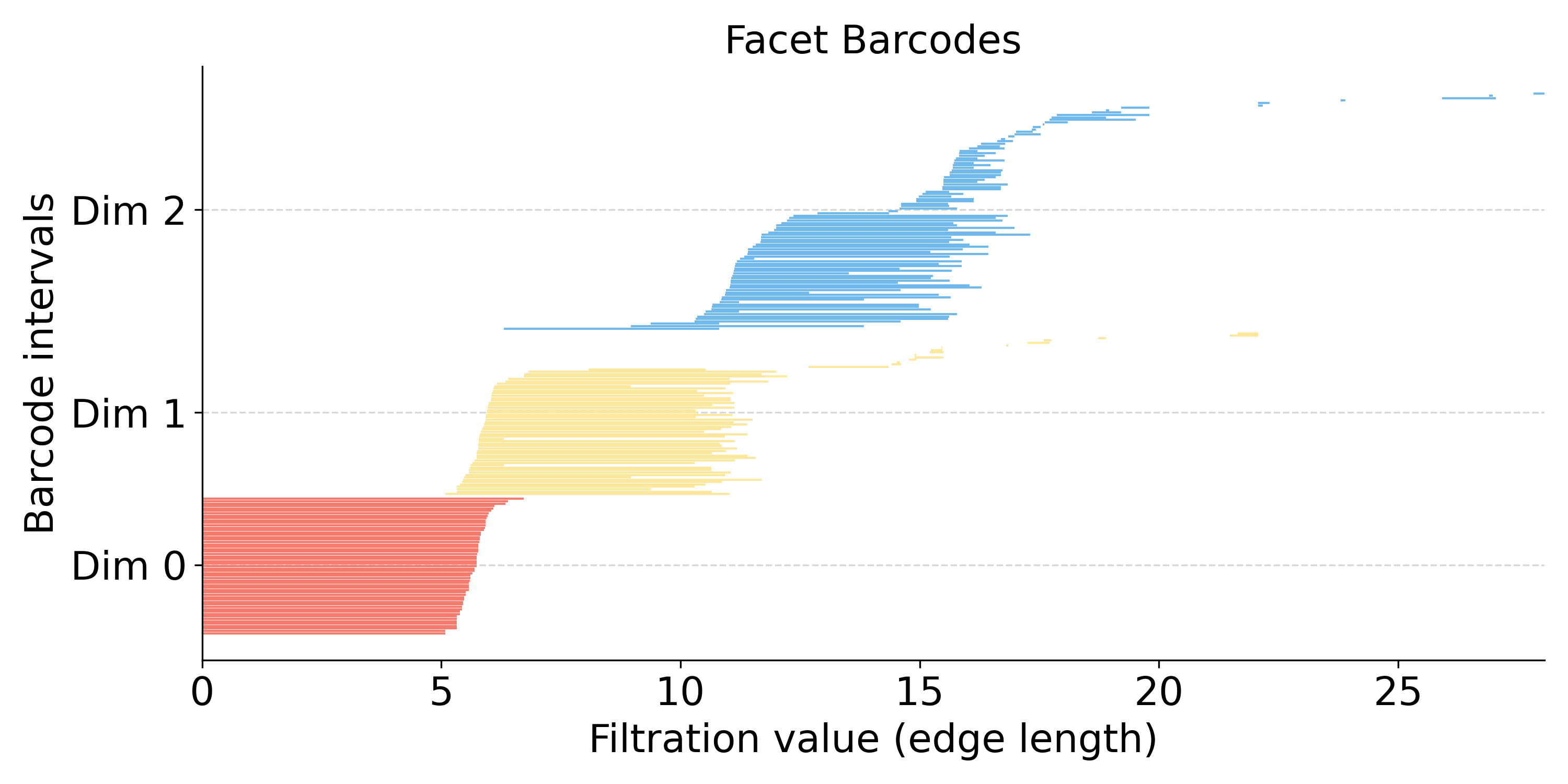}
    \subcaption*{(b) Facet persistence}
  \end{subfigure}

  \vspace{0.6em}

  \begin{subfigure}{0.49\linewidth}
    \centering
    \includegraphics[width=\linewidth]{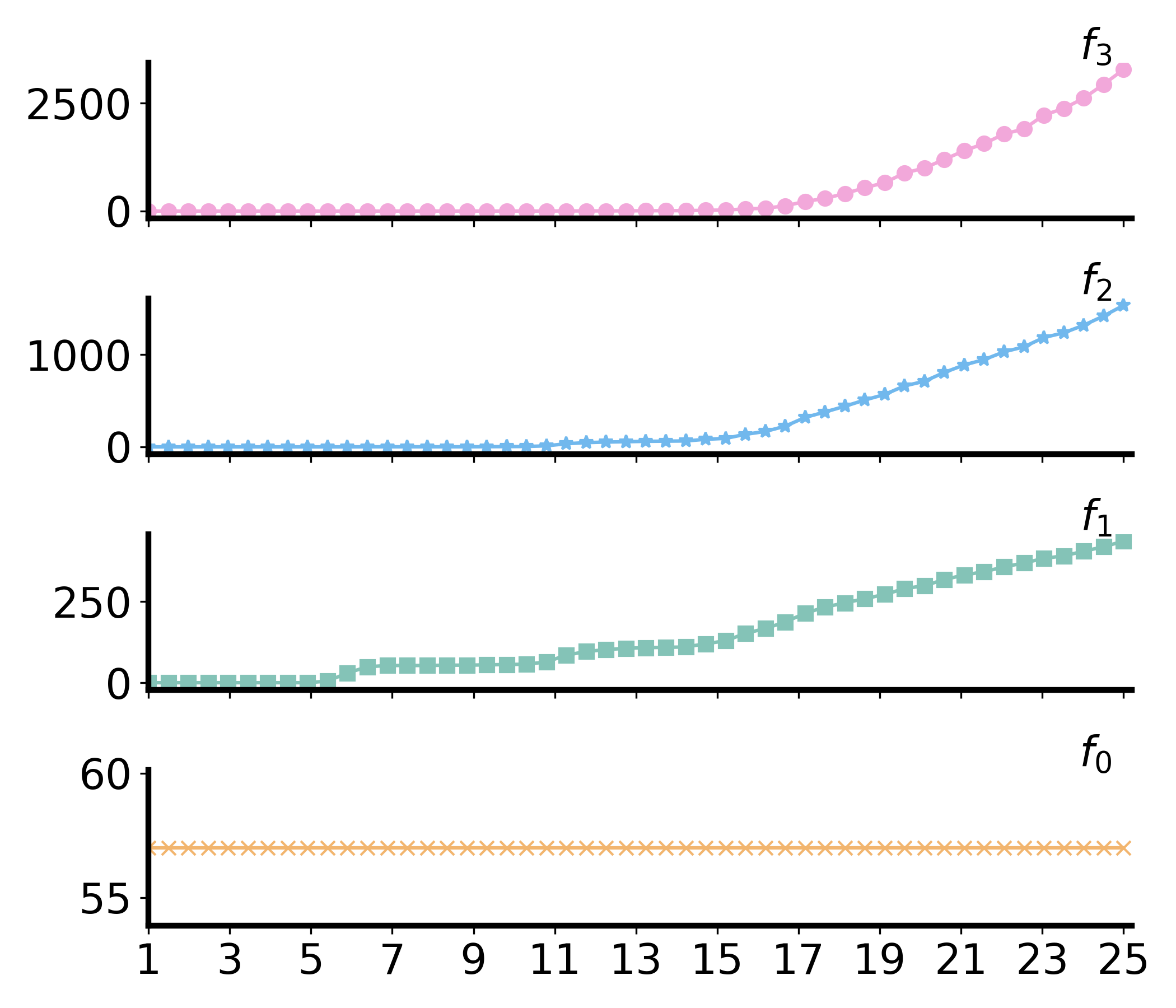}
    \subcaption*{(c) $f$-vector curves}
  \end{subfigure}\hfill
  \begin{subfigure}{0.49\linewidth}
    \centering
    \includegraphics[width=\linewidth]{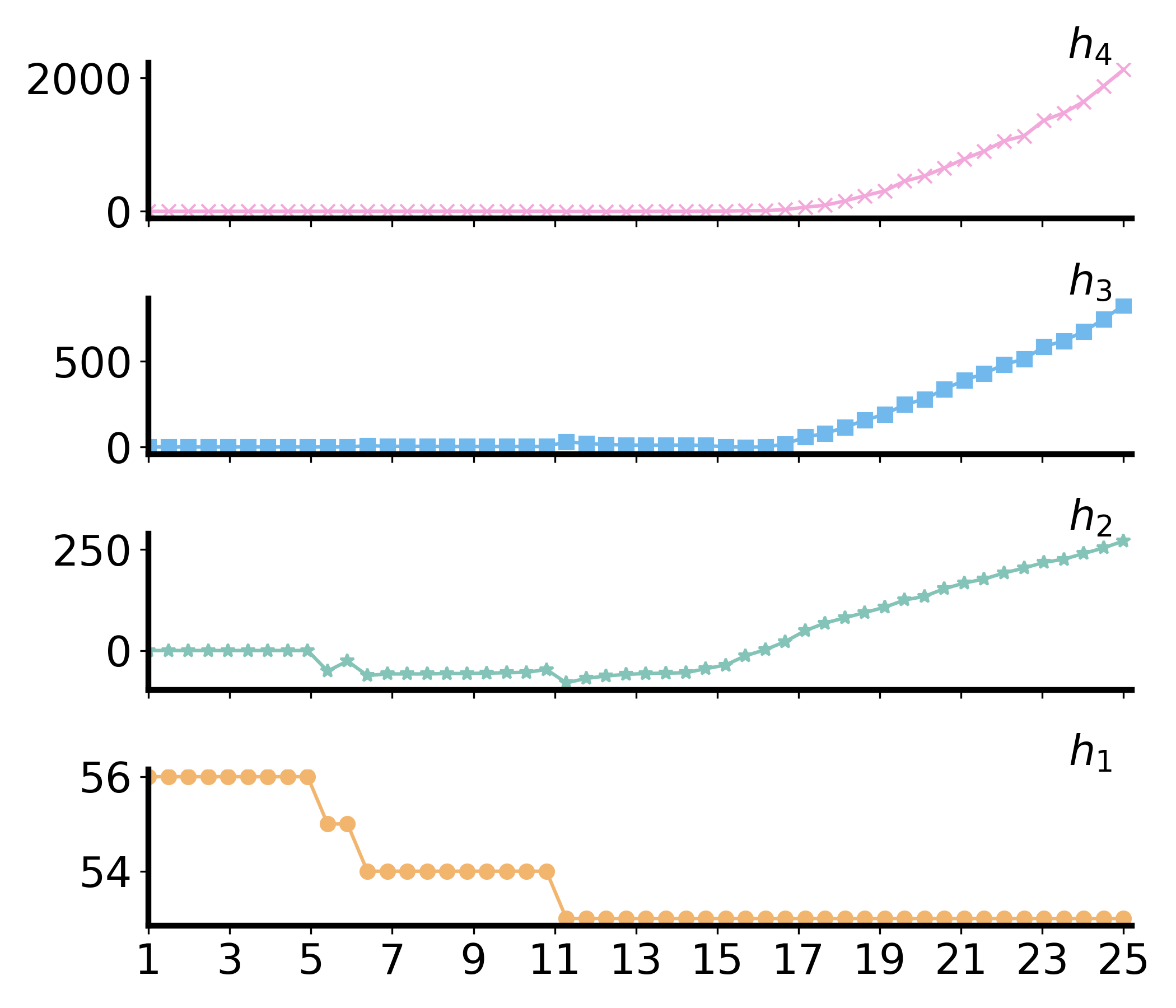}
    \subcaption*{(d) $h$-vector curves}
  \end{subfigure}

  \caption{Persistent commutative algebra analysis of the Protein--RNA complex (PDBID: 1URN) using a Rips complex-based filtration process.}
  \label{fig:pca_1urn}
\end{figure}

The protein–RNA complex (PDB ID: 1URN) contains an RNA A–form helical stem capped by a short loop, with the protein bound asymmetrically along one side and several glycerol molecules stabilizing the interface. The RNA fragment contributes \(57\) phosphate (P) atoms. As shown in Fig.~\ref{fig:pca_1urn}(a), purple spheres mark the RNA P atoms, the protein is rendered as a peach ribbon, and the RNA backbone is traced in blue. The RNA adopts a hairpin that wraps across a \(\beta\)–sheet surface of the protein, producing local bends, spatial loops, and a shallow pocket at the binding site. 

PH and PL for the phosphate cloud in 1URN are summarized in Fig.~\ref{fig:phl_1urn}. In dimension \(0\), \(57\) initial components merge first within strands at \( 5\,\text{\AA}\), consistent with the along–strand P–P spacing. The structure is fully connected by about \( 18\,\text{\AA}\), as cross–strand and inter–segment edges join stems and loop regions into a single component. In dimension \(1\), \(\beta_{1}\) bars occur mainly for \(15\)–\(25\,\text{\AA}\) and reflect loops formed by long links that bridge across hairpin stems or between neighboring RNA segments. These features are short–lived because added edges quickly supply the missing faces. After the first edges appear, the spectrum \(\lambda_{1}\) drops and remains small for \(6\)–\(18\,\text{\AA}\) once those edges merge into long, open paths, and then rises beyond \(\sim 19\,\text{\AA}\) as triangles accumulate and cycles become clamped.
In dimension \(2\), one \(\beta_{2}\) bar appears only at larger radii, \(\varepsilon \gtrsim 15\,\text{\AA}\), corresponding to a transient cavity in tightly packed regions of the RNA. Consistently, \(\lambda_{2}\) exhibits a shallow minimum near \(16\)–\(18\,\text{\AA}\), when triangular shells are most prevalent, and then increases as nearby phosphates form tetrahedra and cap the shells. Overall, the Betti intervals capture the topological transitions from disconnected strands to a single connected backbone with short–lived loops and rare voids, while the nonharmonic spectra quantify the accompanying geometric stiffening of the RNA phosphate network observed in the protein–RNA complex.

For the protein–RNA complex 1URN, the large number of phosphate atoms makes the computation of the graded–Betti table computationally prohibitive. 
Instead, we summarize the structural organization through facet persistence and the associated \(f\)– and \(h\)–vectors, as shown in Fig.~\ref{fig:pca_1urn}. 
The facet persistence barcode in Fig.~\ref{fig:pca_1urn}(b) reveals a clear hierarchical organization of the RNA backbone. 
All zero–dimensional facets persist up to approximately \(6\,\text{\AA}\), indicating the scale at which isolated phosphates merge into continuous backbones. 
One–dimensional yellow facets appear near \(5\,\text{\AA}\), representing local covalent connections along individual RNA strands. 
Most of these edges disappear between \(10\) and \(12\,\text{\AA}\) as they are absorbed into triangular faces, while a subset persists to \(15\)–\(22\,\text{\AA}\), capturing cross–strand and inter–hairpin edges. 
Two–dimensional blue facets emerge near \(6\,\text{\AA}\) and accumulate across \(12\)–\(16\,\text{\AA}\), forming planar phosphate triplets within and between neighboring helical stems. 
Many of these triangles persist across a broad range of scales, extending to nearly \(27\,\text{\AA}\). 
This wide distribution of birth and death scales, particularly between \(6\,\text{\AA}\) and  \(27\,\text{\AA}\), reflects the geometric heterogeneity of the 1URN complex, which combines compact helices, flexible loops, and extended inter–hairpin connections stabilized by the surrounding protein. 
The corresponding \(f\)– and \(h\)–vectors, presented in Fig.~\ref{fig:pca_1urn}(c) and (d), provide quantitative summaries of this multiscale structure by encoding the combinatorial balance among vertices, edges, and higher–dimensional simplices throughout the filtration.

\section{Discussion }

After applying PH, PL, and PCA to a range of representative examples, we compare these three frameworks in terms of what they measure, the information they capture, and their trade-offs between interpretability, geometric sensitivity, and computational scalability in Table~\ref{Tab:comp}. Each method analyzes the same filtration process but interprets it through a different mathematical lens, topological, spectral, or algebraic.

PH remains the most widely used due to its conceptual simplicity, computational efficiency, and intuitive outputs. It tracks the birth and death of topological features such as independent connected components, loops, and voids through persistence barcodes or diagrams, providing a purely topological summary of data shape. PH’s main limitation is its insensitivity to geometric intensity: it treats all edges as binary (present or absent) and therefore cannot distinguish between dense and sparse regions of connectivity. Despite this abstraction, PH scales efficiently to large datasets and provides a stable, interpretable foundation for topological analysis across diverse scientific applications.

PL extends PH by introducing a spectral perspective that integrates both topological and geometric information. The harmonic spectra (zero eigenvalues) correspond to PH’s Betti numbers, ensuring consistency with homological invariants, while the non-harmonic spectra (positive eigenvalues) encode geometric stiffness, local edge density, and redundancy of cycles. PL therefore quantifies not only whether a feature exists, but also how strongly it is geometrically supported. However, this added sensitivity comes at a computational cost: solving eigenvalue problems at each filtration scale is significantly heavier than matrix reductions in PH, limiting PL’s scalability on large or dense complexes.

PCA offers four types of features, and they have different computational costs. Rooted in Persistent Stanley–Reisner theory, PCA reinterprets filtration through algebraic invariants such as facet ideals, $f$– and $h$–vectors, and graded Betti numbers, thereby capturing both combinatorial and algebraic dependencies within the complex. Its strength lies in its ability to unify topological, combinatorial, and algebraic perspectives, providing a structured hierarchy of dependencies beyond what PH or PL can represent. However, its computational complexity and scalability depend strongly on the chosen invariant. For example, $f/h$ analyses are efficient and scale with the number of faces, whereas graded Betti computations involve minimal free resolutions and computing syzygies between generators, which can become impractical for large complexes or data. Moreover, the raw counts of the algebraic invariants can become extremely large, necessitating normalization and scaling for consistent comparison across datasets. Despite these challenges, PCA complements PH and PL by providing algebraic interpretability and multiscale structural detail.

\noindent\textbf{Computational complexity.}
Let $\Delta_\varepsilon$ be the Vietoris--Rips complex at scale $\varepsilon$ built up to dimension $D$ with $s_k(\varepsilon)$ $k$–simplices and $S(\varepsilon)=\sum_{k\le D}s_k(\varepsilon)$. 
PH forms sparse boundary matrices $B_k(\varepsilon)\in\mathbb{F}^{\,s_{k-1}\times s_k}$ over a finite field $\mathbb{F}$ and performs one reduction per sampled scale, with number of nonzeros $\,\mathrm{nnz}(B_k)\!=\!O((k{+}1)s_k)$, practical costs are near
\[
\text{time}_{\mathrm{PH}}\;\approx\; \sum_{\varepsilon}\sum_{k\le D} {O}\!\big(\mathrm{nnz}(B_k(\varepsilon))\big)\]
while the theoretical worst case is $O\!\big(\sum_{\varepsilon}\sum_{k\le D} s_k(\varepsilon)^{\,3}\big)$.
PL at each sampled scale solves small‑eigenpair problems for
\[
L_k(\varepsilon)=B_k(\varepsilon)^{\!\top}B_k(\varepsilon)+B_{k+1}(\varepsilon)\,B_{k+1}(\varepsilon)^{\!\top}\ \in \mathbb{R}^{\,s_k\times s_k}.
\]
if full spectra were computed, the cubic bound $O\!\big(\sum_{\varepsilon}\sum_{k} s_k(\varepsilon)^3\big)$ would apply. It uses floating–point (real/complex) arithmetic and is typically slower per unit work than finite‑field PH. 
PCA has two very different regimes: the facet/$f$–/$h$–curve side and the graded Betti side. Counting faces and converting to $h$–vectors scales roughly with the total number of simplices present,
\[
\text{time}_{\mathrm{PCA}(f/h)}\;\approx\;\sum_{\varepsilon}\sum_{k\le D}{O}\!\big(s_k(\varepsilon)\big),\]
and is usually inexpensive. By contrast, graded Betti numbers require minimal free resolutions of the Stanley–Reisner ideal, which can exhibit super-polynomial (and in practice, exponential) complexity in the number of variables, generators, and their degrees. Facet persistence, on the other hand, is more tractable. It only requires updating and tracking maximal simplices across the filtration, which is combinatorial in nature and typically less costly than full free resolutions, though still more demanding than simple $f-$ or $h-$ vector counts.

In practice, PH remains the fastest and most scalable for fixed $D$ and moderate filtration density. PL is costlier due to repeated eigenvalue solves but yields valuable geometric information. PCA is efficient for $f/h$–based analyses but can become intractable for graded Betti computations on dense or high-dimensional complexes. Together, these methods form a hierarchy of trade-offs between algebraic richness, geometric sensitivity, and computational scalability.

\begin{center}
\renewcommand{\arraystretch}{1.15}
\begin{tabular}{p{3.2cm}p{4.2cm}p{4.2cm}p{4.2cm}}
\toprule
 & \textbf{PH} & \textbf{PL} & \textbf{PCA} \\
\midrule
\textbf{Object tracked}
& Birth and death of homological features (Betti numbers) 
& Harmonic  and non-harmonic spectra of $L^{\,i,j}_k$
& Facet ideals, $f/h$–curves, and graded Betti numbers \\
\midrule
\textbf{Primary signal}
& Purely topological invariants
& Topology + geometric stiffness / connectivity
&  Topological and algebraic invariants + combinatorial traits.\\
\midrule
\textbf{Strengths}
& Fastest, widely scalable; Stable barcodes and persistence diagrams; Intuitive visualization 
& Unifies topology and geometry; Rich Geometric Sensitivity
& Multiscale topology, algebra and combinatorics; Rich interpretability\\
\midrule
\textbf{Limitations}
&Insensitive to local geometry and connection strength
& Computationally expensive for large data, less scalable 
&Complexity and scalability depend strongly on the choice of invariants;  Scaling may be needed\\
\midrule
\textbf{Typical software}
& Ripser\cite{bauer2021ripser}, GUDHI\cite{gudhi2015gudhi}, Dionysus\cite{Dionysus}
& HERMES\cite{wang2021hermes}, PETLS\cite{jones2025petls}
& Macaulay2\cite{grayson2002macaulay2}\\
\midrule
\end{tabular}
\captionof{table}{Summary comparison of three persistent frameworks: PH, PL, and PCA.}
\label{Tab:comp}
\end{center}

\section{Conclusion}

Artificial intelligence (AI) has revolutionized science, engineering, and technology in the past decade. However, current AI faces challenges in interpretability, generalizability, and transparency, among other issues. Mathematical AI, a rapidly growing field, offers explainable AI (xAI) and enables generalizable AI models for a wide variety of data formats. Some of the most promising mathematical AI models are based on persistent homology (PH). Recently, persistent Laplacians (PLs), a technique rooted in spectral theory, and persistent commutative algebra (PCA) have demonstrated superior performance in mathematical AI and its real-world applications.

In this work, we examine the interpretability and representability of PH, PL, and PCA. To this end, we utilize the same filtrations and datasets—ranging from geometric examples and synthetic complexes to fullerenes and biomolecular assemblies—to clarify what each framework measures, how its summaries relate to geometry, and where practical limitations arise.

PH remains one of the most transparent representations of multiscale topology, delivering stable barcodes of topological invariants that can be interpreted as independent components, loops, and voids. It has already proved to be effective in molecular and bimolecular applications with strong computational efficiency and clear visualization. PL retains the homological content of PH in its harmonic spectra and adds non-harmonic spectra that quantify geometric reinforcement. This aligns with the growing evidence that spectral features improve predictive modeling across scientific domains. PCA is a new nonlinear algebraic approach that extends the representative frontier by encoding facet structure, $f-$ and $h-$ vectors, and graded Betti tables, thereby exposing higher order combinatorial dependencies and giving algebraic invariants of how simplices assemble across scales. 

Each method comes with characteristic costs and benefits. PH is typically the fastest and most scalable, but it is largely insensitive to local geometric shape evolution. PL increases geometric sensitivity through eigen analysis, which is numerically heavier and requires careful interpretation of spectral variation. PCA offers the richest structural accounting. In particular, $f-$ and $h-$ summaries are efficient, while the computation of graded Betti numbers may become prohibitive for large complexes, and raw algebraic counts may require normalization across filtration scales for fair comparison. Together, these three frameworks form a hierarchy of mathematical representations, from topological abstraction via PH, to geometric quantification via PL, to algebraic generalization via PCA, offering a unified foundation for interpretable, explainable, and generalizable AI.

\paragraph{Acknowledgments}
This work was supported in part by NIH grant R35GM148196, National Science Foundation
grant DMS2052983, Michigan State University Research Foundation, and Bristol-Myers Squibb
65109.

\bibliographystyle{abbrv}

\bibliography{reference}

\end{document}